\UseRawInputEncoding
\documentclass[12pt,thmsa,emstex]{article}
\usepackage{amsfonts}
\usepackage{t1enc}
\usepackage{amssymb}
\usepackage{epsfig}
\usepackage[french,english]{babel}
\usepackage{amsmath}
\usepackage{graphics}
\usepackage{layout}
\usepackage{latexsym}
\usepackage{color}
\usepackage{verbatim}

\def\FF{{\cal F}}

\def\mylabel#1{\label{#1}}

\newtheorem{theorem}{Theorem}[section]
\newtheorem{lemma}[theorem]{Lemma}
\newtheorem{corollary}[theorem]{Corollary}
\newtheorem{proposition}[theorem]{Proposition}

\newtheorem{definition}[theorem]{Definition}
\newtheorem{exercise}[theorem]{Exercise}

\newtheorem{remark}{Remark}
\newtheorem{example}{\bf{Example}}
\newtheorem{assumption}{\bf{Hypothesis}}

\def\bit{\begin{itemize}}
\def\eit{\end{itemize}}

\reversemarginpar   

\def\bc{\begin{center}}
\def\ec{\end{center}}

\def\bthm{\begin{theorem}}
\def\ethm{\end{theorem}}

\def\bcor{\begin{corollary}}
\def\ecor{\end{corollary}}

\def\bprop{\begin{proposition}}
\def\eprop{\end{proposition}}

\def\blem{\begin{lemma}}
\def\elem{\end{lemma}}

\def\bex{\begin{example}}
\def\eex{\end{example}}

\def\bexo{\begin{exercise}}
\def\eexo{\end{exercise} }

\def\brem{\begin{remark}}
\def\erem{\end{remark}}
\def\prf{{\bf Proof: }}

\def\bdes{\begin{description}}
\def\edes{\end{description}}

\def\ita{\item[(a)]}
\def\itb{\item[(b)]}

\def\itc{\item[(c)]}

\def\iti{\item[(i)]}
\def\itii{\item[(ii)]}

\def\itiii{\item[(iii)]}
\def\itiv{\item[(iv)]}

\def\itv{\item[(v)]}

\def\beq{\begin{equation}}
\def\eeq{\end{equation}}

\def\ben{\begin{enumerate}}
\def\een{\end{enumerate}}

\def\beqar{\begin{eqnarray}}
\def\eeqar{\end{eqnarray}}

\def\beqarr{\begin{eqnarray*}}
\def\eeqarr{\end{eqnarray*}}

\def\qed{\hfill $\Box$ \\[2ex]}
\def\prf{{\bf Proof: }\hspace{.1in}}

\newcommand{\LA}{\mathcal{L}}

\catcode`\@=11

\newcommand{\B}{\mathcal{B}}

\newcommand{\M}{\mathcal{M}}

\newcommand{\F}{\mathcal{F}}

\newcommand{\PR}{\mathcal{P}}
\newcommand{\DA}{\mathcal{D}}
\newcommand{\defn}{\stackrel{def}{=}}
\newcommand{\ua}{\uparrow}

\newcommand{\ep}{\epsilon}

\newcommand{\E}{\mathsf{E}}

\def\Ind{{\mathbf 1}}
\def\RR{{\mathbb R}}  
\def\Rp{{\mathbb R}_+}   

\def\NN{{\mathbb N}}
\def\QQ{{\mathbb Q}}

\def\Pr{{\mathsf P}}
\def\Gam{{\mathbf \Gamma}}
\def\D{\mathsf{diag}}
\def\rar{\rightarrow}
\def\eps{\varepsilon}

\begin{document}


\title{Stochastic Persistence\footnote{Revised version of the unpublished preprint "Stochastic Persistence, April 2014" } \\ (Part I)}

\author{ Michel Bena\"{i}m\thanks{Institut de Math\'{e}matiques, Universit\'{e} de Neuch\^{a}tel, Rue Emile-Argand, Neuch\^{a}tel, Suisse-2000. (michel.benaim@unine.ch).}}
\date{Updated version May 2019}
\maketitle
\begin{abstract}
Let $(X_t)_{t \geq 0}$ be a continuous time Markov process on some  metric space $M,$ leaving invariant a closed  subset
 $M_0 \subset M,$ called the {\em extinction set}.
We give general conditions ensuring either

{\em Stochastic persistence (Part I) :} Limit points of the occupation measure   are invariant probabilities over $M_+ = M \setminus M_0;$ or

{\em Extinction (Part II) :} $X_t \rar M_0$ a.s.

In the persistence case we also discuss conditions ensuring the a.s convergence (respectively exponential convergence in total variation) of the occupation measure (respectively the distribution) of $(X_t)$ toward a unique probability on $M_+.$

 These results extend and generalize previous results obtained for various stochastic models in population dynamics, given by stochastic differential equations, random  differential equations, or pure jump processes.
\end{abstract}

\pagestyle{myheadings}
\thispagestyle{plain}
\tableofcontents
\newpage
\section{Introduction} \mylabel{sec:intro}
 An important issue in  mathematical ecology and population biology is to find out  under which conditions a collection of interacting species can coexist over long periods of time. A similar question, in mathematical models of disease dynamics, is to understand whether or not  a disease will be endemic (i.e~persist in the population) or go extinct.
The mathematical investigation of these types of questions began with the early work of Freedman and Waltman \cite{FW77}, Gard \cite{Gard80b, Gard80a},  Gard and Hallam \cite{GHallam}, Schuster Sigmund and Wolff \cite{SSW79}, among others,
 in the late 1970s,
laying the foundation of what is now called the (deterministic) {\em mathematical theory of persistence}.
The theory  developed rapidly the past 35 years using the  available tools from dynamical system theory.
 The recent books by Smith and Thieme \cite{ST11}; Zhao and  Borwein  \cite{zhao2017dynamical} provide a comprehensive introduction to the theory as well as numerous examples and references.

For (most of)  deterministic models, {\em persistence} amounts to say that there exists an attractor bounded away from the {\em extinction states} (i.e~the subset of the states space where the abundance of one or group of the species vanishes).
  When this attractor is global, meaning that its basin of attraction includes all non-extinction states, the system is called {\em uniformly persistent} or
 {\em permanent} \cite{SSW79,H81}.

Beside biotic interactions, environmental fluctuations  play a key role in population dynamics. In order to take into account these fluctuations and to understand how they may affect
persistence, one approach  is the study of uniform persistence for non-autonomous difference or differential equations \cite{T00, MSZ04, ST11}. Another is to consider {\em systems subjected to environmental random perturbations}.    Classical examples include {\em ecological stochastic differential equations} (see e.g~ the classical paper by Turreli \cite{turelli1977random} or  \cite{lande2003stochastic}).
 \beq
\label{eq:ecol2}
dx_i = x_i [F_i(x) dt  + \sum_{j = 1}^m\Sigma_{i}^j(x) dB_t^j], \, i = 1 \ldots n
\eeq
where $(B^1_t,\ldots,B^m_t)$ is a standard $m$-dimensional Brownian motion;
and  {\em ecological stochastic equations driven by a Markov chain}
\beq
\label{eq:ecol3}
\frac{d x_i}{dt} = x_i(t) F_i^{J(t)}(x(t)), \, i = 1 \ldots n
\eeq
where $J(t) \in \{1, \ldots, m\}$ is a continuous time Markov chain - or more generally, a continuous time Markov chain controlled by $(x(t))$ -
 taking values in a finite set  representing different possible environments.
Both (\ref{eq:ecol2}) and (\ref{eq:ecol3}) are Markov processes defined on  $M = \Rp^n$ (respectively  $\Rp^n \times \{1,\ldots, m\}$) and describe the evolution of $n$ interacting species characterized by their abundances  $x_1, \ldots, x_n.$
The extinction set is the boundary  $M_0  = \partial \Rp^n$ (respectively $\Rp^n \times \{1,\ldots,m\}.$

Generalizing upon these models we will consider here a continuous time Markov process $(X_t)$ living in some metric space $M$ and leaving invariant a closed subset $M_0 \subset M,$ called the {\em extinction set}. That is $$X_0 \in M_0 \Leftrightarrow X_t \in M_0  \mbox{ for all } t \geq 0.$$ Observe that, when $X_0 \in M^+ := M \setminus M_0,$ $(X_t)$ is never absorbed by $M_0$ and extinction can only occur {\em asymptotically}.
The long term behavior of the process is then completely different from  the behavior of a process that would be absorbed (or killed) in {\em finite time} (see e.g the beautiful survey by Villemonais and Méléard \cite{Mel12} for a discussion of such processes). While extinction occurs in finite time
 for  most "realistic" finite population models, this extinction may be proceeded by long-term term transients when habitat sizes are sufficiently large. Hence, under this assumption, one can ignore the effects of {\em demographic stochasticity} (i.e. finite population effects) and focus on models with only {\em environmental stochasticity} where extinction can only be asymptotic. The recent survey paper by Schreiber \cite{schreiber2017coexistence} discusses these distinctions. Since the early observation by Hutchinson \cite{Hutchinson61} that temporal fluctuations of the environment can favor coexistence of species despite very limited resources, the effect of environmental stochasticity has been widely explored in the ecology literature, especially through the influence of Chesson and his coauthors \cite{chesson-warner-81, chesson-ellner-89,chesson-94, Chesson2000}.
\\

For deterministic models given by ecological differential equations - that is equation (\ref{eq:ecol2}) with $\Sigma_i^j = 0$ or (\ref{eq:ecol3}) with $m = 1$ -  general  sufficient conditions ensuring permanence or extinction (and generalizing many of the existing results), were derived by  Hofbauer, Schreiber, and their co-authors in a series of papers \cite{S00, GH03, HS04}. They rely on the existence of  a suitable  {\em average Lyapunov function}, a powerful notion introduced by Hofbauer  \cite{H81} in the early 1980s.

The central idea of the present paper is to define a similar object for Markov processes. First attempts in this directions include
\cite{BHW} dealing with small random perturbations of deterministic systems (i.e (\ref{eq:ecol2}) with small $\Sigma_i^j$) and later \cite{SBA11} for more general systems on compact state spaces (see also \cite{BS09} and \cite{jmb-14} for discrete time models). The results in \cite{BHW, SBA11} have been recently generalized by Hening and Nguyen \cite{HN18b} allowing to treat (\ref{eq:ecol2}) in full generality provided the diffusion term is non-degenerate.

In rough terms, our key assumption will be that there exist real valued continuous functions $V$ and $LV$ defined on $M^+$ with $V \geq 0$ (and typically $V(x) \rar \infty$ as $x \rar M_0$) such that
\bdes
\ita The process $$ M_t = V(X_t) - V(X_0) - \int_0^t LV(X_s) ds, t \geq 0$$ is a martingale for all $X_0 \in M^+;$
\itb $LV$ extends continuously to a function $H$ defined on all $M.$
\edes
In the deterministic case where $X_t$ is solution to an ordinary differential equation, say $\dot{X} = G(X),$ then $LV = \langle G, \nabla V \rangle, M_t = 0,$ and we recover Hofbauer's notion of average Lyapunov function.

Associated to $(V,H)$ are the {\em $H$-exponents}
$$\Lambda^-(H) = - \sup \int H(x) \mu(dx), \Lambda^+(H) = - \inf \int H(x) \mu(dx)$$ where the supremum (respectively  infimum) is taken over the set of ergodic measures for $(X_t)$ supported by $M_0.$ The sign of these exponents determine the behavior of the process near the extinction set.   We will show   that (under certain technical assumptions):
\begin{itemize}
\item (Part I). If $\Lambda^-(H)$ is positive, then
\begin{itemize}
\item The process is {\em stochastically persistent}, meaning that every limit point $\Pi$ of its empirical measure $$\Pi_t = \frac{1}{t}  \int_0^t \delta_{X_s} ds$$ is almost surely an invariant measure on $M^+.$ That is $\Pi(M^+) = 1.$

\item Under further irreducibility condition, such an invariant measure is unique  and the law of $(X_t)$ converges to $\Pi$ possibly at an exponential rate.
\end{itemize}
\item (Part II). If $\Lambda^+(H)$ is negative, then $X_t \rar M_0$ at rate
$$\liminf_{t \rar \infty} \frac{V(X_t)}{t} \geq - \Lambda^+(H)$$
\end{itemize}

This paper is a fully revised and extended version of the unpublished notes \cite{Ben14}, accompanying the Bernoulli lecture given by the author at the Centre Interfacultaire Bernoulli in october 2014.  Some of the ideas contained in these notes, have been already used in a few papers (\cite{BL16, HN18b,BStr17,HenStri17}  devoted to the analysis of certain ecological models. The present version has greatly benefitted from these papers. In particular, the beautiful analysis of the ecological sde (\ref{eq:ecol2}) conducted by Hening and Nguyen \cite{HN18b} has helped  to formulate conditions to deal with the situation where the extinction set is noncompact. Joint work with Edouard Strickler \cite{BStr17} has helped to understand how the general results here can be applied to the situation where the extinction set is no longer the boundary of the state space but an equilibrium point (a situation which naturally occurs in epidemic model), which after a natural change of variables, becomes a sphere.   Discussions with Joseph Hofbauer and Sebastian Schreiber over the recent years have been particularly influential.

\paragraph{Outline} The organization of Part I is as follows. Section \ref{sec:notation} introduces the notation and the main assumptions, ensuring in particular tightness of empirical measures. Section \ref{sec:examples} describes some motivating examples. Section \ref{sec:persistence} contains the main results: the persistence theorem (Theorem \ref{th:persistence}), conditions ensuring uniqueness of a persistent measure, convergence to this measure (Proposition \ref{th:petite} and Theorem \ref{th:OreyTH}), and under additional assumptions, exponential convergence (Theorems \ref{th:expoconvcompact} and \ref{th:expoconvnoncompact}). Section \ref{sec:sde2} applies these results to ecological SDEs (equation (\ref{eq:ecol2})) including degenerate ones. As an illustration, Section \ref{sec:RMAsde} analyzes a Rosenzweig-MacArthur model where the prey variable (but not the predator variable) is subjected to some small Brownian perturbation. Section \ref{sec:pdmpii} considers random ODEs driven by a Markov Chain (equation (\ref{eq:ecol3})) and, as an illustration, fully analyzes  in Section \ref{sec:MayLeonard} a $3$-dimensional process obtained by random switching between two May and Leonard vector fields. This provides an example for which the extinction set is not simply the boundary of the state space, but here the union of this boundary and an invariant line. The stochastic persistence results   combined with known results on competitive systems (in particular the theory of carrying simplices) allow to give precise conditions ensuring the existence of a unique persistent measure, absolutely continuous with respect to Lebesgue, and to characterize its topological support as the cell bordered by the carrying simplices of the two vector fields.
Section \ref{sec:proofposiv} contains the proof of the persistence Theorem and Section \ref{sec:rate} the proof of the exponential convergence results. Section  \ref{sec:append} is an appendix gathering some folklore results and their proofs.
\section{Notation and  hypotheses}
\label{sec:notation}
Let $(M,d)$ be a locally compact Polish space (e.g $\RR^n, \RR^n_+$ with the usual distance metric), equipped with its Borel $\sigma$-algebra $\B(M).$ We denote by $(\M_b(M),||\cdot||)$ the Banach space of all real-valued bounded measurable functions on $M$ under the {\em sup-norm} metric $||\cdot||$
  and  $C_b(M)$ (respectively $C_0(M)$) the Banach (sub)space of real-valued bounded continuous functions on
  $M$ (respectively real valued continuous functions vanishing at infinity). For any set $A \subset M$ we let $\Ind_A$ denote the  indicator function of $A.$
A generic nonnegative constant is noted  $cst.$
    We let
  $\PR(M)$ denote  the space of probability measures on $\B(M)$ equipped with the the topology of weak convergence.
  For $\mu \in \PR(M)$ and $f \in \M_b(M)$,
  we write $\mu f = \int_M f(x) \mu(dx)$.
 Recall that a sequence $(\mu_n)_{n \geq 1} \subset \PR(M)$ is said to converge weakly to $\mu \in \PR(M)$, written
   $\mu_n \Rightarrow \mu$, if for all $f \in C_b(M),\ \mu_n f \rar \mu f$.

Throughout the paper, we assume given  a  probability space $(\Omega, {\cal F},  \Pr),$   a complete right continuous filtration $({\cal F}_t),$   and a
family of {\em cad-lag} Markov processes $\{(X^x_t)_{t \geq 0}, \, x \in M\}$ on $(\Omega, {\cal F}, ({\cal F}_t)_{t \geq 0}, \Pr).$ By this we mean that
\bdes
\iti For all $x \in M$ $X^x_t$ is a $M-$valued $\FF_t$ measurable random variable, $X_0^x = x$ $\Pr$ a.s, and $t \rar X_t^x$ is {\em cad-lag} (i.e~ right-continuous with left-hand limits);
\itii For each $f \in \M_b(M)$ the mapping
\beq
\label{defPt}
(t,x) \in \RR^+ \times M \rar P_t f(x) = \E(f(X_t^x))
\eeq
is measurable, and
\beq \label{smgrpmrkvlnk}
\E\left[f(X^x_{t+s})|\F_t\right] = (P_s f)(X^x_t), \ \Pr ~ a.s.
\eeq
\edes
Equation (\ref{defPt}) defines a semigroup $(P_t)_{t\geq 0}$ of contractions on $\M_b(M).$ That is   $P_t \circ P_s f = P_{t+s} f$ and $\|P_t f \| \leq \|f\|.$

We sometimes let $\mathbb{P}_x$ denote the law of $(X_t^x)$ on the Skorokhod space $D(\Rp, M).$ That is  $\mathbb{P}_x(\cdot) = \Pr( \omega  \in \Omega \: : (X_t^x(\omega))_{t  \geq 0} \in \cdot).$
\\
Our main assumption is the following:
\begin{assumption}[Standing assumption] \label{hyp:standing}
{\rm
There exists a closed set $M_0 \subset M$ called the {\em extinction set} of $(P_t)_{t \geq 0}$ which is invariant under $(P_t)_{t \geq 0}:$
\[ \forall t \geq 0 \: P_t \Ind_{M_0} = \Ind_{M_0}.\]}
\end{assumption}
We let $M_{+} = M \setminus M_0$ denote the {\em non extinction set}. Note that $M_{+}$ is open and invariant (i.e~ $P_t \Ind_{M_+} = \Ind_{M_+}$).
\medskip

In addition to  Hypothesis \ref{hyp:standing} we make certain regularity and tightness assumptions (Hypotheses  \ref{hyp:feller} and \ref{assump02} below) that will be needed throughout.
\begin{assumption}[$C_b(M)$-Feller continuity]
\label{hyp:feller}
 For each $f \in C_b(M)$ the mapping $(t,x) \in \Rp \times M  \rar P_t f(x)$ is continuous.
\end{assumption}
\brem {\rm For further reference we will call such a semigroup  a  $C_b(M)$-{\em Feller Markov semigroup}. This terminology
 is chosen to avoid confusion
 with the usual definition of Feller Markov semigroups (see e.g \cite{EK1} or
 \cite{Legall2}) which assumes that $P_t$
  maps $C_0(M)$ into itself and induces a  strongly continuous semigroup on $C_0(M).$ Note that every Feller semigroup is $C_b(M)$-Feller.
 When $M$ is compact, all  the examples considered here are Feller (in the usual sense). However,
ecological stochastic differential equations on non compact spaces  are usually not, as shown in  the next example.}  \erem
\bex[Logistic SDE]
\label{ex:logistic}
{\rm
Consider the {\em logistic stochastic differential equation} on $\RR^+$ 
$$dx = x((1-x) dt + \sigma dB_t).$$ Then, for all $t > 0,$
$$X_t^x = \frac{x e^{ (1- \frac{\sigma^2}{2}) t + \sigma B_t}}{1 +   x  \int_0^t e^{(1- \frac{\sigma^2}{2}) s + \sigma B_s} ds} \rar X^{\infty}_t : = \frac{ e^{ (1- \frac{\sigma^2}{2}) t + \sigma B_t} }{\int_0^t e^{(1- \frac{\sigma^2}{2}) s + \sigma B_s} ds}$$
 as $x \rar \infty.$ It easily follows that  the induced semigroup doesn't preserve $C_0(\RR^+)$ nor that it is strongly continuous on $C_b(\RR^+).$ However, it is a $C_b(\RR^+)$ Feller Markov semigroup. }
\eex
\brem
\label{rem:strongmark}
{\rm The cad-lag continuity of the paths and Hypothesis \ref{hyp:feller} make $(X_t^x)$ a strong Markov process (see e.g~  Theorem 6.17 in  \cite{Legall2} stated for Feller (in the usual sense) Markov processes but the proof only requires cad-lag continuity and $C_b(M)$ Feller continuity).}
\erem

 We let  $\LA$ denote the {\em generator} of $(P_t)_{t \geq 0}$ on $C_{b}(M)$ and   $\DA(\LA) \subset C_{b}(M)$ its {\em domain.}
 Here, following \cite{priola99} (see also \cite{dpr11})
 $\DA(\LA)$ is defined as the set of $f \in C_b(M)$ for which
 \bdes
 \iti  $\LA f(x): = \lim_{t \rar 0}  \frac{P_t f(x) - f(x)}{t}$ exists for all $x \in M;$
 \itii $\LA f \in C_b(M);$
 \itiii $\sup_{0 < t \leq 1} \frac{1}{t} \| P_t f - f \| < \infty.$
 \edes
It is easily seen (see e.g~ Proposition 3.3 in \cite{dpr11}) that
 for all $f \in \DA(\LA)$ and $t \geq 0,$ $P_t f \in \DA(\LA),$ and that  for all $x \in M, t \mapsto P_t f(x)$ is $C^1$ and satisfies
 \beq
 \label{eq:sg}
\frac{d}{dt} P_t f(x)= \LA (P_t f) (x) = P_t (\LA f)(x).
\eeq
\brem
\label{rem:feller1} {\rm
In case $(P_t)$ induces a strongly continuous semigroup on a Banach set $E \subset C_b(M),$ (for instance $C_b(M)$ or $C_0(M)$)
 the set $\{(f,g) \in E \times E \: f \in \DA(\LA) , g = \LA f\}$ equals the graph
of the infinitesimal generator (defined in the usual sense) of $(P_t)$ restricted to $E.$ }
\erem
\brem {\rm For all $f \in C_b(M)$ and $\eps > 0$ let $f_{\eps} = \frac{1}{\eps}\int_0^{\eps} P_s f ds.$ Then $f_{\eps} \in \DA(\LA), \LA(f_\eps) = \frac{1}{\eps}(P_{\eps}f - f)$ and $\lim_{\eps \rar 0} f_{\eps} = f$ (pointwise).}
\erem
Let $\cal{M}$ denote one of the set $M$ or $M^+.$ We define the {\em extended generator} of $(P_t)_{t  \geq 0}$ on $\cal{M}$ as the set of  (possibly unbounded) continuous  maps $(f,g) : {\cal M} \mapsto \RR^2$  such that for all $x \in {\cal M}$ the process $(M_t^{f}(x)_{t \geq 0})$ defined as
\beq \label{martdefn}
M_t^{f}(x) = f(X_t^x) - f(x) - \int_0^t g(X_s^x) ds, t \geq 0
\eeq is a $({\cal F}_t, \mathbb{P}_x)$ martingale. If furthermore, $$\lim_{t \rar \infty} \frac{M_t^{f}(x)}{t} = 0$$ $\mathbb{P}_x$ a.s for all $x \in {\cal M}$ we say that $(f,g)$ {\em satisfies the strong law of large numbers}.

By (\ref{eq:sg}) for every $f \in {\cal D}({\cal L})$ $(f,\LA f)$ lies in the  extended generator of $(P_t)$ on $M.$  The next proposition (Proposition \ref{extendeddomain}) is a convenient tool to ensure that a given pair $(f,g)$ is in the extended generator and satisfies the strong law. It is one of the key tools that will be used throughout.  We first recall the definition of the {\em carr\'e du champ}. Let $\DA^2(\LA)$ denote the set of $f \in C_b(M)$ such that both $f$ and $f^2$ lie in $\DA(\LA).$
 If $f \in \DA^2(\LA)$ we let
 \beq
 \label{carre}
 \Gamma(f) = \LA f^2 - 2 f \LA f
 \eeq denote the {\em carr\'e du champ} of $f.$ Note that $\Gamma(f) = \lim_{t \rar 0} \frac{1}{t} (P_t f^2 - (P_t f)^2)$ so that $\Gamma(f) \geq 0.$

\bprop
\label{extendeddomain}
\label{martconv}
Let $\cal{M}$ be one of the set $M$ or $M^+.$ Let $(f, g) : {\cal M} \mapsto \RR^2$ be a continuous function. Assume that
for every compact set $K \subset {\cal M}$ there exists $f_K \in \DA^2(\LA)$ such that
\bdes
\ita $f|_K = f_K$ and $\LA (f_K) |_K = g |_K,$
\itb  $\forall x \in {\cal M}\; \sup\{ \frac{1}{t} \int_0^t P_s(\Gamma (f_K))(x) ds : \; t \geq 1, K \subset {\cal M} \; \mbox{ compact }\}  < \infty.$
\itc ***FIX JUMPS***
 \edes
 Then $(f,g)$ lies in the extended generator of $(P_t)$ on $\cal{M}$ and satisfies the strong  law of large numbers. Furthermore, $(M_t^{f}(x))_t$ is a $L^2$ martingale and for each compact set $K \subset {\cal M}$
 $$\langle M^f(x) \rangle_{t \wedge \tau_K} = \int_0^{t \wedge \tau_K} \Gamma(f_K)(X^x_s) ds,$$ where $(\langle M^f(x) \rangle_t)$ stands for previsible quadratic variation of $(M_t^f(x))$ and $\tau_K = \inf \{ t \geq 0 \: X^x_t \in K^c \}.$
 \eprop
The proof of this proposition is given in appendix, Section \ref{sec:appendmartconv}

\subsection{Empirical, invariant and ergodic probabilities}
We denote the sequence of {\em empirical occupation measures} $(\Pi^x_t)_{t\geq 0}$ of the process $(X^x_t)_{t \geq 0}$ as
\beq \label{empoccmeas}
\Pi^x_t(B) = \frac{1}{t} \int_{0}^{t} \Ind_{\{X^x_s \in B\}} ds, \ \forall B \in \B(M).
\eeq
Hence, $\Pi^x_t(B)$ is the proportion of time spent by the process
 in $B$ up to time $t.$

A probability measure $\mu \in \PR(M)$ is called {\em stationary} or {\em invariant} if
$$\mu P_t = \mu$$ for all $t \geq 0,$
 or equivalently, $\mu (P_t f) = \mu f $ for all $f \in \M_b(M)$ and all $t \geq 0.$
  We denote the set of invariant probability measures of $(P_t)_{t \geq 0}$ by $\PR_{inv}(M)$.  We also let \[\PR_{inv}(M_0) = \{\mu \in \PR_{inv}(M): \: \mu(M_0) = 1 \},\]
and \[\PR_{inv}(M_+) = \{\mu \in \PR_{inv}(M): \: \mu(M_+) = 1 \}.\]
A set $B \in {\cal B}(M) $ is called {\em invariant} if $P_t \Ind_B = \Ind_B$ for all $t \geq 0.$
Invariant probability $\mu \in \PR_{inv}(M)$ is called {\em ergodic} if every invariant set has $\mu-$measure $0$ or $1.$
Equivalently, 
$\mu$ is ergodic if and only if it is extremal, meaning that it cannot be written as a nontrivial convex combination $\mu = \ep \mu_1 + (1-\ep) \mu_0$ with $0 <\ep < 1$ of two other distinct invariant measures $\mu_0, \mu_1 \in \PR_{inv}(M)$. \\

 Given a set $S \subset M$ (typically $M, M_+$ or $M_0$) we denote by $$\PR_{erg}(S) = \{\mu \in \PR_{inv}(M), \, \mu(S) = 1, \mu \mbox{ ergodic} \}$$ the set of ergodic probability measures on $S$.

In order to control the behavior of the process at infinity and to ensure the tightness of $(\Pi_t^x)_{t \geq 0}$ (when $M$ is noncompact) we shall assume the existence of a convenient {\em Lyapunov function}.

Recall that a continuous map $W : M \mapsto \RR$ is called {\em proper} provided $\{x \in M  : W(x) \leq R\}$ is compact for all $R > 0.$

\begin{assumption} \label{assump02}
\label{hyp:tightpi} {\rm
There exist  proper maps  $W, \tilde{W} : M \mapsto \RR_+,$ and a continuous function $LW : M \mapsto \RR$  enjoying the following properties:
\bdes
\iti $(W,LW)$ is in the extended generator of  $(P_t)$ on $M$ and satisfies the strong law of large numbers;
\itii $LW  \leq - \tilde{W} + C$ for some  $C \geq 0.$
\edes}
\end{assumption}

\brem {\rm If $M$ is compact, Hypothesis \ref{hyp:tightpi} is automatically satisfied, say with $W = LW = \tilde{W} = 0.$} \erem
The next result ensures that, under Hypotheses \ref{hyp:feller} and \ref{hyp:tightpi}, the empirical occupation measure $(\Pi_t^x)$ is almost surely relatively compact and that its limit points are invariant. The proof is given in the appendix Section \ref{sec:appendtight}. Note that some versions of this results (for stochastic differential equations)  are already proved in \cite{SBA11} and \cite{EvansHeningSchreiber15}.
\bthm \label{tightlypnv}  Assumes Hypotheses \ref{hyp:feller} and \ref{hyp:tightpi}. Then
\bdes
\iti For all $x \in M$
$$0 \leq P_t W(x) + \int_0^t P_s(\tilde{W})(x) ds \leq W(x) + Ct.$$
\itii For all $x \in M,$ $\Pr$ almost surely, $\limsup_{t \rar \infty} \Pi_t^x \tilde{W} \leq C,$ $(\Pi_t^x)$ is  tight, and every limit point of  $(\Pi_t^x)_{t \geq 0}$ lies  in ${\cal P}_{inv}(M).$
Furthermore  ${\cal P}_{inv}(M)$ is   compact and $\mu \tilde{W} \leq C$ for all $\mu \in {\cal P}_{inv}(M).$
\itiii In case $\tilde{W} = \alpha W$ for some $\alpha > 0,$  $$P_t W \leq e^{-\alpha t}( W - C/\alpha) + C/\alpha.$$
\edes
\ethm
\brem {\rm Note that, while (by Theorem \ref{tightlypnv}) both ${\cal P}_{inv}(M)$ and ${\cal P}_{inv}(M_0)$ are non-empty, $\PR_{inv}(M+)$ may be empty.} \erem

\section{Motivating Examples}
\label{sec:examples}
\subsection{Pure jump ecological processes}
The simplest examples are given by pure jump processes.

Let $M = \Rp^n = \{ x \in \RR^n \: : x_i \geq 0\}, (E, {\cal E}, \nu)$ a probability space (representing the {\em environment}) and for each $i = 1, \ldots, n,$  $R_i : M \times E \mapsto \RR^*_+$
 a positive measurable mapping, continuous in the first variable.

  Vector  $x = (x_1, \ldots, x_n) \in M$ represents the state (abundances) of $n$ interacting species and  $R_i(x,e)$  the {\em fitness} of population $i$ in environment $e.$

  Let $(e_k)_{k \geq 1}$ be a sequence of i.i.d random variables distributed according to $\nu,$  and $(Y_k)_{k \geq 1}$ a discrete time Markov chain defined by
  $$Y_{k+1} = G(Y_k, e_{k+1})$$ where $$G(x,e) = (x_1 R_1(x,e), \ldots, x_n R_n(x,e)).$$
  Such  discrete time models of interacting populations  in a fluctuating environment are  analyzed in \cite{SBA11}.

 Let now  $(N_t)$ be a Poisson process with parameter $1$ independent of  $(e_k).$
  The process $$X_t = Y_{N_t}$$ is a jump Markov process on $M.$ The associated semigroup  is strongly continuous  on $C_b(M)$ (as well as on ${\cal M}_b(M)$) and writes $P_t f = e^{t \LA} f$ where $\LA$ is the bounded operator on $C_b(M)$ defined by
 $$\LA f (x) =  \int (f (G(x,e)) -  f(x)) \nu(de).$$
 Here $\DA(\LA) = \DA^2(\LA) = C_b(M)$ and $$\Gamma(f)(x) =  \int [f (G(x,e))- f(x)]^2 \nu(de).$$
For any given subset  $I \subset \{1, \ldots, n \},$ let
\beq
\label{defM0jump}
M_0^I = \{x \in M \:  : \prod_{i \in I} x_i = 0\}
\eeq be the set corresponding to the extinction of at least one of the species $i \in I.$
Hypothesis \ref{hyp:standing} is clearly satisfied with $M_0 = M_0^I.$ Hypothesis \ref{hyp:feller} is  satisfied by strong continuity of $(P_t).$
A  sufficient condition ensuring Hypothesis \ref{hyp:tightpi} is given by the existence of  suitable continuous Lyapunov function
 $V : M \mapsto \Rp$ for the discrete chain $(Y_n).$ For $f$ measurable and nonnegative,
 set$Kf(x) = \int (f (G(x,e)) \nu(de)$ and $L f(x) = Kf(x) - f(x).$

\bprop Assume there exists a continuous and proper map $V : M \mapsto \Rp$ such that
$KV \leq \rho V + C$ for some $0 \leq \rho < 1.$ Then Hypothesis  \ref{hyp:tightpi} is satisfied with $W = \sqrt{V}$ and $\tilde{W} = \sqrt{\rho V}.$
\eprop
\subsection{Ecological SDEs}
\label{ex:sde}
Consider  a stochastic differential equation having  the form
\begin{equation}
\label{eq:sde}
dx_i = x_i^{\alpha_i} [F_i(x)dt +   \sum_{j = 1}^m\Sigma_i^j(x)dB^j_t], i = 1, \ldots, n.
\end{equation}
where
 $F_i, \Sigma_i^j$ are real valued localy Lipschitz maps on $\RR^n,$ $\Sigma_i^j$ is bounded\footnote{This assumption is chosen here for simplicity and can be relaxed under other conditions as shown in \cite{HN18b}.}
$(B^1_t,\dots,B^m_t)$ is
an $m$-dimensional standard Brownian motion, and
$\alpha_i \in \NN$ (the set of nonnegative integers). The state space of (\ref{eq:sde}) is the set
$$M  = \{x \in \RR^n \: : \alpha_i > 0 \Rightarrow  x_i \geq 0 \}.$$

A variable $x_i$ for which $\alpha_i \neq 0$  typically represent the abundance of a  certain species, while a variable $x_i$ for which $\alpha_i \neq 0$ represents  a "feedback" or "abiotic" variable.
 The Brownian term    $(B^1_t,\dots,B^m_t)$ models the environmental noise.

This type of process includes Brownian perturbations of Lotka-Volterra processes
  as considered in \cite{EvansHeningSchreiber15}, \cite{HN18a,HN18c,HN18d} as well as general stochastic ecological equation that have been
recently considered by Hening and Nguyen in \cite{HN18b}.  The recent paper \cite{BenSch19} fully analyzes discrete time models having both internal (biotic) and external (abiotic) variables.

Let $I \subset \{1, \ldots, n \}$  be a subset of species. That is $\alpha_i > 0$ for all $i \in I.$ Let
\beq
\label{defM0sde}
M_0^I = \{x \in M \:  : \prod_{i \in I} x_i = 0\}
\eeq denote the {\em extinction set} corresponding to the extinction of at least one of the species $i \in I.$

We let $a(x)$ denote the  positive semi definite matrix defined by
\beq
\label{eq:defa}
 a_{ij}(x) = \sum_{k = 1}^m \Sigma_i^k(x) \Sigma_j^k(x).
 \eeq
For all $f : M \mapsto \RR, C^2,$ we let
\begin{equation}
\label{eq:defL1sde2} L f (x) = \sum_i  x_i^{\alpha_i} F_i(x) \frac{\partial f}{\partial x_i}(x) +
\frac{1}{2}\sum_{i,j} x^{\alpha_i}_i x^{\alpha_j}_j a_{ij}(x) \frac{\partial^2 f}{\partial x_i x_j}(x) \eeq
and
\beq
\label{eq:defGammasde} \Gamma_L(f) (x) = \sum_{i,j} x^{\alpha_i}_i x^{\alpha_j}_j a_{ij}(x) \frac{\partial f}{\partial x_i}(x) \frac{\partial f}{\partial x_j}(x).\eeq
The next proposition gives conditions ensuring that hypotheses  \ref{hyp:standing}, \ref{hyp:feller}, \ref{hyp:tightpi} hold. Its proof uses standard arguments given, for completeness, in appendix  Section \ref{sec:append}.

Recall that the maps $F$ and $\Sigma^j$ are  locally Lipschitz with $\Sigma^j$  bounded.
\bprop
\label{prop:ecosde}
 Assume that there exist a $C^2$  proper\footnote{i.e~$\lim_{\|x\| \rar \infty} U(x) = \infty$} map  $U : M  \mapsto [1, \infty[,$   a continuous function $\varphi : M \mapsto \Rp,$ and constants $\alpha > 0, \beta \geq 0$ and $0 \leq \eta < 1$ such that
\beq
\label{lyapsde}
L U \leq - \alpha U (1+\varphi) + \beta,
\eeq and
\beq
\label{gammasde}
\Gamma_L(U) \leq cst ( U^{2 + \eta}). \eeq

 Then
 \bdes
 \iti For each $x \in M$ there exists a unique (strong) solution $(X_t^x)_{t \geq 0} \subset M$ to
 (\ref{eq:sde})  with initial condition $X_0^x = x$ and $X_t^x$ is continuous in $(t,x);$ In particular Hypothesis \ref{hyp:feller} holds.
 \itii $\sup_{t \geq 0} \E(U(X_t^x)) \leq cst (1 + U(x)).$
 \itiii Let  $C^2_c(M)$ be the set of $C^2$ maps $f : M \mapsto \RR$ with compact support\footnote{By this we mean that
  is the restriction to $M$ of a $C^2$ function  $f : \RR^n \mapsto \RR$ with compact support.}. Then $C^2_c(M) \subset {\cal D}^2(\LA)$
 and for all $f \in C^2_c(M)$
$$
\LA f(x) = Lf(x)
\mbox{ and }
\Gamma(f)(x) = \Gamma_L(f)(x).$$
\itiv Hypothesis \ref{hyp:standing} holds true with $M_0 := M_0^I$
\itv  Hypothesis \ref{hyp:tightpi} holds with $$W = U ^{\frac{1-\eta}{2}}$$ and $\tilde{W} = (1 + cst)  W (1 + \varphi).$
\edes
\eprop

\brem[The Hening Nguyen condition]
{\rm
Set $\tilde{U} = \log(U).$ Then
$$LU = e^{\tilde{U}} ( L \tilde{U} + \frac{1}{2} \Gamma_L(\tilde{U})) \mbox{ and } \Gamma_L(\tilde{U}) = \frac{1}{U^2} \Gamma_L(U)$$ so that the above conditions on $U$ are  equivalent to the conditions

\beq
\label{eq:HN0}
\limsup_{\|x\| \rar \infty} L \tilde{U} + \frac{1}{2} \Gamma_L(\tilde{U}) + \alpha (1 + \varphi) < 0
\eeq
 for some $\alpha > 0$ and
\beq
\label{eq:HN0gamma} \Gamma_L(\tilde{U}) \leq cst(\exp{\eta \tilde{U}})
 \eeq for $\eta \geq 0.$

In particular, if $\hat{U} \geq 0$ is any $C^2$ proper function such that
\beq
\label{eq:HN}
      \limsup_{\|x\| \rar \infty} L(\hat{U}) +  \alpha (1 + \varphi)  <  0
      \eeq
      and
      \beq
      \Gamma_L(\hat{U})  \leq   cst
\eeq
Then the conditions ((\ref{eq:HN0}), (\ref{eq:HN0gamma})) are satisfied for  $\tilde{U} = \theta \hat{U}$ (i.e $U = e^{\theta \hat{U}}$), $\theta$ small enough and $\alpha$ replaced by $\alpha \theta.$

In case $\hat{U}(x) = \log (1 + \sum_i c_i x_i)$ with $c_i > 0,$  this condition is the one assumed in \cite{HN18b}.
}
\erem

\bex[Competitive Lotka-Volterra systems]
\label{ex:LVex}
{\rm
Consider the general model given by (\ref{eq:ecol2}) under the assumptions that $\alpha_i = 1$ for all $i = 1,\ldots, n$ and $$F_i(x) \leq f_i(x_i)$$ where $f_i : \RR \mapsto \RR$ is continuous and
\beq
\label{eq:Lotkacond}
x_i > R \Rightarrow p (f_i(x_i) + \frac{(p-1)}{2} a_{ii}(x)) < - \alpha
\eeq
for some positive numbers $R, \alpha$ and $p \geq 1.$
Then the conditions  of Proposition \ref{prop:ecosde} are satisfied with
$U(x) = 1 + \sum_i x_i^p, \varphi = 0$ and $\eta = 0$ (the verification is easy and left to the reader).
A particular case is given by the class of {\em competitive Lotka-Volterra} systems for which
\beq
\label{eq:compLV}
F_i(x) = r_i - \sum_{j} b_{ij} x_j
\eeq with $b_{ij} \geq 0$  and $b_{ii} > 0.$ Here, it suffices to chose $f_i(x) = r_i -b_{ii}x_i.$ Other examples include Lotka-Volterra mutualism systems as considered in \cite{Guo-Hu}.
}
\eex
\brem [Ecological SDEs on the simplex]
{\rm In numerous models occurring in ecology, population dynamics and game theory, $x_i$ represents the  proportion of species $i$ rather that its abundance. The state space is then the unit simplex
$$\Delta^{n-1} = \{x \in \Rp^n \: : \sum_i x_i = 1\}.$$  In this case,  to insure invariance of $\Delta^{n-1}$ by (\ref{eq:sde}),  one assumes that the drift and diffusion vector fields are tangent to $\Delta^{n-1}.$ That is
$$\sum_{i = 1}^n x_i^{\alpha_i} F_i(x) = \sum_{i = 1}^n  x_i^{\alpha_i} \Sigma_i^j(x) = 0.$$ Under these conditions, the processes (\ref{eq:sde}) induces a Feller (in the usual sense) Markov process on a compact metric space, $M = \Delta^{n-1}$. In particular, Hypotheses \ref{hyp:feller} and \ref{hyp:tightpi} hold, while Hypothesis  \ref{hyp:standing} obviously holds with $M_0$ defined by (\ref{defM0sde}).

Such model have been considered by Foster and Young \cite{Foster-Young}, Fudenberg and Harris \cite{Fudenberg-Harris},  Hofbauer and Imhof \cite{Hofbauer-Imhof}. A first general analysis of their persistence was first given by Benaim {\em et al.} \cite{BHW} and generalized in  Schreiber {\em et al.} \cite{SBA11}.
}
\erem

\subsection{Random ecological ODEs}
\label{sec:PDMP}
Let $\{G^j\}_{j = 1, \ldots, m}$ be a family of $m$ vector fields on $\RR^n$
having the form
$$G_i^j(x) = x_i^{\alpha_i} F_i^j(x); \;  i = 1, \ldots n$$ where $\alpha_i \in \NN$ and $F_i^j$ is $C^1.$

Set $$\RR^n_{\alpha} = \{x \in \RR^n \: : \alpha_i > 0 \Rightarrow  x_i \geq 0 \}.$$ As in Section \ref{sec:sde2}, the variables $x_i$ for which $\alpha_i > $ can be viewed as species abundances while the other variables are feedback variables. We let $\Phi^j = \{\Phi^j_t\}$ denote the  local flow on $\RR^n_{\alpha}$ induced by the ordinary differential equation
$\dot{x} = G^j(x).$

We  assume here for simplicity that there exists a compact set $B \subset  \RR_{\alpha}$  positively invariant under each $\Phi^j.$ That is
$\Phi^j_t(B) \subset B$ for all $t \geq 0.$

Let $$M = B \times \{1,\ldots,m\}.$$
For each $(x,j) \in M,$ let $(X_t^{x,j} = (x(t),J(t)))_{t \geq 0}$ be the process
on $M$ starting from $(x,j)$ (i.e~$X_0^{x,j} = (x,j)$) defined by
\beq
\label{eq:PDMP}
\left \{ \begin{array}{l}
 \displaystyle \frac{dx(t)}{dt}   =  G^{J(t)}(x(t)),\\
 \\
  \Pr (J(t+s)  =  k | \FF_t, J(t) = j) = a_{j k}(x(t)) s + o(s)
  \end{array} \right.
\eeq
where
$\forall \; j,k \in \{1,\ldots m\}, \, a_{j k} : \RR^n_{\alpha} \mapsto \Rp$ is continuous  nonnegative, $a_{j j} = 0,$ and $\{a_{j k}(x)\}_{j, k}$
is irreducible for all $x.$

This type of process belongs to the larger class of {\em Piecewise deterministic Markov Processes}, a term coined by Davis \cite{Dav84}.
Their ergodic  properties have recently been the focus of much attention in the literature
 (\cite{bakhtin&hurt},  \cite{ecp},  \cite{BMZIHP}, \cite{CH13}, \cite{bakhtin&hurt&matt}, \cite{BCL16}, \cite{2017arXiv170801390B}, \cite{BenHurthStric}). We refer the reader to the
     recent overview by Malrieu \cite{Mal15}.


Let $C^1(M)$ be  the set of maps $f : M \mapsto \RR, (x,j) \mapsto f(x,j)$ which are  $C^1$  in the  $x$ variable.
 It follows from Proposition 2.1 in \cite{BMZIHP} that $(X_t^{x,j})_{t \geq 0}$ is  Feller,  $C^1(M) \subset \DA^2(\LA)$ and
 for all
$f \in C^1(M)$
$$\LA f (x,j) = \langle \nabla_x f(x,j), G^j(x) \rangle + \sum_{k = 1}^m a_{j k}(x) (f(x,k) - f(x,j))$$
and
$$\Gamma(f)(x,j) = \sum_{k = 1}^m a_{j k}(x) (f(x,k) - f(x,j))^2.$$

Let $I \subset \{1, \ldots, n\}$ be a set of species (i.e $i \in I \Rightarrow \alpha_i > 0$). Then  Hypothesis \ref{hyp:standing} holds true with
 $$M_0 = M_0^I = \{ (x,j) \in M \: : \prod_{i \in I} x_i = 0\}.$$
\section{Stochastic Persistence and $H$-Exponents}
\label{sec:persistence}
The following definition, inspired by the seminal work of Chesson \cite{C78, C82}, follows from  Schreiber \cite{Sch12}.
\begin{definition} \label{defstochpers}
{\rm
The family $\{(X^x_t)_{t \geq 0} : x \in M_+\}$ is called {\em stochastically persistent (with respect to $M_0$)}  if for all $\eps > 0$ there exists a compact set $K_{\eps} \subset M_+$ such that for all $x \in M_+$:
\beq \label{stochpersdefn}
\Pr(\liminf_{t \rar \infty} \Pi_t^x(K_{\eps}) \geq 1-\eps) = 1.
\eeq
}
\end{definition}
If $M_0$ is unambiguous we  simply say that $\{(X^x_t)_{t \geq 0} : x \in M_+\}$ is stochastically persistent.

In  models of population dynamics, the interpretation of stochastic persistence is that all the species, initially present, persist (stay away from the extinction set) over arbitrary long periods of time.
\brem
\label{rem:decompM}
 {\rm
Suppose that $M_0 = M_0^1 \cup M_0^2$ where $M_0^{1,2}$ are closed and invariant under $(P_t)_{t \geq 0}.$ If the process if stochastically persistent with respect to $M_0^1$ and $M_0^2$ then it is stochastically persistent with respect to $M_0.$  Note that, however, the converse is false, as shown by the following deterministic example}
\erem
\bex
\label{ex:RMA}
{\rm
Consider the Rosenzweig MacArthur \cite{RMA63} prey predator model
\beq
\label{eq:RMA}\left \{ \begin{array}{l}
  \frac{dx_1}{dt} = x_1( 1- \frac{x_1}{\kappa} - \frac{x_2}{1+x_1})\\
  \frac{dx_2}{dt} = x_2 (- \alpha + \frac{x_1}{1+x_1}) \end{array} \right.
 \eeq
 on the state space $M = \Rp^2$ where $\alpha, \kappa$ are  positive parameter. Set $M_0^1 = \Rp \times \{0\}, M_0^2 = \{0\} \times \Rp$ and $M_0 = M_0^1 \cup M_0^2.$
Every trajectory on $M_0^2$ converges to the origin, so that the system is never persistent with respect to  $M_0^1.$
Assume $\alpha < \frac{\kappa}{1+ \kappa}.$ Then (see e.g~\cite{SMITHonRMA}) the  system admits an  equilibrium $p  \in M_+.$ If
 $\alpha < \frac{\kappa-1}{\kappa + 1}$ $p$ is a source and there is
and a limit cycle $\gamma  \subset M_+$ surrounding $p$ whose basin is $M_+ \setminus \{p\}.$ If $\alpha \geq \frac{\kappa-1}{\kappa + 1}$ every positive trajectory converges to $p.$ This makes the system persistent with respect to $M_0.$}
\eex

Proving or disproving stochastic persistence requires to control the behavior of $(X_t)$ near the extinction set. This will be done by assuming the existence of another suitable type Lyapunov function.
\begin{assumption}
\label{hyp:H}
There exist  continuous maps $V : M_{+} \mapsto  \RR_{+}$ and  $H : M \mapsto \RR$  enjoying the following properties:
 \bdes
 \iti  The pair $(V,H|_{M_+})$ lies in the extended generator of $(P_t)$ on $M_+$ and satisfies the strong law;
\itii The map $\frac{\tilde{W}}{1 + |H|}$ is proper, where  $\tilde{W}$ is like in Hypothesis \ref{hyp:tightpi}.
\edes
\end{assumption}
We will sometimes assume the stronger version of {\bf (ii)}:
\bdes
\item[(ii)'] $|H|^q \leq cst (1 + \tilde{W})$ for some $q > 1.$
\edes
Note that by condition $(ii)$ above and Theorem \ref{tightlypnv}, $H \in L^1(\mu)$ for all $\mu \in \PR_{inv}(M).$ The following definition then makes sense.
\begin{definition}[H-exponents]{\rm
 If $V$ and $H$ are like in Hypothesis \ref{hyp:H} we let
 $$\Lambda^-(H) = - \sup \{\mu H \: : \mu \in \PR_{erg}(M_0) \},$$ and $$\Lambda^+(H) = - \inf \{\mu H \: : \mu \in \PR_{erg}(M_0) \}$$ denote the {\em $H$-exponents} of $(X_t).$
 }
 \end{definition}

 \brem
 \label{rem:whenlambdais0}
 {\rm The key point here is that while $H$ is defined on {\bf all} $M,$ $V$ is defined {\bf only} on $M_+$ and typically $V(x) \rar \infty$ when $x \rar M_0.$

 Actually, if $V$ can be  defined {\bf on all} $M$ with condition $(i)$ of Hypothesis \ref{hyp:H} valid on $M$ then (see Remark \ref{rem:muHzero})  $$\Lambda^-(H) = \Lambda^+(H) = 0.$$
 }
 \erem
\begin{definition} \label{defHpers}
{\rm
We call  $\{(X^x_t)_{t \geq 0} : x \in M_+\}$  $H$-{\em persistent} if there exists $(V,H)$ like in Hypothesis \ref{hyp:H} such that  $\Lambda^{-}(H) > 0.$
}
\end{definition}

\bex[Logistic SDE, continuation of example \ref{ex:logistic}]
 \label{ex:logistic2}
 {\rm Consider the logistic equation given in Example \ref{ex:logistic}. Here  $M = \Rp$ and  $M_0 = \{0\}.$ Let $V : ]0, \infty[ \mapsto \Rp$ be any smooth function with bounded support (say, $V(x) = 0$ for $x \geq 1$) and coinciding with $-\log(x)$ on a neighborhood of $0.$ Then the map $H(x) = V'(x) x(1-x) + \frac{\sigma^2}{2} x^2 V''(x)$  extends continuously to $\Rp$ and coincide with $x - 1 +  \frac{\sigma^2}{2}$ on a neighborhood of $0.$
Clearly $V$ and  $H$ satisfy Hypothesis \ref{hyp:H}, $\PR_{erg}(M_0) = \{0\}$ and $$\Lambda^+(H) = \Lambda^{-}(H) = 1- \frac{\sigma^2}{2}.$$
Here $H-$ persistence simply writes
$$1- \frac{\sigma^2}{2} > 0.$$ More sophisticated examples will be studied later.}
\eex
\brem
{\rm
By the ergodic decomposition theorem, compactness of $\PR_{inv}(M_0)$ (Theorem \ref{tightlypnv}) and continuity of $\mu \mapsto \mu H$ (Lemma \ref{lem:tightpi} (ii))
 combined with condition  $(ii)$ in Hypothesis  \ref{hyp:H}, the following conditions are equivalent :
 \ben
 \ita  $\Lambda^{-}(H) > 0,$
 \itb $\mu H < 0$ for all $\mu \in \PR_{erg}(M_0),$
 \itc  $\mu H < 0$ for all $\mu \in \PR_{inv}(M_0).$
 \een
  }
  \erem
Similar to Remark \ref{rem:decompM} is the following
\brem
{\rm
Suppose that $M_0 = M_0^1 \cup M_0^2$ where $M_0^{1,2}$ are closed and invariant under $(P_t)_{t \geq 0}.$ Let $M_+^i = M \setminus M_0^i,$ and
 $V^i : M_+^i \mapsto \Rp, \,H^i : M \mapsto \RR$
be as in Hypothesis \ref{hyp:H}. Let $V$ be defined on $M_+$ by
$V = V^1 + V^2$ and let $H = H^1 + H^2.$

Then for all $\mu \in \PR_{erg}(M_0)$ either
\bdes
\iti
$\mu(M_0^1) = 0$ and $\mu H = \mu H^2,$ or
\itii $\mu(M_0^2) = 0$ and $\mu H = \mu H^1,$  or
\itiii
 $\mu(M_0^1 \cap M_0^2) = 1$ and $\mu H = \mu H^1 + \mu H^2.$
 \edes In particular
if the process is $H$-persistent with respect to $M_0^i, i = 1, 2$ it is $H$-persistent with respect to $M_0.$}
\erem

\subsection{H-Persistence implies Stochastic Persistence}
From now on, hypotheses \ref{hyp:standing} to \ref{hyp:H} are implicitly assumed.
The main result of this section is given by the following theorem whose  proof is postponed to Section \ref{sec:proofposiv}.
\begin{theorem} \label{posinvthm} \label{th:persistence}
 Assume that $\{(X^x_t)_{t \geq 0} : x \in M_+\}$ is $H$-persistent. Then
 \bdes
 \iti For all $x \in M_+,$ every weak limit point of $(\Pi_t^x)_{t \geq 0}$ lies in $\PR_{inv}(M_+)$ a.s.
\itii $\{(X_t^x)_{t \geq 0} : x \in M_+\}$ is stochastically persistent.
\edes
\end{theorem}
This theorem has the following immediate consequence.
\bcor
\label{cor:posinv}
 Assume that $\{(X^x_t)_{t \geq 0} : x \in M_+\}$ is $H$-persistent and  that ${\cal P}_{inv}(M_+)$ has cardinal at most one.
  Then, ${\cal P}_{inv}(M_+)$  has cardinal one, and letting ${\cal P}_{inv}(M_+) = \{\Pi\}$, for all $x \in M_+,$
  $\Pi_t^x \Rightarrow \Pi$ a.s. as $t \rar \infty.$
  \ecor
  For further references,  the probability $\Pi$ in Corollary \ref{cor:posinv} is called the {\em persistent measure}. In ecological models, the persistence measure describes the long term behavior of coexisting species.
\subsection{Proof of Theorem \ref{posinvthm} }
\label{sec:proofposiv}
 Since $M$ is  locally compact and separable there exists a sequence $\{C_n\}_{n \geq 1}$ of compact sets with $C_n \subset \mathrm{int}(C_{n+1})$ such that $M = \cup_{n \geq 1} C_n$. Throughout we let
\beq
\label{defKn}
K_n = \{x \in M : d(x,M_0) \geq \frac{1}{n}\} \cap C_n.
 \eeq
 Note that $$M_+  = \bigcup_{n \geq 1} K_n$$
and $$M_0 = \bigcap_{n\geq 1} {K_n^c} = \bigcap_{n\geq 1} \overline{K_n^c}$$
where the later equality follows from the inclusion  $C_n \subset \mathrm{int}(C_{n+1}).$

The proof of the following Lemma is similar to the proof of Proposition 1 in \cite{SBA11}.

\blem \label{muHlemma}
Assume that $\{(X_t^x)_{t \geq 0} : x \in M_+\}$ is $H$-persistent. Then
\bit
\iti For all $\mu \in \PR_{inv}(M),\ \mu H \leq 0;$ and $\mu H = 0  \Leftrightarrow\mu \in \PR_{inv}(M_+)$.
\itii  $\PR_{inv}(M_+)$ is tight : $\forall \eps > 0, \, \exists K \subset M_+$ compact such that $\inf \{\mu(K)  \: :\mu \in \PR_{inv}(M_+)\} \geq 1-\eps .$
\eit
\elem
\prf
 $(i).$ By Hypothesis \ref{hyp:H} $(ii)$ and Theorem \ref{tightlypnv}  $H \in L^1(\mu)$ for all $\mu \in  \PR_{inv}(M).$ Let $\mu \in \PR_{inv}(M_{+}).$ We claim that $\mu H = 0.$ By the ergodic decomposition theorem it suffices to prove the result for $\mu$  ergodic. By Birkhoff ergodic Theorem, for $\mu$ almost all $x$ and $\mathbb{P}_x$ almost surely $$\lim_{t \rar \infty} \Pi_t^x H = \mu H.$$ Hence, by Hypothesis \ref{hyp:H},
 $$\lim_{t \rar \infty} \frac{V(X_t^x)}{t} = \mu H.$$  Since $\mu(M_0) = 0$ there exists $n \geq 1$ such that  $\mu(K_n) \geq 1/2$, so that, by  Birkhoff ergodic Theorem again, $(X_t^x)_{t  \geq 0}$ visits $K_n$ infinitely often for $\mu$ almost all $x,$ $\mathbb{P}_x$ almost surely. Since $V$ is bounded on $K_n$ this proves that $\mu H = 0.$

Let now $\mu \in \PR_{inv}(M) \setminus \PR_{inv}(M_{+}).$
We can write, by Hypothesis \ref{hyp:standing}, $\mu = (1-t)\mu_0 + t \mu_1,\ 0 \leq t < 1$, with $\mu_0 \in \PR_{inv}(M_0)$ and $\mu_1 \in \PR_{inv}(M_+)$. Thus
$\mu H = (1-t) \mu_0 H < 0$

 $(ii).$ Suppose not. Then there exists some $\ep > 0$ such that for each $n \geq 1$ there exists some $\mu_n \in \PR_{inv}(M_+)$ with $\mu_n(K_n) < 1- \ep.$
Thus,  $\mu_n(\overline{K_m^c}) > \ep$ for all $m < n$ as, by definition, $\overline{K_{n+1}^c} \subset \overline{K_n^c}$. Let $\mu$ be a limit point of $(\mu_n)$ for the weak* topology. Then $\mu \in \PR_{inv}(M)$ as $\PR_{inv}(M)$ is tight and,
 by application of Portemanteau
 $\mu(\overline{K_m^c}) \geq \ep$ for all $m \geq 1.$
Since  $M_0 = \cap_{m \geq 1}\overline{K_m^c}$ this implies $\mu(M_0) \geq \ep.$
Now, by  by part (i) $\mu_n H = 0$ implying   $\mu H = 0$
and $\mu \in \PR_{inv}(M_+)$ again by part (i). A contradiction.
\qed
\brem
\label{rem:muHzero}
{\rm The proof of Lemma \ref{muHlemma} $(i)$ also shows that when $M_0 = \emptyset$, then $\mu H = 0$ for all $\mu \in {\cal P}_{inv}(M).$}
\erem
We now prove Theorem \ref{posinvthm}.

 $(i).$ By Theorem \ref{tightlypnv}, for every $x \in M_+$, every weak limit point of $(\Pi_t^x)_{t \geq 0}$ lies $\mathbb{P}_x$-a.s.\ in $\PR_{inv}(M).$ Let $\mu = \lim_{n \rar \infty} \Pi_{t_n}^x \in \PR_{inv}(M)$ be such a weak limit point. Then,  by  Proposition \ref{martconv} again, $\lim_{n \rar \infty} \frac{V(X_{t_n}^x)}{t_n} = \mu H \geq 0$ (because $V \geq 0$) and the result follows from assertion $(i)$ of  Lemma \ref{muHlemma}.

$(ii).$
Suppose that this is not true implying that there exists some $\eps > 0$ and a sequence $(x_n)  \subset M_+$
such that
$$\Pr(\liminf_{t \rar \infty} \Pi_t^{x_n} (K_n^c) \geq \eps) > 0.$$
By assertion $(i)$, with probability 1 there exists a subsequence $\{t^{n}_m\}_{m\geq 1} \ua \infty$ and $\mu_n \in \PR_{inv}(M_+)$
such that $\liminf_{t \rar \infty} \Pi_t^{x_n} (K_n^c) = \lim_{m \rar \infty} \Pi_{t^{n}_m}^{x_n} (K_n^c)$ and $\Pi_{t^{n}_m}^{x_n} \Rightarrow \mu_n$ as $m \rar \infty.$ Thus, by Portemanteau theorem, $\mu_n(\overline{K_n^c}) \geq \eps$ on the event $\liminf_{t \rar \infty} \Pi_t^{x_n} (K_n^c) \geq \eps.$
By tightness of $\PR_{inv}(M_+),$ (Lemma \ref{muHlemma}, $(ii)$) for $n$ large enough $\mu_n(\overline{K_n^c}) < \eps.$ A contradiction.

\subsection{Support and Irreducibility}
\label{sec:support}
Unlike in deterministic models where persistence equates the existence of an attractor bounded away from the extinction set, the support of the persistent measure  may well have nonempty intersection with  $M_0.$

  The nature of this support  provides useful information on the dynamics. For specific models (see for instance \cite{BL16} section 4,
   and  \cite{2016arXiv160108151M} section 5) it can be computed by using  some elementary control theory type arguments that we now briefly discuss. The general definitions given here will be rephrased in Sections  \ref{sec:sde2} and \ref{sec:pdmpii} in terms of deterministic control systems.

Point $y \in M$ is said {\em accessible} from $x \in M$ if for every neighborhood $U$ of $y$ there exists $t \geq 0$ such that $P_t(x,U) = P_t \Ind_U (x) > 0.$ We let ${\mathbf \Gamma}_x$ denote the set of points $y$ that are accessible from $x.$ For  $A \subset M$ we let
$$\Gam_A = \cap_{x \in A} \Gam_x$$ denote
the (possibly empty) closed set of points that are accessible from every $x \in A.$
\bcor
\label{cor:supportofpi}
Under the assumptions of Corollary \ref{cor:posinv}, $$supp(\Pi) = \Gam_{M_+} = \Gam_x$$ for all $x \in \Gam_{M_+} \cap M_+.$
\ecor
\prf Let $O \subset M$ be an open set such that $\Pi(O) > 0.$ Then, by Fatou lemma, Corollary \ref{cor:posinv}   and Portmanteau theorem,  for all $x \in M_+$
$$\liminf_{t \rar \infty} \frac{1}{t} \int_0^t P_s \Ind_O(x) ds \geq \E(\liminf_{t \rar \infty} \Pi_t^x(O)) \geq \Pi(O) > 0.$$ This proves that $supp(\Pi) \subset \Gam_{M_+}.$

Conversely, let $R(x,dy)$ be the resolvent kernel defined by
$Rf = \int_0^{\infty} e^{-s} P_s f ds.$  We claim that for every $y \in \Gam_{M_+},$ $O$ a neighborhood of $y$ and $x \in M_+$ $R(x,O) > 0.$
By accessibility, there exists $t > 0$ such that $P_t(x,O) > 0.$ Thus, by right continuity and Fatou Lemma
$$\liminf_{s \rar t, s > t} P_s(x,O) \geq \E(\liminf_{s \rar t, s > t} \Ind_O(X_s^x)) \geq P_t(x,O) > 0.$$ This proves that $s \rar P_s(x,O)$ is positive on some interval $[t,t+\eps],$ hence $R(x,O) > 0.$

Now, by invariance, $\Pi = \Pi R.$  Therefore,
$\Pi(O) = \int_{M_+} \Pi(dx) R(x,O) > 0.$ This proves that $\Gam_{M_+} \subset supp(\Pi).$

We now prove the last assertion.
By definition $\Gam_{M_+} \subset \Gam_x$ for all $x \in M_+.$ It then suffices to show that $\Gam_x \subset \Gam_{M_+}$ for $x \in \Gam_{M_+} \cap M_+.$
Let $y \in \Gam_x$ and $O$ a neighborhood of $y.$ Then for some $t \geq 0$ and $\delta > 0$ $P_t(x,O) > \delta.$ By $C_b(M)$ Feller continuity and Portmanteau Theorem, the set $V = \{z \in M_+ \: : P_t(z,O) > \delta \}$ is an open neighborhood of $x.$ But since $x \in \Gam_{M_+}$ for all $z \in M_+$ there is some $s > 0$ such that $P_s(z,V) > 0$
Thus $P_{t+s}(z,O) \geq \delta P_s(z,V) > 0.$
\qed
\brem
\label{rem:supportofpi}  {\rm
The preceding proof also shows that  $\Gam_{M^+} \subset supp(\Pi)$ for all $\Pi \in {\cal P}_{inv}(M^+)$ even if ${\cal P}_{inv}(M^+)$ has cardinal greater than $1.$ However, in this case, $\Gam_{M^+}$ may be empty.}
\erem
A sufficient (although non-necessary) condition ensuring that ${\cal P}_{inv}(M_+)$ has cardinal at most one (hence one when the process is stochastically persistent) is given by $\psi$-irreducibility, in the sense of Meyn and Tweedy (see  \cite{MT1} or \cite{duf00}). A practical condition (implying $\psi$-irreducibility) is the existence of  an {\em accessible weak Doeblin point}.

We say that $x^* \in M$ is a {\em weak Doeblin point} if there exists a neighborhood $U$ of $x^*,$ a non zero  measure $\xi$ on $M$, and a probability measure $\gamma$ on $\Rp$ such that for all $x \in U$
$$R_{\gamma}(x,\cdot) := \int P_t(x,\cdot) \gamma(dt) \geq \xi(\cdot).$$
\bprop
\label{th:petite}
Assume there exist a weak Doeblin point $x^* \in \Gam_{M_+}$ (in particular $\Gam_{M_+} \neq \emptyset$). Then ${\cal P}_{inv}(M_+)$ has cardinal at most one. If furthermore, $(X_t)$ is $H$-persistent, then ${\cal P}_{inv}(M_+) = \{\Pi\}, \Pi_t^x \Rightarrow \Pi$ a.s.  and
$$\lim_{t \rar \infty} \Pi_t^x f = \Pi f$$ a.s. for all $f \in L^1(\Pi)$ and $x \in M^+.$
\eprop
\prf As shown in the proof of the last corollary,  accessibility implies that for all  $x \in M_+$ $R(x,U) > 0,$ where $R$ is the resolvent.
Thus $R R_{\gamma}(x,\cdot) \geq \int_U R(x,dy) R_{\gamma}(y,\cdot) \geq R(x,U) \xi(\cdot).$
This proves that $R R_{\gamma}$ is $\xi$ irreducible on $M_+$. Thus, $R R_{\gamma},$ and therefore  $P_t,$
has at most one invariant probability on $M_+$ (see e.g~\cite{MT1} or \cite{duf00}). If furthermore $(X_t)$ is  $H$ persistent, then $\Pi_t^x \Rightarrow \Pi$ a.s by Corollary \ref{cor:posinv}. It remains to prove the last assertion. The set  $U$ is, by assumption, a {\em petite set} in the sense of Meyn and Tweedy \cite{MT1}. By Portmanteau theorem, accessibility and Corollary \ref{cor:supportofpi}, $\liminf \Pi_t^x(U) \geq \Pi(U) > 0$ proving that $U$ is recurrent. Now, the existence of a petite and recurrent set  makes $(X_t)$ {\em Harris recurrent} on $M^+$, and since ${\cal P}_{inv}(M_+)$ is nonempty $(X_t)$ is positively recurrent on $M^+.$ \qed
\brem
\label{rem:petite} {\rm In many cases there exists a measure $\lambda$ on $M$  such that  $\lambda P_t \ll \lambda.$ A typical situation is when $M \subset \RR^n,$ $\lambda$ the is the Lebesgue measure on $\RR^n$ and $x \mapsto X_t^x(\omega)$ is a $C^1$ diffeomorphism for all (or $\Pr$ almost all) $\omega$. If in addition,  the assumptions of Proposition \ref{th:petite} are satisfied with
$\xi \ll \lambda,$ then $$\Pi \ll \lambda.$$
Indeed, by Lebesgue decomposition Theorem $\Pi = \Pi_{ac} + \Pi_s$ with $\Pi_{ac} \ll \lambda$ and $\Pi_s \perp \lambda.$ The proof of Proposition \ref{th:petite} easily implies that  $\xi \ll \Pi.$ Thus $\Pi_{ac} \neq 0$ because   $\xi \ll \lambda.$ By invariance $\Pi_{ac} + \Pi_s = \Pi_{ac} P_t + \Pi_s P_t.$ Thus $\Pi_{ac} \geq \Pi_{ac} P_t$ by uniqueness of Lebesgue decomposition. This shows that  $\frac{\Pi_{ac}}{\Pi_{ac}(M)}$ is excessive, hence invariant. That is $\Pi = \Pi_{ac}.$}
\erem
 \subsection{Convergence}
 The next result shows that if the measure $\gamma$ in Proposition \ref{th:petite} can be chosen to be a dirac mass then the law of $X_t^x$ converges in total variation  to $\Pi$ whenever $x \in M_+.$

Recall  that the {\em total variation distance} between two probabilities $\alpha, \beta \in {\cal P}(M)$ is defined as
$$|\alpha - \beta|_{TV} = \sup \{ |\alpha(f) - \beta(f)| ~:~  f \in {\cal M}_b(M),~ \|f\|_{\infty} \leq 1\}.$$
We say that $x^* \in M$  is a {\em Doeblin point} if there exist a neighborhood $U$ of $x^*,$ a non zero  measure $\xi$ on $M$,  and $t^* > 0$  such that for all $x \in U$
\beq
\label{eq:Doeblin}
P_{t^*}(x,\cdot)  \geq \xi(\cdot).
\eeq
If a Doeblin point is accessible, the minorization condition (\ref{eq:Doeblin}) extends to every compact space. More precisely
\blem
\label{lem:compactsmall}
Let $x^* \in \Gam_{M_+}$ be a Doeblin point and $x_0 \in supp(\xi) \cap M_+$ where $\xi$ is like in (\ref{eq:Doeblin}). Then there exist a neighborhood $A \subset M_+$ of $x_0,$
a probability $\nu$ on $A$ (i.e $\nu(A) = 1$) and positive numbers $T, c$ such that:
\bdes
\iti For all $x \in A$ $P_T(x, \cdot) \geq c \nu(\cdot);$
\itii For every compact set $K \subset M_+$ there exist $n_K \in \NN, c_K > 0$ such that $P_{n_k T}(x,\cdot) \geq c_K \nu(\cdot)$ for all $x \in K;$
\edes
\elem
\prf  Let  $t^*$ and $U$ be like in (\ref{eq:Doeblin}). By accessibility there exist  $t_0, \delta > 0$  such that $P_{t_0}(x_0,U) > \delta.$  By $C_b(M)$-Feller continuity (Hypothesis \ref{hyp:feller}) and Portmanteau's theorem there exist $\eps > 0$ and an open neighborhood $A$ of $x_0$ such that $P_t(x,U) > \delta$ for all $x \in A$ and  $|t-t_0| < \eps.$  Set $T =  t_0 + t^*, \nu(\cdot) = \frac{\xi(\cdot \cap A)}{\xi(A)}$ and $c = \delta \xi(A).$   Then, $P_{\tau}(x,\cdot) \geq  c \nu(\cdot)$ for  $|\tau -T| < \eps.$    This proves $(i).$

Let $O(t,r) = \{ x \in M_+ \: : P_{t}(x,U) > r\}.$ The family $\{O(t, r)\}_{t \geq 0, r > 0}$  is an open (by $C_b(M)$-Feller continuity) covering (by accessibility)  of $M_+.$ Thus, for $K \subset M_+$ compact, $K \subset \cup_{i = 1}^l O(t_i,r)$ for some $l \in \NN, t_1, \ldots, t_l \geq 0$ and $r > 0.$
 Choose $k$ large enough so that $\frac{t_i + t^*}{k} < \eps$ for all $i = 1, \ldots, l$ and set $\tau_i = T-\frac{t_i + t^*}{k}.$
Then for $x \in O(t_i,r)$ $$P_{kT}(x, \cdot ) = P_{t_i + t^* + k \tau_i}(x,A) \geq r \xi(A) (c \nu(A))^{k-1}\nu(\cdot).$$
This proves $(ii)$.
\qed
\bthm
  \label{th:OreyTH} Assume that $\{(X^x_t)_{t \geq 0} : x \in M_+\}$ is $H$-persistent and that there exists a Doeblin point  $x^* \in \Gam_{M_+}.$  Then  $${\cal P}_{inv}(M_+) = \{\Pi\}$$ and
for all $x \in M_+$
$$\lim_{t \rar \infty} |P_t(x,\cdot) - \Pi|_{TV} = 0.$$
\ethm
\prf We use the notation of Lemma  \ref{lem:compactsmall}. Let $Y$ be   the discrete chain on $M_+$ whose transition kernel is $P_T$ (restricted to  $M_+$). By  Lemma  \ref{lem:compactsmall} $(i)$,
$A$ is a {\em small} set for $Y$   and, by $(ii),$ it is accessible for $Y$ from every point in $M_+.$  In addition, by Proposition \ref{th:petite}, $Y$ has an invariant probability implying that $A$ is {\em recurrent}. By application of Orey's theorem (see e.g Theorem 8.3.18 in  \cite{duf00}), the existence of a {\em small accessible recurrent set}  imply that
$$\lim_{n \rar \infty} |\mu P_T^n - \Pi|_{TV} = 0$$
for all $\mu \in {\cal P}(M_+).$ Now, writing $t = n_tT + r_t$ with $0 \leq r_t < T, n_t \in \NN,$
$$\lim_{t \rar \infty}  |\mu P_t - \Pi|_{TV} = \lim_{t \rar \infty} |\mu P_{n_t T}  P_{r_t} - \Pi P_{r_t} |_{TV} \leq \lim_{n \rar \infty} |\mu P_{nT}  - \Pi |_{TV} = 0.$$
\qed
\subsection{Rate of Convergence}
Under certain additional assumptions, the rate of convergence in Theorem \ref{th:OreyTH} can be shown to be exponential.

Throughout this section we will assume the following strengthening of Hypothesis \ref{hyp:H}:
\begin{assumption}[strong version of Hypothesis \ref{hyp:H}]
\label{hyp:H+}
$V$ and $H$ are like in Hypothesis \ref{hyp:H}, and in addition:
\bdes
\ita The jumps of $V(X_s)$ are almost surely bounded. That is $$|V(X_s) - V(X_{s-})| \leq \Delta V,$$ where $0 \leq \Delta V < \infty.$
\itb There exists $\gamma <\infty$ such that
$$\langle M^V(x)\rangle_t \leq \gamma t$$ for all $x \in M_+, t \geq 0.$
\edes
\end{assumption}
Note that, by Proposition \ref{martconv}, a sufficient condition ensuring assertion $(b)$ is that $$\gamma : = \sup\{ \|\Gamma(V_K)\| \: : K \subset M_+, K \mbox{ compact } \} < \infty$$
For further reference, we will say that $\{(X_t^x)_{t \geq 0} \: x \in M_+\}$ is {\em $H$-persistent, strong version}
 if it satisfies Hypothesis \ref{hyp:H+}
 and is $H$-persistent.
 If additionally condition $(ii)'$ in Hypothesis \ref{hyp:H} is verified,
 we will say that it is  {\em $H$-persistent, strong version'}.


\subsection*{The case $M_0$ compact}
We first consider the situation where $M_0$ is compact. We let  $$M_0^{\delta} = \{x \in M_+ \: : d(x,M_0) < \delta\}$$ denote the $\delta$ neighborhood of $M_0.$

\bthm
\label{th:expoconvcompact}
Assume  that  $\{(X^x_t)_{t \geq 0} : x \in M_+\}$ is $H$-persistent (strong version), $M_0$ is compact,  $\tilde{W} = \alpha W$ for some $\alpha > 0$ (where $W$ and $\tilde{W}$ are like in Hypothesis \ref{hyp:tightpi}) and
that there exists a Doeblin point $x^* \in \Gam_{M_+}.$ Then, there exist $\lambda, \theta > 0, cst$ such that for all $x \in M_+$ and $f : M^+ \mapsto \RR$ measurable,
$$|P_t f(x) - \Pi f| \leq cst(1 + W_{\theta}(x)) e^{-\lambda t} \|f\|_{W_{\theta}},$$
where $W_{\theta}$ is continuous,  lies in $L^1(\Pi),$ and coincide with $V_{\theta} = e^{\theta V}$ on $M_0^{\delta}$ and with $W$ on $\{W > R\}$ for some $R > 0.$ Here $$\|f\|_{W_{\theta}} = \sup_{x \in M^+} \frac{|f(x)|}{1 + W_{\theta}(x)}.$$ In particular,
$$|P_t(x,\cdot) - \Pi|_{TV} \leq cst(1 + W_{\theta}(x)) e^{-\lambda t}$$
\ethm
\paragraph{Proof of Theorem \ref{th:expoconvcompact}. }
The following lemma follows from Proposition \ref{hajeklem} proved in Section \ref{sec:ratecompact}. Assertion $(ii)$ follows from \ref{pakeslem} and Remark \ref{rem:pakes}.
\blem
\label{hajeklem0}
There exist positive numbers $\theta, \kappa, T_0 < T_1,$ $0 < \rho < 1,$ and a continuous function $V_{\theta} : M_+ \mapsto \Rp$ such that
\bdes
\iti $V_{\theta} = e^{\theta V}$ on  $M_0^{\delta}$ for some $\delta > 0.$
\itii $\lim_{x \rar M_0} V_{\theta}(x) = \infty$ and  $V_{\theta}$ is bounded on $M_+ \setminus M_0^{\delta};$
\itiii For all $T \in [T_0,T_1], \, P_T V_{\theta} \leq \rho V_{\theta} + \kappa.$
\edes
 \elem
The proof of Theorem \ref{th:expoconvcompact} is now a consequence of a classical result  often refereed as "Harris's theorem" which proof can be found in numerous places
(e.g \cite{duf00}, \cite{MT1}). Here we rely on the following version  given (and proved) by Hairer and Mattingly \cite{Hairer-Mattingly}. Let ${\cal P}$ be a Markov kernel
on a measurable space $E.$
Assume that:
\bdes
\iti
 There exists a map ${\cal W}: E \mapsto [0, \infty[$ and constants $0 < \gamma < 1, \tilde{K} \geq 0$ such that
${\cal P} {\cal W} \leq \gamma {\cal W} + \tilde{K}$
\itii
 For some $R > \frac{2 \tilde{K}}{1-\gamma}$ there exists a probability measure $\nu$ and a constant $c$ such that
${\cal P}(x,.) \geq c \nu(.)$ whenever  ${\cal W}(x) \leq R.$
\edes
Then there exists a unique invariant probability $\Pi$ for ${\cal P},$  and constants  $0 \leq \tilde{\gamma} < 1, cst \geq 0$ such that for every measurable map
$f : E \mapsto \RR$ and all $x \in E$
$$|{\cal P}^n f(x)
 - \Pi f| \leq cst \: \tilde{\gamma}^n (1 + {\cal W}(x)) \|f
 \|_{\cal W}.$$ Here $\|f\|_{\cal W} = \sup_{x \in E} \frac{|f(x)|}{1 + {\cal W}(x)}.$
To apply this result, set $E = M_+$ and, using the notation of  Lemma \ref{hajeklem0}, ${\cal W} = V_{\theta} + W.$
Proposition \ref{hajeklem0} combined with the fact that $P_t W \leq e^{-\alpha t} W + C/{\alpha}$ (see Theorem \ref{tightlypnv} $(iii)$)  yield
\beq
\label{PTWleq}
P_{n T} {\cal W} \leq  \tilde{\rho}^n {\cal W} + \tilde{K}
\eeq
for all $n \in \NN$ and $T_0 \leq T \leq T_1,$
with $\tilde{\rho} = \max \{\rho, e^{-\alpha T_0}\}$ and $\tilde{K} = \frac{\kappa}{1-\rho} + \frac{C}{\alpha}.$

Choose $R > \frac{2 \tilde{K}}{(1-\tilde{\rho})^2}.$ The set ${\cal W}_R = \{x \in M_+ : \: {\cal W} \leq R\}$ is a compact subset of $M_+$ and by Lemma \ref{lem:compactsmall}, there exist some constants $T_R, c_R > 0$ depending on $R$ and a probability measure $\nu$ on $M_+$
- which we can assume to be supported by ${\cal W}_R$-
such that $P_{T_R}(x,.) \geq c_R \nu(.)$ for all $x \in {\cal W}_R.$ By iteration, this gives $P_{m T_R}(x,.) \geq c_R^m \nu(.)$ for all $x \in {\cal W}_R$ and $m \in \NN^*.$
Choose now $T \in [T_0, T_1]$ such that $T_R/T$ is rational, and positives integers $m,n$ such that $m/n = T_R/T.$
Thus $P_{n T} = P_{m T_R} = {\cal P}$ verifies  conditions $(i),(ii)$ above of Harris's theorem with $\gamma = \tilde{\rho}^n.$
The end of the proof is similar to the end of the proof of Theorem \ref{th:OreyTH}.

\qed
\subsection*{The case $M_0$ noncompact}
This section is strongly inspired by the recent beautiful work of Hening and Ngyuen \cite{HN18b} on Kolmogorov systems.
When $M_0$ is noncompact, the existence of a Lyapunov function $W$ controlling the behavior of the process at infinity, doesn't seem to be sufficient to ensure an exponential rate of convergence and one need to control the behavior of $H$ at infinity.

 We say that $\{(X^x_t)_{t \geq 0} : x \in M_+\}$ is {\em  $H$-persistent with respect to $M_0$ {\bf and} at infinity}, it is $H$-persistent and the maps $(V,H)$ of  Definition \ref{defHpers} satisfy the two following additional properties: \bdes
  \ita $V$ is proper;
  \itb There exists a compact $C \subset M$ such that
$$\sup_{x \in M \setminus C} H(x)  < 0.$$
\edes
\bthm
\label{th:expoconvnoncompact}
Assume  that
 $\{(X^x_t)_{t \geq 0} : x \in M_+\}$ is $H$-persistent (strong version') with respect to $M_0$ and at infinity and that there exists  a Doeblin point $x^* \in \Gam_{M_+}.$
Then there exists $\lambda > 0, \theta, cst$ such that for all $x \in M_+$ and $f : M^+ \mapsto \RR$ measurable,
$$|P_t f(x) - \Pi f| \leq cst(1 + W_{\theta}(x)) e^{-\lambda t}\|f\|_{W_{\theta}}.$$ In particular,
$$|P_t(x,\cdot) - \Pi|_{TV} \leq cst(1 + W_{\theta}(x)) e^{-\lambda t}$$
Here $W_{\theta} = e^{\theta V}$ and $\|f\|_{W_{\theta}}$ is like in Theorem  \ref{th:expoconvcompact}.
\ethm
The proof is given in section \ref{sec:ratenoncompact}.

The "construction" of $V$ (and $H$) ensuring persistence is (at least in all the examples we have in mind) dictated by our knowledge of the behavior of the process near the extinction set and there is, in general, no reason that the additional conditions ($a)$ and $(b)$ above ensuring persistence at infinity  are equally valid. The following simple result is a useful trick to get around this problem.
\bprop
\label{prop:HandtildeH}
Let
 $\{(X^x_t)_{t \geq 0} : x \in M_+\}$ be $H$-persistent (strong version') with respect to $M_0.$ Assume  that there  exist continuous functions
$(\tilde{V}, \tilde{H})$   satisfying   Hypothesis \ref{hyp:H+} with  condition $(ii)'$ of Hypothesis \ref{hyp:H} and such that:
  \bit
  \iti $\tilde{V}$ is defined  on {\bf all} $M$ and proper;
  \itii Conditions $(i)$ in Hypothesis \ref{hyp:H} (respectively $(b)$ in  Hypothesis \ref{hyp:H+}) are valid for {\bf every} compact subset of $M;$ \itiii  $\limsup_{x \rar \infty} \eps H(x) +   \tilde{H}(x) < 0,$ for some $\eps > 0.$
  \eit Then the process is $H$-persistent (strong version') at $M_0$ and at $\infty.$
\eprop
\prf First assume that for all $K \subset M^+$ compact, $V_K + \tilde{V_K} \in {\cal D}^2({\cal L}).$ Then
 $(\eps V + \tilde{V}, \eps H + \tilde{H})$ satisfies Hypothesis  \ref{hyp:H+} and  condition $(ii)'$ of Hypothesis \ref{hyp:H}. This easily follows from the  linearity of ${\cal L}$ and the property $\Gamma(f+g) \leq (\sqrt{\Gamma(f)} + \sqrt{\Gamma(g)})^2$ valid for $f,g \in {\cal D}^2({\cal L})$ with $fg \in {\cal D}({\cal L})).$ By remark   \ref{rem:whenlambdais0}, $\Lambda(\eps H + \tilde{H}) = \eps \Lambda(H) > 0.$
Hence the result.

In general (if we cannot argue that $V_K + \tilde{V_K} \in {\cal D}^2({\cal L})$)  note that (with the notation of Lemma \ref{vmartdefn}) $M_t^V, M_t^{\tilde V}$ being square integrable martingales, the same is true for $M_t^{V + \tilde{V}} = M_t^V + M_t^{\tilde{V}}$ and $\langle M^{V + \tilde{V}} \rangle_t \leq (\sqrt{\langle M^V_t \rangle} +  \sqrt{\langle M^{\tilde{V}} \rangle_t})^2$ and the proof goes through.
\qed

\section{Application to Ecological SDEs (ii)}
\label{sec:sde2}
 Consider the ecological SDE defined by (\ref{eq:sde}).
 Let $I \subset \{1, \ldots, n\}$ and
 $$M_0 := M_0^I = \{x \in M \:  : \prod_{i \in I} x_i = 0\}$$
be  the set corresponding to the extinction of at least one of the species $i \in I.$

Following \cite{Sch12}, \cite{SBA11}, \cite{BHW}, define the {\em invasion rate} of species $i$ with respect to $x$ as
\beq
\label{defsdelambda}
\lambda_i(x) = F_i(x) -  \alpha_i \frac{a_{ii}(x)}{2} x_i^{\alpha_i - 1}
\eeq
and the {\em invasion rate} of species $i$ with respect to $\mu \in \PR_{erg}(M_0^I)$ as
\beq
\label{defsdelambdamu}
\mu \lambda_i =  \int \lambda_i d\mu
\eeq
provided $\lambda_i \in L^1(\mu).$

The following result asserts that if a weighted combination of  the invasion rates $\{\mu  \lambda_i\}_{ i \in I}$ is positive for all $\mu \in \PR_{erg}(M_0^I),$ then the process is $H$- (hence stochastically) persistent. This criterion goes back to the early work of Hofbauer \cite{H81} (see also \cite{S00} and \cite{GH03}) but has been shown to apply also for SDEs, only recently, first in \cite{BHW} (for small noise), then in \cite{SBA11} (on compact state spaces) and recently in \cite{HN18b} (on $\RR^n_+$ for nondegenerate noise).

Note here that there is no assumption that the diffusion matrix (defined by (\ref{eq:defa})) is nondegenerate. We then retrieve Hofbauer's criterion, and - more importantly - this allows to handle situations where the "noise" only affect certain variables. Examples will be given in Section \ref{sec:degensde}.

\bthm
\label{persSDE} Let $U, \varphi$ and $\eta$ be as in Proposition \ref{prop:ecosde}.
 Assume that
\beq
\label{eq:integraSDE1}
\limsup_{x \rar \infty} \frac{  U^{\frac{1-\eta}{2}}(x) (1 + \varphi(x))}{1 + \sum_{i \in I} |F_i(x)|} = \infty
\eeq
and
\beq
\label{eq:integraSD2}
\sum_{i=1}^n x_i^{\alpha_i -1} \leq   \mbox{cst}\sqrt{U(x)}
\eeq
Then
\bdes
\iti
For all $\mu \in  \PR_{erg}(M_0^I)$ and $i \in I$  $\lambda_i \in L^1(\mu)$  and $$\mu \lambda_i \neq  0 \Rightarrow supp(\mu) \subset M_0^i = \{x \in M \: : x_i = 0\}.$$
\itii If there exist positive numbers  $\{p_i\}_{i \in I}$ such that for all $\mu \in  \PR_{erg}(M_0^I)$
\beq
\label{eq:weightpos}
\sum_{i \in I} p_i (\mu \lambda_i) > 0;
\eeq
Then the process is   $H$-persistent with respect to $M_0^I.$

\itiii If furthermore  $\eta = 0, \alpha_i = 1,$
 \beq
\label{eq:integraSDE2+}
1 + \varphi + \eps F_i  \geq 0
\eeq
for some $\eps > 0$ and $i \in I,$ and (\ref{eq:integraSDE1}) is strengthened to
 \beq
\label{eq:integraSDE1+}
\mid\frac{LU}{U}\mid^q + \sum_{i \in I} \mid F_i\mid^ q \leq cst \sqrt{U}
\eeq
for some $q > 1.$
Then, under  condition (\ref{eq:weightpos}), the process is $H$-persistent with respect to $M_0^I$ (strong version') and persistent at infinity.
\edes
\ethm
\prf $(i)$ Condition  (\ref{eq:integraSDE1})  combined with Theorem \ref{tightlypnv} and Proposition \ref{prop:ecosde} (v) imply that $\lambda_i \in L^1(\mu)$ for all $\mu \in \PR_{inv}(M).$ The second assertion will be proved after the proof of assertion $(ii).$

 $(ii)$ For all $i \in I$ let  $h_i(u) =  \log(\frac{1}{u})$ if $\alpha_i = 1,$ and  $h_i(u) = \frac{u^{1- \alpha_i}}{\alpha_i-1}$ if $\alpha_i > 1.$
Let $v : \RR \mapsto \Rp$ be a smooth function with bounded first and second derivatives such that $v(t) = t$ for $t \geq 1.$
Set
$$V(x) =  v(\sum_{i \in I } p_i h_i (x_i))$$  and
\beq
\label{eq:HforSde}
H(x) = v'(\sum_{i \in I } p_i h_i (x_i)) ( - \sum_i p_i \lambda_i(x)) + \frac{1}{2} v''(\sum_{i \in I } p_i h_i (x_i)) \langle a(x) p, p \rangle
 \eeq for $x \in M_+.$ Then $V$ (respectively $H$) coincide with $\sum_{i \in I } p_i h_i$ (respectively $\sum_{i \in I }   - p_i \lambda_i$ on the set $\{x \in M_+ \: : \sum_{i \in I } p_i h_i(x) > 1 \}$ and $H$ extends continuously to $M_0.$
Furthermore $$|H| \leq cst ( \sum_{i \in I }  p_i |\lambda_i| +  \sum_{i \in I } p_i^2)$$ so that, condition (\ref{eq:integraSDE1}), imply that condition $(ii)$ of Hypothesis \ref{hyp:H} is satisfied.
Let $b : \Rp \mapsto [0,1]$ be a smooth  function such that $b(t) = 1$ for $t \leq 1$ and  $b(t) = 0$ for $t \geq 2,$ and let $B(t) = \int_0^t b(u) du.$
For all $n \geq 1,$  set $$V_n(x) = n B(\frac{V(x)}{n}) b( \frac{\log( 1 + \sum_i x_i)}{n})$$  for $x \in M_+$ and
 $$V_n(x) = n B(2)  b( \frac{\log( 1 + \sum_i x_i)}{n})$$  for $x \in M_0.$ Then $V_n \in \DA^2(\LA)$ (because $V_n$ is  smooth with compact support),  $V_n = V$ and $\LA(V_n) = H$ on the set $K_n = \{x \in M_+ :\:  V(x) \leq n, \log(  1 + \sum_i x_i) \leq n\}.$
Furthermore $$\Gamma(V_n)(x)  \leq cst (1 + \frac{\sum_i x_i^{\alpha_i}}{1 + \sum_i x_i})^2 \leq
cst ( 1 + \sum_i x_i^{\alpha_i -1})^2 \leq cst (1 + U(x)) $$ so that assumption $(i)$ of Hypothesis \ref{hyp:H} is satisfied in view of assertion $(ii)$ of Proposition \ref{prop:ecosde}. This proves assertion $(ii).$ Also, by Lemma \ref{muHlemma}, $\mu H > 0 \Rightarrow \mu(M_0^I) > 0$ which concludes the proof of $(i).$

$(iii)$ When $\alpha_i = 1,$ $\Gamma(V_n)(x) \leq cst$ and Hypothesis \ref{hyp:H+} is satisfied. By condition (\ref{eq:integraSDE1+}), condition $(ii)'$ in  Hypothesis \ref{hyp:H} also holds, so that  the process is  $H$-persistent (strong version').

We will now apply Proposition \ref{prop:HandtildeH}.
Let $\tilde{V} = \log(U)$ and $\tilde{H} = L \tilde{U}.$
Then $$\tilde{H} = \frac{LU}{U} - \frac{1}{U^2} \Gamma_{L}(U) \leq -\alpha (1 + \varphi) + \frac{\beta}{U}.$$
Using the right hand side equality,   the assumptions on $U,$ and condition (\ref{eq:integraSDE1+}), it is easily checked that $(\tilde{V}, \tilde{H})$ satisfies Hypothesis \ref{hyp:H+} and condition $(ii)'$ of Hypothesis \ref{hyp:H}. It suffices to set $\tilde{V}_n = n B(\frac{\tilde{V}}{n})$ and to argue as previously.
From the left hand side inequality we get that
$$\limsup_{x \rar \infty}\tilde{H} + \alpha (1 + \varphi) \leq 0.$$
From condition (\ref{eq:integraSDE2+}) we get that
$$\sum_i p_i (1 + \varphi) + \eps \sum_i p_i F_i \geq 0.$$ Hence, replacing $\eps$ by a sufficiently smaller $\eps,$
$$2 \alpha (1 + \varphi) - \eps H \geq 0$$
where $H$ is defined by (\ref{eq:HforSde}).
This proves that $\limsup_{x \rar \infty}  \eps H + \tilde{H} < 0$ and Proposition \ref{prop:HandtildeH} applies.
\qed

\bex
\label{ex:Fbounded}
{\rm
Consider the system (\ref{eq:sde}) with $\alpha_i = 1.$ Suppose that the conditions of Proposition \ref{prop:ecosde} hold with $$U(x) = 1 + \sum_i x_i^p$$ for some $p > 1$ and  $\varphi = \eta = 0.$ Assume also that
$$\mid F_i(x)\mid^q \leq cst \sqrt{U}$$ for some $q > 1,$ and that $F_i$ is bounded from below (i.e~$\liminf_{x \rar \infty} F_i(x) > -\infty$). Then the assumptions (\ref{eq:integraSDE2+}) and (\ref{eq:integraSDE1+}) of Theorem \ref{persSDE} hold true.
}
\eex
\bex[continuation of example \ref{ex:LVex}]
\label{ex:LVex2}
{\rm
Using the notation of example \ref{ex:LVex}, assume that
$$\mid F_i(x) \mid \leq cst (1 + \parallel x \parallel)$$
and
$$f_i(x_i) = r - b x_i$$ for some $b > 0.$
Then the conditions (\ref{eq:integraSDE2+}) and (\ref{eq:integraSDE1+}) of Theorem \ref{persSDE} hold with
$$U(x) = 1 + \sum_i x_i^p$$ for all $p > 2$ and
$$\varphi(x) = \|x\|.$$
}
\eex
Call a point $x^* \in M$ {\em non-degenerate} if the matrix $a(x^*)$ is non-degenerate  or, equivalently, if $\{\Sigma^1(x^*), \ldots, \Sigma^m(x^*)\}$ span $\RR^n$.
\bcor
 \label{cor:ellipsde}
 Assume   $I = \{1, \ldots, n\}$ (so that  $M_0 = \partial \Rp^n$ and $M_{+} = int(\Rp^n)$) and that  the conditions (\ref{eq:integraSDE1}), (\ref{eq:integraSD2}), (\ref{eq:weightpos}) in Theorem \ref{persSDE} are satisfied. Assume furthermore that  there exists  $x^*  \in \Gam_{M_+} \cap M_+$  which is non-degenerate. Then
   \bdes
   \iti  $${\cal P}_{inv}(M_+) = \{\Pi\}$$ and
for all $x \in M_+$
$$\lim_{t \rar \infty} |P_t(x,\cdot) - \Pi|_{TV} = 0.$$
\itii Under the stronger conditions (\ref{eq:integraSDE2+}) and  (\ref{eq:integraSDE1+}), there exist positive constant $\eps, \theta, cst$ such that
$$|P_t f(x) -  \Pi f| \leq cst (1 + W_{\theta}(x)) e^{-\lambda t} \|f\|_{W_{\theta}}$$
with $$W_{\theta}(x) = \frac{U^{\theta}(x)}{(\prod_{i \in I} x_i^{p_i})^{\theta \eps}}.$$
\itiii If all the points in $M_+$ are non degenerate, then $\Gam_{M_+} \cap M_{+} = M_{+}.$
\edes
 \ecor
\prf $(i), (ii).$ The non-degeneracy of $a(x^*)$ makes $x^*$ a Doeblin point. Indeed, by Theorems 3.6 and
 3.7 in  Durrett \cite{Durret96}, Chapter7, relying on Dynkin \cite{Dynkin65},   there exist a open ball $D$ centered at $x^*$ and   a positive function $q_t(x,y)$ continuous in $t > 0, x, y \in D,$ such that if $f : \overline{D} \mapsto \RR$ is continuous with $f|_{\partial D} = 0$,
  \beq
  \label{doeblinSDE}
  \E_x(f(X_t) \Ind_{\tau > t}) = \int_{D} q_t(x,y) f(y) dy
   \eeq
   for all $t > 0, x \in D,$  where $\tau = \inf \{t \geq 0 \: : X_t \not \in D\}.$ Hence,  (\ref{eq:Doeblin}) holds with $\xi(dy) = \eps \Ind_{U}(y) dy,$ for some  $\eps > 0$ and $U \subset D$ a closed ball around $x^*$.  The result then follows from Theorems \ref{th:OreyTH} (respectively \ref{th:expoconvnoncompact}) and \ref{persSDE}.

   $(iii).$ For all $x \in M_+, \Gam_x$ is a closed set containing $x.$ If every point in $M_+$ is non-degenerate,  the proof of $(i)$ shows that $\Gam_x \cap M_+$ is open. Hence, by connectedness of $M_+,  \Gam_x \cap M_+ = M_+.$ \qed

\bex[Two dimensional systems]
\label{ex:2d}
{\rm
To illustrate the results above we consider here a simple model involving two species in interaction having the form
\beq \label{2dsde}
\left \{ \begin{array}{l}
  dx_1 = x_1(F_1(x) + \sigma_1(x) dB_t^1) \\
  dx_2 = x_2(F_2(x) + \sigma_2(x) dB_t^2)
  \end{array}\right.
\eeq
where $F_i, \sigma_i$ are smooth, $\sigma_i$ is positive and bounded,  $(B_t^1), (B_t^2)$  are two independent Brownian motions.
We furthermore assume that the assumption of Proposition \ref{prop:ecosde} as well as the conditions (\ref{eq:integraSDE1}, \ref{eq:integraSD2}) (or \ref{eq:integraSDE1+}, \ref{eq:integraSDE2+})  are satisfied (see for instance the examples \ref{ex:Fbounded} and \ref{ex:LVex2}).

 We let $M_+ = \{x \in M \: : x_1 > 0, x_2 > 0\}$
 and $M_0 = \{x \in M : \: x_1 x_2 = 0\}.$

We let
$$\lambda_i(x) = F_i(x) - \frac{\sigma^2_i(x)}{2}.$$
On the invariant face $x_2 = 0$  the process admits one ergodic probability  given by the dirac at the origin $\delta_{0,0}$ and an invariant measure (non necessarily a probability)  $h_1(x_1) dx_1  \delta_0(dx_2)$ where
$$h_1(x_1) = \frac{2}{x_1^2 \sigma_1^2(x_1,0)} \exp {\{\int_r^{x_1} \frac{2 F_1(u,0)}{u \sigma_1(u,0)^2} du\}} \Ind_{x_1 > 0} $$ and $r > 0$ is an arbitrary number.
This invariant measure is finite if only if the integrability condition
$$\frac{2 F_1(0,0)}{\sigma_1^2(0,0)} > 1 \Leftrightarrow \lambda_1(0,0) > 0$$ holds true.
Observe that this condition is exactly the persistence condition (\ref{persSDE}) of the process restricted to the face $\{x \in M :\: x_2 = 0\}$ with invariant set $\{0,0\}.$
In this later case we let $\mu_1$ denote the ergodic probability obtained by normalizing $h_1.$ That is
$$\mu_1(dx_1 dx_2) = \frac{ h_1(x_1)}{\int_0^{\infty} h_1(u) du} dx_1  \delta_0(dx_2).$$
We define $\mu_2$ similarly.
In summary,
 $$\Pr_{erg}(M_0) = \left \{ \begin{array}{l}
                      \{ \delta_{0,0} \} \mbox{ if } \lambda_1(0,0) < 0 \mbox{ and } \lambda_2(0,0) < 0,\\
                       \{\delta_{0,0}, \mu_1 \} \mbox{ if } \lambda_1(0,0) > 0 \mbox{ and } \lambda_2(0,0) < 0,\\
                        \{\delta_{0,0}, \mu_2 \} \mbox{ if } \lambda_1(0,0) < 0 \mbox{ and } \lambda_2(0,0) > 0,\\
                       \{\delta_{0,0}, \mu_1, \mu_2 \} \mbox{ if } \lambda_1(0,0) > 0 \mbox{ and } \lambda_2(0,0) > 0.
                    \end{array} \right.$$
 Therefore, the persistence condition (\ref{eq:weightpos}) is satisfied in one of the three following cases:
\bdes
\iti
 $\lambda_1(0,0) > 0,  \lambda_2(0,0) < 0,$ and  $\mu_1(\lambda_2) > 0,$ or
\itii $\lambda_1(0,0) < 0,  \lambda_2(0,0) > 0,$ and   $\mu_2(\lambda_1) > 0,$ or
\itiii  $\lambda_1(0,0) > 0,  \lambda_2(0,0) > 0,  \mu_1(\lambda_2) > 0,$ and  $\mu_2(\lambda_1) > 0.$
\edes
In each case the conclusions of Corollary \ref{cor:ellipsde} hold. (Compare to section 4.3 of \cite{SBA11} and to Section 2.1 of \cite{HN18b}).
}
\eex
\bex[Example \ref{ex:2d} continued, Randomness promotes persistence]
{\rm
Suppose that  the noise term  $\sigma_i$ in equation (\ref{2dsde}) takes the form
$$\sigma_i = \eps s_i$$ for some $0 < \eps << 1$ with $s_i > 0$ and bounded.

Let $V_1, V_2 : ]0, \infty[ \mapsto \RR$ be the maps defined by $$V_1(t) = \int_r^t \frac{F_1(u,0)}{u s^2_1(u,0)} du \mbox{ and }
V_2(t) = \int_r^t \frac{F_2(0,u)}{u s^2_2(0,u)} du.$$
 If $F_i(0,0) > 0,$ $V_i$ achieves its maximum at a point $x^*_i  > 0.$
 Assume for simplicity that such a maximum is unique. Note that $F_1(x^*_1,0) = 0$ (similarly $F_2(0,x^*_2) = 0$) and that
 $\frac{\partial F_1}{\partial x_1}(x_1^*,0) \leq 0$ (similarly $\frac{\partial F_2}{\partial x_2}(0,x_2^*) \leq 0$).

 If $F_1(0,0) > 0$ (respectively $F_2(0,0) > 0$) the probability $\mu_1$ (respectively $\mu_2$) converges when $\eps \rar 0$ (in the weak * sense) toward the dirac measure at $(x_1^*,0)$ (respectively $(0,x_2^*)$.)
 Therefore, using the results described in example \ref{ex:2d}, one see that the process is stochastically persistent for every $\eps > 0$ sufficiently small, provided one of  the following conditions hold:
 \bdes
 \iti $F_1(0,0) > 0, F_2(0,0) < 0 $ and $F_2(x_1^*,0) > 0,$ or
 \itii $F_1(0,0) < 0,  F_2(0,0) > 0$ and $ F_1(0,x_2^*) > 0,$ or
 \itiii $F_1(0,0) > 0,  F_2(0,0) > 0, F_2(x_1^*,0) > 0$ and $F_1(0,x_2^*) > 0$
 \edes
An interesting consequence of this result is that an arbitrary small random perturbation of a non persistent deterministic ecological ODE can be stochastically persistent. Indeed condition $(i)$ above simply means that the origin is a saddle point (for the ode obtained with $\eps = 0$) which stable (respectively unstable) manifold is the axis $x_1 = 0$ (respectively
$x_2 = 0$) and the point $(x_1^*,0)$ a saddle point which stable manifold is the axis $x_2 = 0.$ However there may exist other equilibria on the boundary including sinks or saddle points.
}
\eex
\subsection{Degenerate ecological SDEs}
\label{sec:degensde}
We discuss here  the situation where  (\ref{eq:sde}) is a degenerate SDE. This is motivated by  models for which the noise  only affects certain variables. We assume  throughout the section that  the vector fields $F$ and $\Sigma^j, j = 1, \ldots, m$ are $C^{\infty}$.

 Rewrite the  stochastic differential equation (\ref{eq:sde}) using the Stratonovich formalism  as
\beq
\label{eq:sdestrato}
dx_t = S^0(x_t) dt + \sum_{j = 1}^m S^j(x_t) \circ dB^j_t
\eeq
where for all $j = 1 \ldots m,$
$$S^j_i(x) = x_i^{\alpha} \Sigma_i^j(x), i = 1 \ldots n$$
and
$$S^0_i(x) = x_i^{\alpha_i} F_i(x) - \frac{1}{2} \sum_{j = 1}^m \sum_{k = 1}^n \frac{\partial S_i^j}{\partial x_k}(x) S_k^j(x), i = 1 \ldots n.$$
Associated to this system is the deterministic control system
\beq
\label{controlsde}
\dot{y}(t) = S^0(y(t)) + \sum_{j = 1}^m u^j(t) S^j(y(t))
 \eeq where the control function   $u = (u^1, \ldots, u^m) : \Rp \mapsto \RR^m, $ can be chosen to be piecewise continuous.
Given such a control function, we let $y(u,x, \cdot)$ denote the  maximal solution\footnote{Note that there is no assumption here that the vector fields $S^j$ are globally integrable.} to (\ref{controlsde}) starting from $x$ (i.e~ $y(u,x,0) = x$).
The following proposition easily follows from the
 the celebrated Strook and Varadhan's support theorem \cite{StrVar72} (see also Theorem 8.1, Chapter VI in \cite{Ikeda}.
 Recall that we let  $\Gam_x$ denote the accessible set from $x$ (as defined in section \ref{sec:support}).
 \bprop
 \label{prop:control}
Let $x \in M.$ Point $p \in M$ lies in $\Gam_x$ if and only if for every neighborhood $O$ of $p$ there exist a control $u$ such that $y(u,x,\cdot)$ meets $O$ (i.e $y(u,x,t) \in O$ for some $t \geq 0$).
\eprop
\prf
If the vector fields $S^j$ were bounded with bounded (first and second) derivatives, this would follow directly from the support theorem (see Theorem 8.1, Chapter VI in \cite{Ikeda}). To handle the fact that the $S^j$ are typically unbounded we use a localization argument relying on the existence of the Lyapunov function $U$ assumed in proposition \ref{prop:ecosde}. Let $S^{j,n}$ be a smooth vector field with compact support coinciding with $S^j$ on the set $U_n = \{x \in M \: : U(x) < n\}.$   Let $\mathbb{P}^n_x$ be for the law of the process starting from $x$ solution to the SDE obtained by replacing $S^j$ by $S^{j,n}$ in (\ref{eq:sdestrato}). Let $y^n(x,u,\cdot)$ be defined like $y(x,u,\cdot)$ when $S^j$ is replaced by $S^{j,n}$ in (\ref{controlsde}). Let $\tau_n = \inf\{t \geq 0 :\: U(X_t) \geq n\}.$ The assumptions on $U$ imply that $\lim_{n \rar \infty} \mathbb{P}_x(\tau_n > t) = 1$ (see the proof of Proposition \ref{prop:ecosde}, equation (\ref{eq:nonexplos})). Thus, for every open set $O \subset M$,
$$\mathbb{P}_x( X_t \in O) > 0 \Leftrightarrow \exists n  \: \:   \mathbb{P}_x( X_t \in O; \tau_n > t) > 0 \Leftrightarrow \exists n \: \: \mathbb{P}^n_x( X_t \in O; \tau_n > t) > 0.$$
By the support theorem, $$\mathbb{P}^n_x( X_t \in O; \tau_n > t) > 0 \Leftrightarrow
\exists u \:  y^n(u,x,[0,t]) \subset U_n \mbox{ and } y^n(u,x,t) \in O$$
 $$\Leftrightarrow
\exists u \:  y(u,x,[0,t]) \subset U_n \mbox{ and } y(u,x,t) \in O$$ because  $S^j = S^{j,n}$ on $U_n.$ Therefore
$$\mathbb{P}_x( X_t \in O) > 0 \Leftrightarrow \exists u \: y(u,x,t) \in O.$$ This proves the result.
\qed
The local {\em ellipticity} condition given by the non degeneracy of $a(x^*)$ in Corollary \ref{cor:ellipsde} can be weakened and replaced by a local {\em hypoellipticity} condition.

Recall that the Lie bracket of two smooth vector fields $Y, Z : \RR^n \mapsto \RR^n$ is the vector field defined as
 $$[Y,Z](x) = DZ(x) Y(x) - D Y(x) Z (x).$$
 Given a family $\mathcal{X}$ of smooth vector fields on $\RR^n,$ we let $[\mathcal{X}]_k, k \in \NN,$ and $[\mathcal{X}]$ denote
 the set of vector fields defined  by $[\mathcal{X}]_0 = \mathcal{X},$  $$[\mathcal{X}]_{k+1} = [\mathcal{X}]_k \cup \{  [Y,Z] \: : Y, Z \in [\mathcal{X}]_k\}$$ and  $[\mathcal{X}] = \cup_k [\mathcal{X}]_k.$ We also let $[\mathcal{X}](x) = \{Y(x) \: : Y \in [\mathcal{X}] \}.$
\\

Consider again the SDE (\ref{eq:sde}) (or equivalently \ref{eq:sdestrato}). We say that $x^* \in M$    satisfies the  {\em Hörmander condition} (respectively the {\em strong Hörmander condition})  if $[\{S^0, \ldots, S^m\}](x^*)$ (respectively
 $${\{S^1(x^*),  \ldots, S^m(x^*)\} \cup \{[Y,Z](x^*) :\: Y, Z \in [\{S^0, \ldots, S^m\}] \}})$$ spans  $\RR^n.$

The next corollary just states that the local ellipticity condition in Corollary \ref{cor:ellipsde} can be weakened to a local hypoellipticity.
\bcor
 \label{cor:hyposde}
 Assume   $I = \{1, \ldots, n\},$  so that  $M_0 = \partial \Rp^n$ and $M_{+} = int(\Rp^n)$, and that the conditions (\ref{eq:integraSDE1}, \ref{eq:integraSD2}) and (\ref{eq:weightpos}) of Theorem \ref{persSDE} are satisfied.
   \bdes
   \iti If  there exists  $x^*  \in \Gam_{M_+} \cap M_+$  which satisfies the Hörmander condition, then  $${\cal P}_{inv}(M_+) = \{\Pi\},$$
   where $\Pi << \lambda$ (the Lebesgue measure on $\RR^n$),
$\Pi^x_t \Rightarrow \Pi$ $\Pr$ a.s for all $x \in M_+$,  and
$\lim_{t \rar \infty} \Pi_t^x f = \Pi f$ a.s for all $f \in L^1(\Pi)$ and $x \in M^+.$
\itii If the Hörmander condition at $x^*$ is strengthened to the strong Hörmander condition,
then $(P_t)$ converge to $\Pi$ in total variation (like in  Corollary \ref{cor:ellipsde} $(i)$). Under the stronger conditions (\ref{eq:integraSDE2+}) and (\ref{eq:integraSDE1+}), the convergence is exponential ((like in  Corollary \ref{cor:ellipsde} $(ii)$).
\itiii If for all $x \in M_+$ $[\{S^1, \ldots S^m\}](x)$ spans $\RR^n,$ then $\Gam_{M_+} \cap M_{+} = M_{+}.$
\edes \ecor
\prf $(i)$ Let $D$ be a domain (connected open set) containing $x^*$, relatively compact, and small enough so that $[\{S^0, \ldots, S^m\}](x)$ spans $\RR^n$ for each $x \in \bar{D}.$ First assume that
\bdes
\ita For each $x \in \bar{D}$ $\sum_{i = 1}^m \|S^i(x)\| \neq 0;$
\itb For each $x \in \partial D = \bar{D} \setminus D$ there exists a vector $u$ normal to $\bar{D}$ such that
$\sum_{i = 1}^m \langle S^i(x),u\rangle^2 > 0.$
\edes
Under these assumptions, by  a Theorem of Bony (\cite{Bony}, Theorem 6.1),  there exists a kernel  $G : {\bar D} \times {\bar D}\mapsto \Rp,$ smooth on $  D \times D \setminus \{(x,x) \: : x \in D\}$ such that:  For each $f \in C_b(\bar{D}),$  there exists a unique $g \in C_b(\bar{D})$ solution to   the Dirichlet problem
$$\left \{ \begin{array}{l}
  L g - g  = - f \mbox{ on } D (\mbox{ in the sense of distributions}) \\
  g|_{\partial D}  =  0;
\end{array} \right.$$ and $g(x) = Gf(x):= \int G(x,y) f(y) dy.$ Furthermore, if $f$ is smooth on $D$ so is $g.$

Note that, by continuity of $G$ off the diagonal, there exist disjoint open sets $U, V \subset D,$ with $x^* \in U$ and $\delta > 0$ such that $G(x,y) \geq \delta$ on $U \times V.$

Let $\tau = \inf \{ t > 0 \: X_t \not \in D\}.$ For $f$ smooth on $D,$ Ito's formula shows that, $$(e^{-{t \wedge \tau}} g(X_{t \wedge \tau}) +\int_0^{t \wedge \tau} e^{-s} f(X_s) ds)$$ is a local martingale. Being bounded, it is  a uniformly integrable martingale. Thus,   $$\E_x(\int_0^{\tau} e^{-s} f(X_s) ds) = Gf(x).$$
Let $R(x,\cdot)  = \int_0^{\infty} e^{-t} P_t(x,\cdot).$  It follows that  for all $x \in U$
 $$R(x,dy)  \geq \delta \Ind_V(dy)$$
 and the result follows from Proposition \ref{th:petite}.

 It remains to explain how we can choose $D$ to ensure  that conditions $(a)$ and $(b)$ above are satisfied. We assume here that $n \geq 2.$ For $n = 1$ the proof is left to the reader.
 If $\sum_{i \geq 1} \|S^i(x^*)\| > 0,$ then $(a)$ holds provided $D$ is small enough. If  $\sum_{i \geq 1} \|S^i(x^*)\| =  0,$  set $y^* = \Phi^0_t(x^*)$ where $\{\Phi^0_t\}$ is the local flow induced by $S^0.$  We claim that, for $t > 0$ small enough, $y^* \in D$ and
  $\sum_{i \geq 1} \|S^i(y^*)\| >  0.$ Since $y^*$ is accessible, it then suffices to replace $x^*$ by $y^*$ and  $D$ by a neighborhood of $y^*.$ To prove this claim, assume to the contrary, that $S^i(\Phi_t^0(x^*)) = 0$ for all $0 < t < \eps$ and $i = 1, \ldots, m.$ Then
 $$0 = DS^i(\Phi_t^0(x^*)) \frac{d}{dt} \Phi^0_t(x^*) = DS^i(y^*) S^0(y^*) = [S^0,S^i](y^*).$$ Similarly $Z(y^*) = 0$ for all $Z \in \{[S^0, \ldots, S^m]\} \setminus \{S^0\}.$
  A contradiction.

   For condition $(b),$ we can assume (by condition $(a)$) that $S^1(x^*) \neq 0$ and without loss of generality that $\frac{S^1(x^*)}{\|S^1(x^*)\|}  = e_1$ the first vector in the canonical basis of $\RR^n.$
  Let, for $\eps > 0,$ small enough  $D_{\eps} = \{x \in \RR^n \: \|x-x^*\|_1 < \eps \}$ where $\|u\|_1 = \sum_{i = 1}^n |u_i|.$
  For $x \in \partial D_{\eps}$ let $u_x$ be the vector defined by $u_{x,i} = \frac{x_i - x^*_i}{|x_i - x^*_i|}$ if $x_i \neq x^*_i$ and
  $u_i^* = 1$ otherwise. Vector $u_x$ is normal to $\partial D$ and $\langle e_1, u_x \rangle^2 = 1.$ Hence, for $\eps$ small enough
  $\langle S^1(x), u_x \rangle^2 > 0$ for all $x \in \partial D.$ It suffices to replace $D$ by $D_{\eps}$ and $(b)$ is satisfied.

  $(ii)$
Under the strong Hörmander condition, the law of $(X_t)$ killed at $D$ (see Ichihara and Kunita \cite{Ichihara}) has a density $q_t(x,y)$ which is $C^{\infty}$ in $t > 0, x,y \in D.$  Choose $y^* \in D$ and $t^* > 0$ such that $q_{t^*}(x^*,y^*) > 0$ (such a $(t^*,y^*)$ exists for otherwise $\tau$ would be almost surely $0$ contradicting the continuity of paths). The end of the proof is then identical to the proof of Corollary \ref{cor:ellipsde}.

$(iii)$ follows from Chow's Theorem.
\qed
\brem
\label{rem:hyposde}
{\rm  In case all the points in $M_+$ satisfy the Hörmander condition, the invariant measure $\Pi$ in the previous corollary has a  $C^{\infty}$ density by hypoellipticity of $L^*$ (the formal adjoint of $L$). }
\erem
\subsection{A stochastic Rosenzweig-MacArthur model}
 \label{sec:RMAsde}
As an illustration of the previous result, we
consider here the Rosenzweig-MacArthur model described in Example \ref{ex:RMA},
under the assumption that only the prey-variable is subjected to some small environmental fluctuation. That is
\beq
\label{eq:RMAstoc}
\left\{ \begin{array}{l}
  dx_1 = x_1 (F_1(x_1,x_2) dt  + \eps dB_t)\\
  dx_2  = x_2 F_2(x_1,x_2) dt \end{array}
  \right.
\eeq
where $$F_1(x) = 1- \frac{x_1}{\kappa} -  \frac{x_2}{1+x_1}, F_2(x) = - \alpha + \frac{x_1}{1+x_1}, \: \alpha, \kappa > 0.$$
We let $M_0 = \partial \RR^2_+, M_+ = Int(\RR^2_+).$

For $0 < \eps^2 < 2,$ let $k = \frac{2}{\eps^2}-1, \theta = \frac{\eps^2 \kappa}{2}$ and let $$\gamma_{\eps,\kappa}(x) = \frac{x^{k-1} e^{-\frac{x}{\theta}}}{\Gamma(k) \theta^{k}}$$ be the density of a $\Gamma$ distribution with parameters $k, \theta.$
Set $$\Lambda(\eps,\alpha,\kappa) = \int_0^{\infty} \frac{x}{1+x} \gamma_{\eps,\kappa}(x) dx - \alpha.$$
A rough estimate of $\Lambda(\eps,\alpha,\kappa)$ is
$$ \frac{ \kappa (1-\frac{\eps^2}{2})}{1 + \kappa (1-\frac{\eps^2}{2})} - \kappa \frac{\eps}{4} \sqrt{1-\frac{\eps^2}{2}} -\alpha
\leq \Lambda(\eps,\alpha,\kappa) \leq  \frac{ \kappa (1-\frac{\eps^2}{2})}{1 + \kappa (1-\frac{\eps^2}{2})} -\alpha.$$
The right hand side inequality follows from Jensen inequality and the fact that  $\gamma_{\eps,\kappa}$ has mean $k \theta = \kappa (1-\frac{\eps^2}{2}).$  The left hand side follows from the fact that $x \mapsto \frac{x}{1+x}$ is $1-$ Lipschitz, Cauchy Schwarz  inequality and the fact that  $\gamma_{\eps,\kappa}$ has variance $k\theta^2 = \kappa^2 \frac{\eps^2}{2}(1-\frac{\eps^2}{2}).$
\bthm
\label{th:RMAstoc}  System (\ref{eq:RMAstoc}) behaves as follows:
\bdes
\iti If $\Lambda(\eps,\alpha,\kappa) > 0$ (in particular $0 < \eps^2 < 2$), then
 for all $x \in M^+$ $(\Pi_t^x)$ (respectively $( P_t(x,\cdot) )$ converges almost surely (respectively in total variation) toward a unique measure $\Pi$ (depending on $\eps$). Furthermore,   $\Pi$ has a smooth density (with respect to Lebesgue measure) strictly positive over $M_+.$
\itii If  $\Lambda(\eps,\alpha,\kappa) < 0$ (in particular $0 < \eps^2 < 2$), then $x_2(t) \rar 0$ a.s and for all $x \in M^+$ $x(t) \Rightarrow \gamma_{\eps,\kappa}(dx_1) \otimes \delta_0(dx_2).$
\itiii If $\eps^2 > 2,$ then  $x(t) \rar 0$ almost surely.
\edes
\ethm
\prf We only prove $(i).$ Assertion $(ii)$ and $(iii)$ will be proved in Part II.
Fix $n > 2.$ We first notice that the assumptions of Proposition \ref{prop:ecosde} are satisfied with $U(x) = (x_1 + x_2)^n, \eta = 0$ and $\varphi = 0.$
Indeed,
$$LU(x) =  n[(x_1 + x_2)^{n-1}(x_1(1-\frac{x_1}{\kappa}) - \alpha x_2) + \frac{n-1}{2} (x_1 + x_2)^{n-2} \eps^2 x_1^2]$$
$$\leq n(x_1 + x_2)^{n-1} [x_1(1 + \frac{n-1}{2} \eps^2 -\frac{x_1}{\kappa}) - \alpha x_2]$$
$$\leq  n(x_1 + x_2)^{n-1} (-\alpha (x_1 + x_2)  + \beta)$$ for some $\beta > 0$ .
Thus
$$LU \leq -\alpha U + \tilde{\beta}$$ for some $\tilde{\beta} > 0.$
Also
$$\Gamma_L(U)(x) = (n(x_1 + x_2)^{n-1} \eps x_1)^2 \leq n^2 \eps^2 U^2(x).$$
With such a function $U$ the conditions (\ref{eq:integraSDE1}) and (\ref{eq:integraSD2})    of Theorem  \ref{persSDE} are clearly   satisfied.

Reasoning like in  Example \ref{ex:2d} we see that for $\eps^2 < 2,$ $\Pr_{erg}(M_0) = \{\delta_{0,0},\mu_1\}$
where $\mu_1(dx_1 dx_2) = \gamma_{\eps,\kappa}(x_1) dx_1 \delta_0(dx_2)$ and that the persistence condition (\ref{eq:weightpos}) is given by $\Lambda(\eps,\alpha,\kappa) > 0.$

 The vector fields $S^0, S^1$ (in equation (\ref{eq:sdestrato}))  write
 $$S^0(x_1,x_2) = \left(
                    \begin{array}{c}
                      x_1 (F_1(x_1,x_2) - \eps^2/2))  \\
                      x_2 F_2(x_1,x_2) \\
                    \end{array}
                  \right); \: S^1(x_1,x_2) = \left(
                                \begin{array}{c}
                                  x_1 \eps \\
                                  0 \\
                                \end{array}
                              \right).$$
                               Thus, by a simple computation, $$det( [S^0, S^1](x), S^1(x)) = \eps^2 x_1^2 x_2 \frac{\partial F_2}{\partial x_1} = \frac{\eps^2 x_1^2 x_2}{(1+ x_1)^2}.$$
 This shows that the strong Hörmander condition holds at every point $x \in M_+.$
 To prove the claim, it remains  to show that $\Gam_{M_+}$ contains $M_+.$ Existence and convergence to $\Pi$ will then follow from Corollary \ref{cor:hyposde} $(ii).$ Smoothness of the density from Hypoellipticity of the adjoint $L^*$ (see e.g the reasoning in \cite{Ichihara} before the proof of Proposition 5.1) and positivity of the density from the fact that  $\Gam_{M+}$ is the support of $\Pi.$

 Introduce the new control variable $v = \eps u - \eps^2/2.$ Then, the  control system (\ref{controlsde}) rewrites
 $$\dot{x_1} = x_1 (F_1(x) + v), \dot{x_2} = x_2 F_2(x).$$
 Let $L$ be the line $x_1 = \frac{\alpha}{1-\alpha}$. That is $F_2(x) = 0.$ Let $(P_v)$ be the parabola $(1 + v-\frac{x_1}{\kappa})(1+x_1) = x_2.$ That if $F_1(x) + v = 0.$ Let $z = (z_1,z_2) \in M_+$ and $O_z$ a neighborhood of $z.$ Choose $v^*$ large enough so that $z$ is below $(P_{v^*})$ and the point at which $(P_{v^*})$ reaches
 its maximum is $ > \frac{\alpha}{1-\alpha}.$
 It is not hard to verify that there exists a neighborhood $O$ of the origin such that for all $x \in O \cap M_+$ $t \mapsto y(v^*,x,t)$ crosses $(P_{v^*})$ near $(v^*,0)$ then remains above $(P_{v^*})$ and then crosses $(L).$ In particular, it crosses the line $x_2 = z_2.$
 Given any $x \in M_+$ use a piecewise constant  control $v(t)$ as follows: $v(t) = -1$ until $t \mapsto y(v,x,t)$ enters $O.$ Then  $v(t) = v^*$ until $t \mapsto y(v,x,t)$ crosses the horizontal line $x_2 = z_2.$ Then $v(t) = -R$ where $R$ is large enough so that $t \mapsto y(v,x,t)$ eventually enters $O_z.$ By  Proposition \ref{prop:control} this proves that $z \in \Gam_x.$
\qed
\newpage
Figure \ref{fig:RMA1} obtained in Scilab by Edouard Strickler, illustrates the behavior of the process when $\Lambda(\eps,\alpha,\kappa) > 0$. The red trajectory is a trajectory of the unperturbed system.
\begin{figure}
\centering
\includegraphics[width=16cm]{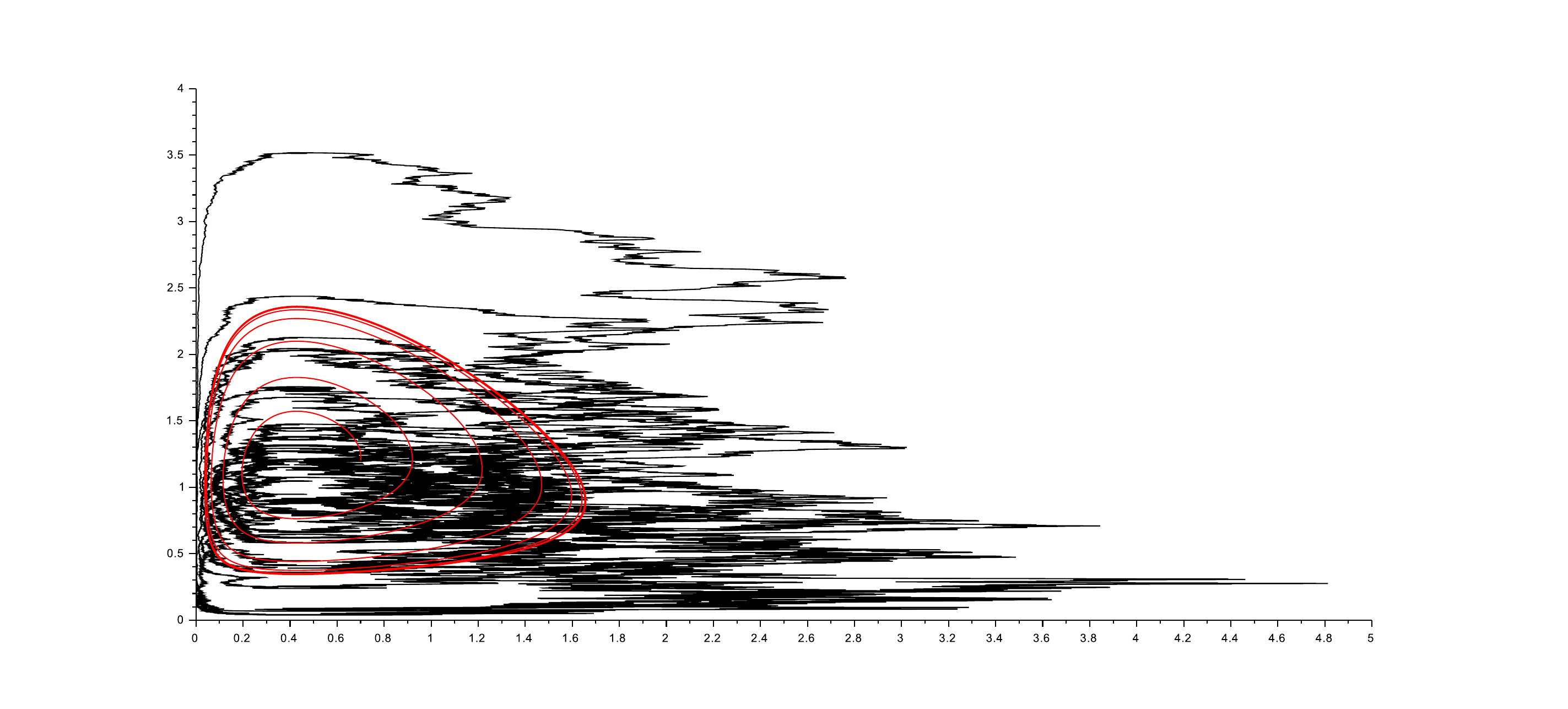}\\
\caption{$\alpha = 0,3; \kappa = 2,5: \eps = 0,6$ \label{fig:RMA1}}
\end{figure}
\section{Application to Random Ecological ODEs (ii)}
\label{sec:pdmpii}
Here we consider the random ecological ODEs introduced in section  \ref{sec:PDMP} and  use the same notation.
Recall that the state space has the form $M = B \times \{1,\ldots,m\},$  with $B$ compact, and that for $I \subset \{1,\ldots,n\}$ we let
$$M_0^I = \{(x,u)  \in M \:: \prod_{i \in I} x_i = 0 \}.$$ The {\em invasion rate} of species $i$ with respect to $(x,u)$ is defined as
$\lambda_i(x,u) = F_i^j(x)$ and
 invasion rate of species $i$ with respect to $\mu \in \PR_{erg}(M_0^I)$ as
\beq
\label{ratePDMP}
\lambda_i(\mu) =  \sum_j \int_{B} F_i^j(x) d\mu^j(x)
\eeq where
 $\mu^j(A)  :=  \mu (A \times \{j\}).$
\bthm
\label{persPDMP}
Assume that  that there exist positive number  $\{p_i\}_{i \in I}$ such that for all $\mu \in  \PR_{erg}(M_0^I)$
$$\sum_{i \in I} p_i \lambda_i(\mu) > 0.$$
Then the process  given by (\ref{eq:PDMP}) is  $H$ persistent with respect to $M_0^I.$
\ethm
\prf The proof is similar to the proof of Theorem \ref{persSDE}. If $(x,u) \in M_+$ (respectively $(x,u) \in M$) let $V(x)$ (respectively $V_n(x)$) be defined exactly as in the proof of Theorem \ref{persSDE}. See $V$ and $V_n$ as functions of $(x,u)$ (i.e set $V(x,u) := V(x), V_n(x,u):= V_n(x)$)  and let $H(x,u)= - \sum_{i \in I} p_i   \lambda_i(x,u).$
Then $V_n \in \DA^2(\LA),$  $V_n = V$ and $\LA(V_n) = H$ on the set $K_n = \{(x,u) \in M_+ :\:  V(x,u) \leq n, \sum_i x_i \leq n\}.$
Furthermore $\Gamma(V_n)(x,u) = 0$ (because $V_n(x,u)$ doesn't depend on $u$) so that assumption $(i)$ of definition \ref{hyp:H} is satisfied
 \qed
Associated to (\ref{eq:PDMP}) is the control  system
\beq
\label{eq:controlode}
\dot{y}(t) = \sum_{j = 1}^m u^j(t) G^j(y(t))
\eeq
where the control function $u = (u^1,\ldots,u^m) : \Rp  \mapsto \RR^n$  can be chosen to be piecewise continuous with values either in $\{e_1, \ldots, e_n\},$ the canonical basis of $\RR^n;$ or in $\Delta^{n-1},$ the unit simplex of $\RR^n.$ The solution to (\ref{eq:controlode}) starting from $x$ is denoted $t \mapsto y(u,x,t).$

The following proposition is analogous to Proposition \ref{prop:control} and follows  from the support theorem established in (\cite{BMZIHP}, Theorem 3.4). Note that this support theorem is phrased in terms of the differential inclusion whose set valued vector field is given by the convex hull of $\{G^1, \ldots, G^m\}$ but the link with the control system (\ref{eq:controlode}) is spelled out, for instance, in (\cite{BCL16}, Theorem 2.2).
 \bprop
 \label{prop:controlode}
Let $(x,i) \in M.$ Point $(p,j) \in M$ lies in $\Gam_{(x,i)}$ if and only if for every neighborhood $O$ of $p$
 there exists a control $u$ such that $y(u,x,\cdot)$ meets $O$ (i.e $y(u,x,t) \in O$ for some $t \geq 0$).
\eprop
By analogy with the terminology used for SDE's in section \ref{sec:degensde}, we say that $(x,i) \in M$ satisfies the {\em Hörmander} or {\em weak bracket} (the terminology coined in \cite{BMZIHP}) condition, respectively the {\em strong  Hörmander} or {\em strong bracket condition} if
$[\{G^1, \ldots, G^m\}](x),$  respectively $$\{G^i(x)-G^j(x) \: : i,j = 1, \ldots, m\} \cup \{[Y,Z](x) :\: Y, Z \in [\{G^1, \ldots, G^m\}] \}$$ spans  $\RR^n.$ These two  conditions are named  $A$ (for the stronger) and $B$ (for the weaker) in \cite{bakhtin&hurt}.
\bcor
 \label{cor:hypoode}
 Assume   $I = \{1, \ldots, n\},$  and that the  condition of Theorem \ref{persPDMP} holds. Assume that there exists  $(x^*,j)  \in \Gam_{M_+} \cap M_+$  which satisfies the Hörmander condition. Then
   \bdes
   \iti  $${\cal P}_{inv}(M_+) = \{\Pi\},$$
   where $\Pi << \lambda$ (the Lebesgue measure on $M$), and
for all $(x,j) \in M_+$
$$\lim_{t \rar \infty} \Pi^{(x,j)}_t =  \Pi$$ $\Pr$ a.s
\itii Assume in addition that, either
\bdes
\ita  the Hörmander condition at  $(x^*,j)$ is strengthened to the strong Hörmander condition, or
\itb There exist $\alpha_j, \ldots, \alpha_m \in \RR$ with $\sum_j \alpha_j = 1$ and $e \in \Gam_{M^+} \cap M^+$ for which $$\sum_{j = 1}^m \alpha_j G^j(e) = 0.$$
\edes Then for all $(x,j) \in M_+$
$$\lim_{t \rar \infty} |P_t((x,j),\cdot) - \Pi|_{TV} \leq cst(1 + W_{\theta}(x))e^{-\lambda t}.$$
\edes for some $\theta, \lambda > 0.$  Here $$W_{\theta}(x) = e^{\theta [\max (\sum_i p_i h_i(x_i), 1)]}$$ where $h_i(u) = -\log(u)$ for $\alpha_i = 1$ and $h_i(u) = \frac{u^{1-\alpha_i}}{1-\alpha_i}$ if $\alpha_i > 1.$
\ecor
\prf
Under condition $(i)$ ${\cal P}_{inv}(M_+)$ has at most cardinal one, as proved in \cite{bakhtin&hurt} for constant rates $a_{ij}$
and in \cite{BMZIHP} for nonconstant rates $(a_{ij}(x)).$ Note that for constant rates \cite{bakhtin&hurt} actually prove that the condition of Proposition \ref{th:petite} are satisfied with $\gamma(dt) = e^{-t} dt.$ The result then follows from  Corollary  \ref{cor:posinv}. Under condition $(ii) a)$, it follows again from \cite{bakhtin&hurt} (for constant rates) and \cite{BMZIHP} for nonconstant rates that $x^*$ is a Doeblin point. Under condition $(ii), b)$ this follows from a result recently proved in \cite{BenHurthStric} strongly inspired by the work of \cite{li17}. The result then follows from Theorem \ref{th:expoconvcompact}.
\qed
\brem
\label{rem:hypoPDMP}
{\rm In contrast with the SDE situation (see remark \ref{rem:hyposde}), the question of the smoothness of the density  of the invariant measure remains largely an open problem (see \cite{bakhtin&hurt&matt}, \cite{2017arXiv170801390B} for some results in dimension one and two).}
\erem
\subsection{A May-Leonard System with random switching}
\label{sec:MayLeonard}
The goal of this section is twofold: Illustrate the preceding results and provides a simple example, albeit non trivial, for which the H-persistence machinery applies to a situation where the extinction set is not just the boundary of $\Rp^n.$

Let $r: \Rp^3 \mapsto \Rp^*$ be a smooth function, $(\alpha,\beta)$ a pair or parameters - called an {\em environment} -satisfying
\beq
\label{def:enviML}
 0 < \beta < 1 < \alpha,
  \eeq
 and let  $G$ be the  vector field on $\Rp^3$ defined as $G_i(x) = x_i F_i(x),$ with
\beq
\label{eq:may} F(x) = r(x)
\left \{ \begin{array}{l}
(1-x_1-\alpha x_2 - \beta x_3) \\
(1-\beta x_1 - x_2 - \alpha x_3) \\
(1-\alpha x_1 - \beta x_2 - x_3)
       \end{array} \right.
\eeq
 When $r:= 1$ we recover the celebrated model introduced by May and Leonard \cite{MayLeonard} in 1975. A nonconstant $r$ has no effect on the phase portrait of $G$ (it only changes the velocity) but will have some  on the persistence properties of the process obtained by random switching of the parameters.

Before considering such a process,  we first recall some basic properties of the dynamics  induced by $G.$
 \paragraph{Background} We let  $\Phi = \{\Phi_t\}$ denote the local solution flow in $\Rp^3$ to the differential equation $\dot{x} = G(x).$ Here, the terminology {\em trajectories, equilibria, limit points} etc. refer to trajectories, equilibria, limit points of $\Phi.$
Throughout we let
$$N(x) = x_1 + x_2 + x_3$$ and $$\Delta = \{x \in \Rp^3 \: : N(x) = 1\}$$ denote the unit simplex.

Vector field $G$ has $5$ equilibria, the origin $0$ which is a source, the canonical basis vectors $e_1,e_2,e_3$ which are saddle points, and the interior equilibrium $$x_*  = \frac{e_1 + e_2 + e_3}{ \alpha + \beta + 1}.$$
The diagonal $$D = \{x \in \Rp^3 \: : x_1 = x_2 = x_3\}$$
 is invariant and  on $D$ every nonzero trajectory converges to $x_*.$

On the face $\Rp \times \Rp^* \times \{0\}$ every trajectory  converges to $e_2.$ In words, {\em species $2$ beats species $1$ in absence of species $3.$} Similarly, species $1$ beats $2$ in absence of $3$ and $3$ beats $2$ in absence of $1.$
 This makes the set $$\Upsilon = W^s(e_1) \cup W^s(e_2) \cup W^s(e_3)$$ an heteroclinic cycle; where $W^s(e_i)$ stands for the stable manifold of $e_i.$

 The global dynamics of $G$ can now be described :
 \bit
 \item If $\alpha + \beta < 2$, $x_*$ is a sink and every trajectory starting from  $\Rp^3 \setminus  \partial \Rp^3$ converges to $x_*.$

 \item If $\alpha + \beta > 2,$ $x_*$ is a saddle whose stable manifold is $D \setminus \{0\}$ and every trajectory starting from $\Rp^3 \setminus  (\partial \Rp^3 \cup D)$ has $\Upsilon$  as omega limit cycle.
\item If $\alpha + \beta = 2,$ $\Delta$  is invariant and attracts every nonzero trajectory. In this case $\Upsilon = \partial  \Delta$ and on $\Delta \setminus (\{x_*\} \cup \Upsilon)$ every trajectory is periodic.
 \eit
All these properties are proved in Section 5.5 of \cite{HS98}.

 When $\alpha + \beta \neq 2,$ $\Delta$ is no longer invariant but general results on competitive systems first developed by Hirsch  (see in particular \cite{Hirsch88}, Theorem 1.7 or  Hirsch and Smith \cite{HirschSmith}, Theorem 3.18) imply that
 \bit
 \item There exists a compact invariant set $\Sigma,$ unordered, homeomorphic to $\Delta$ by  radial projection $x \mapsto \frac{x}{N(x)},$ such that
  $\Phi_t(x) \rar \Sigma$ as $t \rar \infty,$ for all $x \neq 0.$
Also,
 $\Sigma \cap \partial \Rp^3 = \Upsilon$ and $\Sigma \cap (\Rp^3 \setminus \partial \Rp^3)$ is a Lipschitz manifold.
     \eit
     Unordered means that if $x,y \in \Sigma$ and $y-x \in \Rp^3,$ then $x = y.$
The set $\Sigma$ is called the {\em carrying simplex}, a term coined by M. Zeeman \cite{Zeeman02},  and can be characterized as the boundary (in $\Rp^3$) of the basin of repulsion of  $\infty$ or equivalently, the boundary of the basin of repulsion of the origin. That is
$$\Sigma = \partial_{\Rp^3} R(\infty) = \partial_{\Rp^3} R(0)$$ where
$$R(\infty) = \{ x \in \Rp^3 \: : \lim_{t \rar - \infty}  \|\Phi_{t}(x)\| = \infty\} \mbox{ and } R(0) = \{x \in \Rp^3 \: : \lim_{t \rar -\infty} \Phi_t(x) = 0\}.$$ Smoothness properties of $\Sigma$ have been investigated in several papers (see in particular Mierczynski \cite{Mier94, Mier99}). Further properties of the carrying simplex for Lotka Volterra systems are discussed in Zeeman \cite{Zeeman02}.

Clearly
\beq
\label{eq:evolN}
\dot{N} = r(x)(N - (x_1^2 + x_2^2 + x_3^2) - (\alpha + \beta) (x_1x_2 + x_1 x_3 + x_2 x_3))
\eeq
from which it follows  that  for $\alpha + \beta < 2$ (respectively $\alpha + \beta > 2$)
 $\dot{N} > N (1-N) (\mbox{ respectively }  \dot{N} < N (1-N))$ whenever   $x_1 x_2 + x_2 x_3 + x_1 x_3 \neq 0.$ As a consequence,
\bit
\item  If $\alpha + \beta < 2,$ $\Sigma$ is {\em above} $\Delta.$  That is $N(x) > 1$ for all $x \in \Sigma \setminus \{e_1, e_2, e_3\};$
\item If  $\alpha + \beta > 2,$ $\Sigma$ is {\em below} $\Delta.$  That is $N(x) < 1$ for all $x \in \Sigma \setminus \{e_1, e_2, e_3\}.$
\eit
\paragraph{Random switching}  Let  $(\alpha_1,\beta_1)$ and $(\alpha_2,\beta_2)$  be  two environments - as defined by equation (\ref{def:enviML}) -  such that
$$\alpha_1 + \beta_1 > 2 \mbox{ and } \alpha_2 + \beta_2 < 2.$$  For each $j$ we let  $G^j$ denote  the vector field defined like $G$  in environment  $(\alpha^j,\beta^j),$ and  $\Phi^j, x_*^j, \Upsilon^j, \Sigma^j, etc.$
the corresponding flow,  interior equilibrium, heteroclinic cycle, carrying simplex, etc.

In view of the preceding discussion, for each $x \in \Delta$ the line $\RR^+ x$ meets $\Sigma^j$  in a single point $\varsigma_j(x).$ If $x \not \in \{e_1,e_2,e_3\}$
$N(\varsigma_1(x)) < 1 < N(\varsigma_2(x))$ while for $x \in \{e_1,e_2,e_3\}, \varsigma_i(x) = x.$
Define the  {\em cell bordered by $\Sigma^1, \Sigma^2$} as
$$\mathbf{C}(\Sigma^1, \Sigma^2) = \{ t \varsigma_1(x) + (1-t) \varsigma_2(x) \: : 0 \leq t \leq 1, x \in \Delta\}.$$
This set is homeomorphic to the closed unit ball in $\RR^3$ and  its boundary (in $\Rp^3$) is the union of the carrying simplices $\Sigma^1$ and $\Sigma^2.$ It will characterize the support of the persistent measure (when there is such a measure).

Fix $0 < \eta < \frac{\alpha_2 + \beta_2}{6}$ and set $$B = \{x \in \Rp^3 \: : 3 \eta \leq N(x) \leq 3 \}.$$ By equation (\ref{eq:evolN}), $B$ is positively invariant by $G^1$ and $G^2.$

Consider now the Markov process $X_t = (x(t),J(t)) \in M = B \times \{1, 2\}$ induced by (\ref{eq:PDMP}), where
the rate matrix is given as
\beq
\label{ratematrixML}
a = \left ( \begin{array}{ll}
0 & \tau (1-p) \\
\tau p & 0 \\
\end{array} \right )
\eeq
with $0 < p  < 1$ and $\tau > 0.$ In other words, there is a Poisson clock with parameter $\tau$ and each time the clock rings, the process switches from its current environment to the other with probability  $p$ (respectively $1-p$) if the current environment is   $2$  (respectively $1$).

Set $$M_0^{bd} = \{(x,j) \in M:  \: x_1 x_2 x_3 = 0\}, M_0^{D} = \{(x,j) \in M \: : x \in D\},$$  $$M_0 = M_0^{bd} \cup M_0^{D},$$
and $$M_+ = M \setminus M_0.$$
We shall prove here the following result.
\bthm
\label{th:Maypersist}
Assume that  $$p (\alpha_1 + \beta_1 -2) + (1-p) (\alpha_2 + \beta_2 -2) < 0$$ and $$ p r(x_*^1) \frac{\alpha_1 + \beta_1 -2}{\alpha_1 + \beta_1 + 1}
+ (1-p) r(x_*^2) \frac{\alpha_2 + \beta_2 -2}{\alpha_2 + \beta_2 + 1} > 0.$$ Then for $\tau$ sufficiently small,
 there is a unique persistent measure $\Pi.$ Moreover
 \bdes
 \iti $\Pi$ is
 absolutely continuous with respect to the  Lebesgue measure $dx_1 dx_2 dx_3 \otimes d(\delta_1 + \delta_2);$
  \itii $Supp(\Pi) =  \mathbf{C}(\Sigma^1, \Sigma^2) \times \{1,2\};$
  \itiii For all $(x,i) \in M^+,$ $$\Pi_t \Rightarrow \Pi$$  $\mathbb{P}_{x,i}$ almost surely;
   \itiv Suppose $r$ is constant on a open set meeting $\mathbf{C}(\Sigma^1, \Sigma^2);$ Then  $$|P_t((x,i), \cdot) - \Pi| \leq Cst (1 + dist(x,\partial \Rp^3 \cup D)^{-\theta}) e^{-\lambda t}$$
where  $\theta, \lambda$ are positive constants (independent on $(x,i)$).
\edes
\ethm
\brem {\rm The assumption that $r$ is constant on a open set meeting $\mathbf{C}(\Sigma^1, \Sigma^2)$ is an ad-hoc assumption chosen to simplify  the computation of the Lie brackets involved in the verification of the strong bracket condition. We conjecture that the result holds true for any smooth $r.$}
\erem
 Figures \ref{fig:F1May},  \ref{fig:F2May}, \ref{fig:swicthMay} result from simulations  by Edouard Strickler and  illustrate Theorem \ref{th:Maypersist}. Figures \ref{fig:F1May} and  \ref{fig:F2May} picture the phase portraits of $G^1, G^2$ and Figure  \ref{fig:swicthMay} is a realization of the switching process with $\tau = 10, p = 0,5$
and the function $$r(x) = 100 (\sum_{i = 1}^3 \exp [- 200 (x_i - z_*^1)^2 ])^{1/2}.$$ with $z_*^1 = (\alpha_1 + \beta_1 + 1)^{-1}$
\begin{figure}
\centering
\includegraphics[width=16cm]{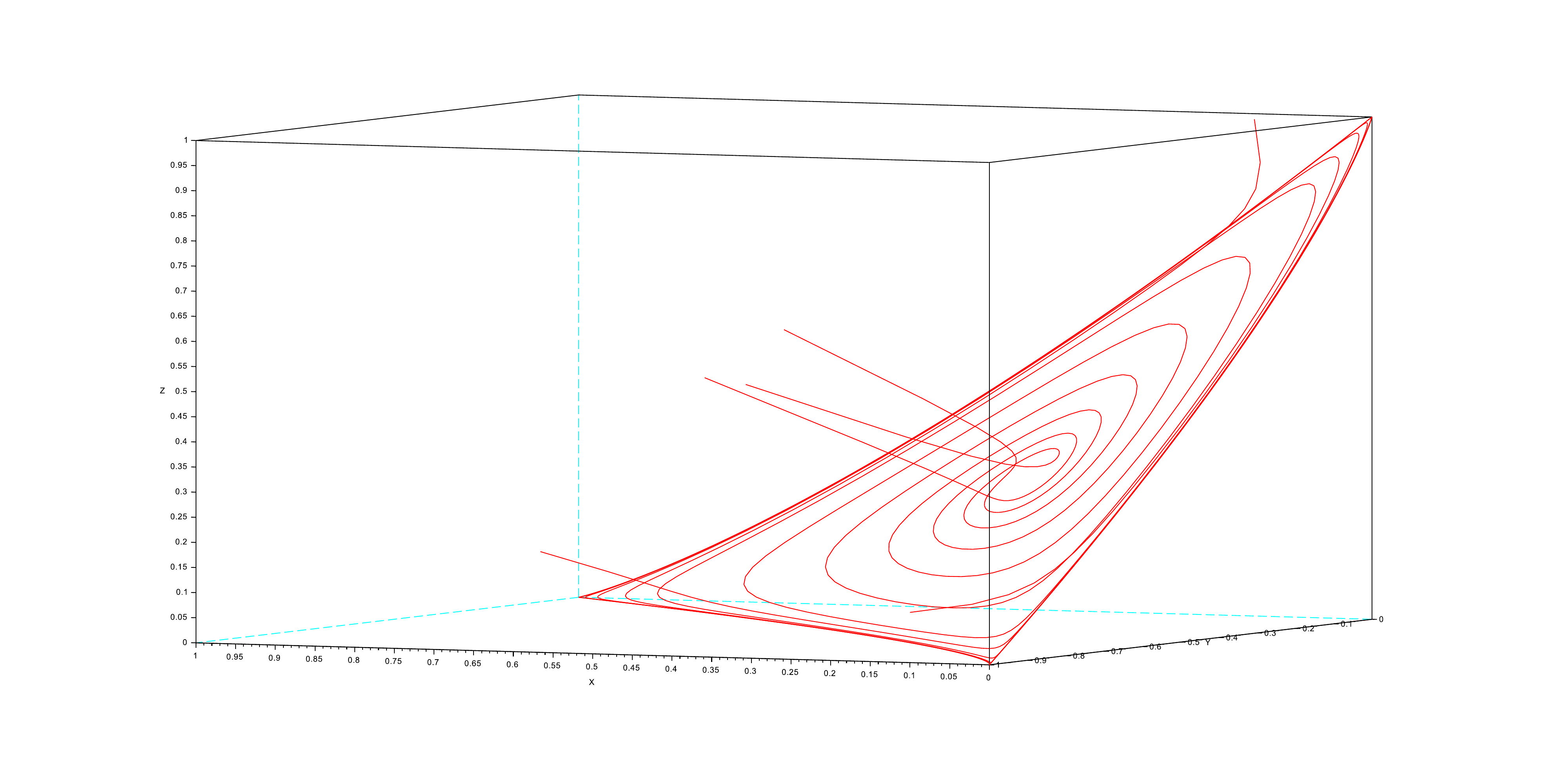}\\
\caption{$\alpha_1 = 1,8; \beta_1 = 0,6$ \label{fig:F1May}}
\end{figure}
\begin{figure}
\centering
\includegraphics[width=16cm]{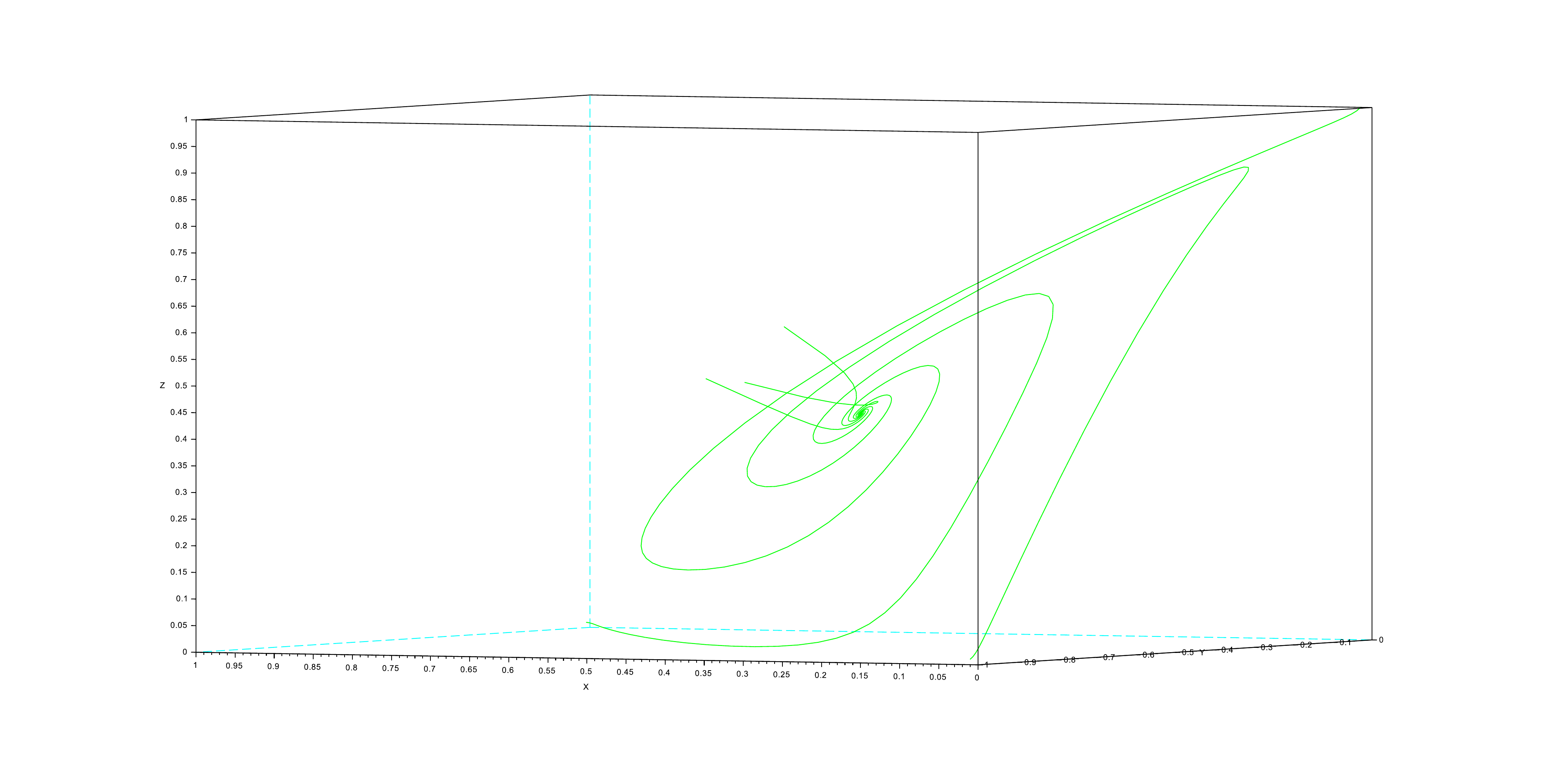}\\
  \caption{$\alpha_2 = 1,1; \beta_2 = 0,2$ \label{fig:F2May}}
\end{figure}

\begin{figure}
\centering
\includegraphics[width=16cm]{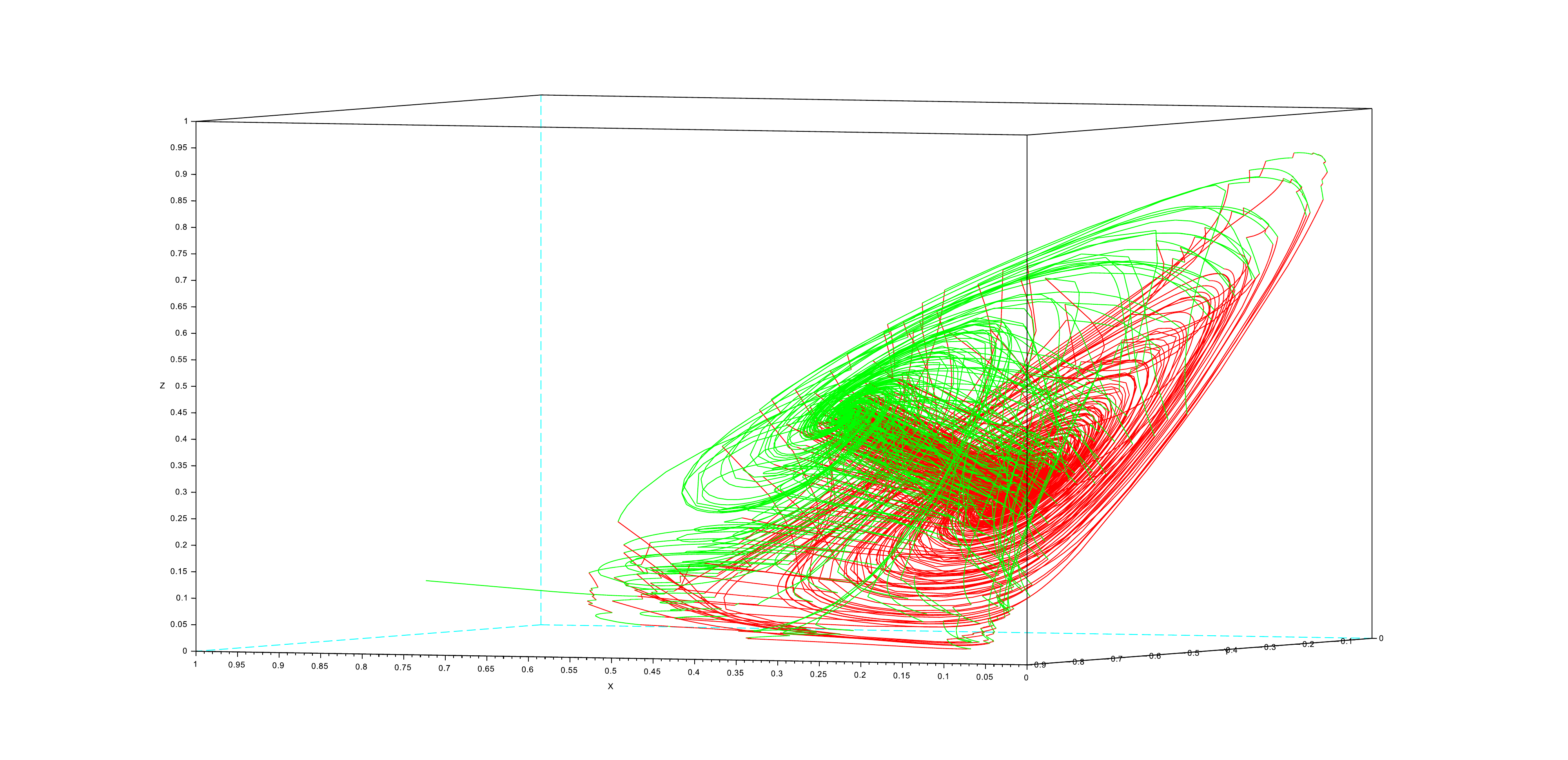}\\
 \caption{The switching process \label{fig:swicthMay}}
\end{figure}

\newpage
\paragraph{Stochastic persistence with respect to $M_0^{bd}$}
 The next proposition shows that permanence of the {\em average vector field} $p G^1 + (1-p) G^2$  implies stochastic persistence with respect to $M_0^{bd}.$
 \bprop
 \label{prop:Maybd} Let $\Lambda^{bd} = p (\alpha_1 + \beta_1 -2) + (1-p) (\alpha_2 + \beta_2 -2).$

If $\Lambda^{bd} < 0,$  then $(X_t)$ is $H$-persistent with respect to $M_0^{bd}.$
 \eprop
\prf It is not hard to prove that on the face $x_i > 0, x_{i+1} = 0,$   $x(t)$ converges to $e_i.$ A (more general) proof can be found in \cite{BL16}, Theorem 3.1. Thus the only ergodic measures on $M_0^{bd}$ are  $\delta_{e_i} \otimes \nu ,i =1,2,3$ where $\nu$ is the Bernoulli measure on $\{1,2\}$ $\nu = p \delta_1 + (1-p) \delta_2.$ The persistence criterion of Theorem \ref{persPDMP}   writes
$\sum_j \nu_j (\alpha^j + \beta^j) < 2$ and the result follows.
\qed
\paragraph{Stochastic persistence with respect to $M_0^{D}.$}
Let $\ell : \RR^3 \mapsto Ker(N) \times \RR$  be the linear change of variable defined by $l(x) = (y,z)$ with
$$y = x - z (e_1 + e_2 + e_3) \mbox{ and } z = \frac{N(x)}{3}.$$
For  $(y,z) \in  Ker(N) \times \RR$ set
 $$(G^j_1(y,z),G^j_2(y,z))= \ell \circ G^j \circ \ell^{-1}(y,z),$$
 $\rho = \|y\|$ (the Euclidean norm of $y$), and
 $\theta = \frac{y}{\rho} \in S^1$ (the unit circle in $Ker(N)$) if $\rho \neq 0$ .

In coordinates $(\rho, \theta, z) \in \Rp^* \times S^1 \times \Rp$, the dynamics of $(X_t)$ in $B \setminus D \times \{1,2\}$ rewrites,
 \beq
 \label{eq:pdmprotheta} \left \{ \begin{array}{lll}
\dot{\rho} & = & \langle \theta,  G^{J(t)}_1(\rho \theta,z) \rangle  \\
 \dot{\theta} & = &  \frac{1}{\rho} [G^{J(t)}_1(\rho \theta,z) - \langle \theta,  G^{J(t)}_1(\rho \theta,z) \rangle]. \\
\dot{z} & = & G^{J(t)}_2(\rho \theta,z)
 \end{array} \right.
 \eeq
  This extends to a dynamics on
  $$\tilde{M} = \{(\rho,\theta,z): \: \ell^{-1}(\rho\theta, z) \in B\} \times \{1,2\},$$ leaving invariant the extinction set
 $\tilde{M}_0 = ( \{0\} \times S^1 \times [\eta, 1] ) \times \{1,2\},$
 whose dynamics on $\tilde{M}_0$ is given by
 \beq
 \label{eq:pdmprothetaonM0}
 \left \{ \begin{array}{lll}
 \dot{\theta} & = &  f^{J(t)}(\theta,z) \\
\dot{z} & = & g^{J(t)}(z)
 \end{array} \right.
 \eeq
 where
 $$f^j(\theta,z) = \partial_{y} G^{j}_1(0,z) \theta - \langle \theta, \partial_{y} G^{j}_1(0,z) \theta \rangle$$
 and
 $$g^j(z) = \hat{r}(z) z (1- z (1 + \alpha_{j} + \beta_{j})).$$
 Here $\hat{r}(z)$ stands for $r(z(e_1+e_2+e_3)).$
\bprop
\label{prop:MayD}
\bdes
\iti On $\tilde{M}_0$ the process $(\theta(t), z(t), J(t))$ (given by (\ref{eq:pdmprothetaonM0})) has a unique invariant measure $\mu = \mu^1 (d\theta dz) \delta_1 + \mu^2(d \theta dz) \delta_2.$
\itii Let $$\Lambda^D = \sum_{j = 1, 2} \int \langle \partial_y G_1^j(0,z) \theta, \theta \rangle \mu^i(d\theta dz).$$  If $\Lambda^D > 0,$ then  $(\rho(t), \theta(t), z(t), J(t))$ is $H$-persistent with respect to $\tilde{M}_0.$
\itiii (slow and rapid switching).  Assume the parameter $p$ (in the definition  of the rate matrix (\ref{ratematrixML})) is fixed and write $\mu^j =\mu^j_{\tau}, \Lambda^D = \Lambda^D_{\tau}$ to emphasize the dependency on $\tau$. Then
\bdes
\ita
$$\lim_{\tau \rar 0} \Lambda^D_{\tau} = \frac{1}{2}[p r(x_*^1) \frac{\alpha_1 + \beta_1 -2}{\alpha_1 + \beta_1 + 1}
+ (1-p) r(x_*^2) \frac{\alpha_2 + \beta_2 -2}{\alpha_2 + \beta_2 + 1}] . $$
\itb $$ \lim_{\tau \rar \infty} \Lambda^D_{\tau}  = \frac{1}{2}[ r(\bar{x}_*) \frac{\bar{\alpha} + \bar{\beta} -2}{\bar{\alpha} + \bar{\beta} + 1}]$$
with $$\bar{\alpha} = p \alpha_1 + (1-p) \alpha_2, \bar{\beta} = p \beta_1 + (1-p) \beta_2$$ and $$\bar{x}_* = \frac{e_1 + e_2 + e_3}{1 + \bar{\alpha} + \bar{\beta}}.$$
\edes
\edes
\eprop
\prf $(i)$ Let $U^j(\theta,z) = (f^j(\theta,z), g^j(z)).$ We claim that for the dynamics induced by $U^1,$ every point in $S^1 \times [\eta,1]$  has  $S^1 \times \{z_*^1\} $ as $\omega$ limit set and that $U^1(\theta, z_*^1)$ and $U^2(\theta,z_*^1)$ are linearly independent. This makes  the point  $(\theta,z_*^j,1)$ an accessible point for the process $(\theta(t),\rho(t),J(t))$ at which the weak bracket condition is satisfied. The result follows from the standard arguments already used in the proof of Corollary \ref{cor:hypoode}.

We now prove the claim.
It is easy to see that
 the Jacobian matrix $DG^j(x_*^j)$ leaves $Ker(N)$ invariant and that $DG^j(x_*^j)|_{Ker(N)},$ hence $\partial_{y} G^{j}_1(0,z_*^j),$ has two non real conjugates eigenvalues $\lambda_j, \overline{\lambda}_j$ with  $$\lambda_j = \frac{r(x_*^j)}{2} ( \frac{\alpha_j + \beta_j -2}{1+ \alpha_j + \beta_j} + i \sqrt{3}(\beta_j - \alpha_j)).$$ Therefore
 $\theta \mapsto f^j(\theta,z_*^j)$ never vanishes and the first part of the claim  easily follows. For the second, note that $det(U^1(\theta, z_*^1), U^2(\theta,z_*^1)) =  f^1(\theta,z_*^1) g^2(z_*^1) \neq 0.$

The $H$-persistence follows by choosing $V(\rho,\theta,z,j) = -\log(\rho)$ (for $0 < \rho \leq  1$) and
 $$H((\rho,\theta,z,j)) = \left \{ \begin{array}{ll} \frac{1}{\rho} \langle \theta,  G^{j}_1(\rho \theta,z) \rangle  \mbox{ if } \rho \neq 0 \\ \langle \partial_y G_1^j(0,z) \theta, \theta \rangle \mbox{ if } \rho = 0
 \end{array} \right.$$

 $(ii)$ The preceding discussion implies  that $U^j$ is uniquely ergodic on $S^1 \times [\eta,1]$ with an invariant probability $\nu^j$ supported by $S^1 \times \{z_*^j\}.$
 Then $$ \int \langle \partial_y G_1^j(0,z) \theta, \theta \rangle \nu^j(d\theta dz) = \lim_{t \rar \infty} \frac{\log (\|\exp{ (t \partial_y G_1^j(0,z_*^j)) }\| )}{t} = \Re(\lambda_j).$$ On the other hand, it is not hard to show that when $p$ is fixed and $\tau \rar 0,$ every limit point of $\{\frac{\mu_{\tau}^1}{p}\}$ (respectively $\{\frac{\mu_{\tau}^2}{1-p}\}$) for the weak $*$ topology is invariant for $U^1$ (respectively $U^2$). Thus $\frac{\mu^1}{p} \Rightarrow \nu^1, \frac{\mu^2}{1-p} \Rightarrow \nu^2,$ as $\tau \rar 0$ and, consequently, $$\lim_{\tau \rar 0} \sum_{j = 1, 2} \int \langle \partial_y G_1^j(0,z) \theta, \theta \rangle \mu^j(d\theta dz)  = p  \Re(\lambda_1) + (1-p) \Re(\lambda_2).$$
 For statement $(b)$, remark that, by a standard averaging result, $\mu_{\tau} \Rightarrow \bar{\nu}$ as $\tau \rar \infty$ where $\bar{\nu}$ is the invariant probability of the average vector field $\bar{U} = p U^1 + (1-p) U^2.$ Thus, reasoning like in $(a),$
 $\sum_{j = 1, 2} \int \langle \partial_y G_1^j(0,z) \theta, \theta \rangle \mu_{\tau}^j(d\theta dz)$ converge, as $\tau \rar \infty,$
   to the real part of the conjugate eigenvalues of $D\bar{G}(\bar{x}_*)|_{Ker(N)};$ Where $\bar{G} = p G^1 + (1-p) G^2$ and $\bar{x}_* = \frac{e_1 + e_2 + e_3}{1 + \bar{\alpha} + \bar{\beta}}$ is the interior equilibrium of $\bar{G}.$
\qed
\paragraph{The accessible set}
We  now characterize the accessible set $\Gam_{M^+}.$
\bprop
\label{prop:Maysupport}
$\Gam_{M^+} = \mathbf{C}(\Sigma^1, \Sigma^2) \times \{1,2\}.$
\eprop
\prf  Relying on Proposition \ref{prop:controlode}, we  say that a point $p \in \Rp^3$ is $(G^i)$-accessible from $x \in \Rp^3$ if for  every neighborhood $O$ of $p$  there exists a control $u$ such that
     the solution $y(u,x,\cdot)$ to the control system (\ref{eq:controlode}) meets $O.$ By Proposition \ref{prop:controlode}, what we need to prove is that the set of points that are $(G^i)$ accessible from any $x \in \Rp^3 \setminus (\partial \Rp^3 \cup D)$ coincide with $\mathbf{C}(\Sigma^1, \Sigma^2).$

We first show that every $p \in \Delta,$ is $(G^i)$-accessible from every $x \in \Rp^3 \setminus (D \cup \partial \Rp^3).$
We can always assume that $p \in \Delta \setminus (\partial \Delta \cup \{\frac{e_1 + e_2 + e_3}{3}\})$ because this latter  set is dense in $\Delta.$

      Let $0 < s < 1$ be such that
 $s (\alpha_1 + \beta_1) + (1-s) (\alpha_2 + \beta_2) = 2$ and let $\bar{G} = s G^1 + (1-s) G^2.$
Note that $\bar{G}$ is the vector field defined by (\ref{eq:may}) in  the environment $(\bar{\alpha}, \bar{\beta}) = s (\alpha_1, \beta_1) + (1-s) (\alpha_2, \beta_2).$

  Let $W(x) = \frac{x_1 x_2 x_3}{N^3(x)}.$ A direct computation (see  \cite{HS98}, Section 5.5) shows that $W$ strictly decreases
   (respectively increases) along trajectories of $G^1$ (respectively $G^2$) in $\Rp^3 \setminus (D \cup \partial \Rp^3)$
    and is constant along trajectories of $\bar{G}.$

 If $W(x) > W(p)$ (respectively $<$) use the flow $\Phi^1,$  that is the control $u^1(t) = 1, u^2(t) = 0$ (respectively $\Phi^2$) to steer $x$ to a point $x'$ at which $W(x') = W(p).$ Then use the flow $\bar{G},$ that is the control $u^1(t) = s, u^2(t) = 1-s,$ until $y(u,x',\cdot)$ meets $O.$
Recall that $\Delta$ is a global attractor for $\bar{G}$ in $\Rp^3 \setminus \{0\}$ and that orbits on $\Delta \setminus \{\partial \Delta\}$ are periodic orbits (given as the level set of $W|_{\Delta}$).

We next show that every point $p \in \mathbf{C}(\Sigma^1, \Sigma^2)$ is $(G^i)$-accessible. Again it suffices to show that this is the case for $p \in  \mathbf{C}(\Sigma^1, \Sigma^2) \setminus (\Sigma^1 \cup \Sigma^2).$ Such a point $p$ lies in an interval
$]\varsigma^1(q), \varsigma^2(q)[ = \{ t \varsigma^1(q) + (1-t) \varsigma^2(q) \: : 0 < t < 1\}$ where $q \in \Delta \setminus \{e_1,e_2,e_3\}.$
If $N(p) = 1$ $p \in \Delta$ and there is nothing to prove. If $N(p) > 1,$ the characterization  $\Sigma^2 = \partial_{\Rp^3} R^2(0)$ implies that $\lim_{t \rar \infty} \Phi^2_{-t}(p) = 0.$ Therefore $\Phi^2_{-t}(p) \in \Delta$ for some $t > 0.$ Point $\Phi^2_{-t}(p)$ is then $G^i$ accessible from $x$ and so is $p$ since $p = \Phi^2_t (\Phi^2_{-t}(p)).$ If $N(p) < 1$ the proof is similar, using the  characterization  $\Sigma^1 = \partial_{\Rp^3} R^2(\infty).$ \qed
\paragraph{Proof of Theorem \ref{th:Maypersist}}

The following result  implies Theorem \ref{th:Maypersist}. We use the notation of Propositions \ref{prop:Maybd},  \ref{prop:MayD}, \ref{prop:Maysupport}.
\bthm
\bdes
\iti If $\Lambda^{bd} > 0$ and $\Lambda^D > 0,$  there exists a unique persistent measure $\Pi$ verifying the conclusions of $(i),(ii),(iii)$ of  Theorem \ref{th:Maypersist} and conclusion $(iv)$ for generic $r.$
\itii If $\Lambda^{bd} < 0$ and $\Lambda^D > 0,$ $x(t) \rar \partial \Rp^3$ almost surely, for all $x(0) \in B \setminus D;$
\itiii If $\Lambda^{bd} > 0$ and $\Lambda^D < 0,$ $x(t) \rar D$  almost surely, for all $x(0) \in B \setminus \partial \Rp^3;$
\itiv If $\Lambda^{bd} < 0$ and $\Lambda^D < 0.$  $x(t) \rar \partial \Rp^3 \cup D$ almost surely  for all $x(0) \in B \setminus (\partial \Rp^3 U D)$ and both events
$x(t) \rar  \partial \Rp^3$ and $x(t) \rar  D$ have positive probability.
\edes
\ethm
\prf
We only prove  the first assertion.
The other ones are a consequence of the extinction results to be described in part II. They can also be proved directly like Theorems 3.1, 3.3 and 3.4 in \cite{BL16}.

In view of Propositions \ref{prop:Maybd}, \ref{prop:MayD}, \ref{prop:Maysupport} and Theorem \ref{th:expoconvcompact}, it suffices to show that there exists a point $x \in \mathbf{C}(\Sigma^1, \Sigma^2)$ at which the weak (respectively strong)  Hörmander condition is satisfied.
Let
$$C^i = \begin{bmatrix}
          -1 & -\alpha_i & -\beta_i \\
          -\beta_i & -1 & -\alpha_i \\
          -\alpha_i & -\beta_i & -1 \\
        \end{bmatrix}.$$
For $x \in \Rp^3$ let $\D(x)$ denote the diagonal matrix whose entries are the components of $x$ and let ${\bf 1} = e_1 + e_2 + e_3.$
Then $G^i(x) = r(x) U^i(x)$ with $U^i(x) = \D(x) ({\bf 1} + C^i x).$ Since the  term $r(x)$ has no incidence on the weak bracket condition, it suffices to verify that it holds for the vector fields $U^1,U^2.$
A straightforward computation show that
$$[U^2,U^1](x) = \D(x) (C^1 \D(x) C^2 x - C^2 \D(x) C^1 x) + U^1(x) - U^2(x).$$
Thus
$Det(U^1(x),U^2(x),[U^2,U^1](x)) = (x_1 x_2 x_3)^3 P(x)$ where
$$P(x) = Det ({\bf 1} + C^1 x, {\bf 1} + C^2 x, C^1 \D(x) C^2 x - C^2 \D(x) C^1 x).$$
Since the function $P$ is a polynomial in the variables $x_1,x_2,x_3,$ it suffices  to show that is is not identically $0$ to deduce that $P(x) \neq 0$ for some $x$ in the interior of $\mathbf{C}(\Sigma^1,\Sigma^2).$ The tedious computation of the coefficients of $P$ becomes a child's play with the help of the mathematical  software Python/Sympy and the great help of Jean Baptiste Bardet who knows how to use it. It appears that the coefficient of the monomial $x_1 x_2 x_3$ is $$P_{1,1,1} =3 \left [(\alpha_1 + \beta_1 -2)(\beta_2 -1) - (\beta_1 -1)(\alpha_2 - \beta_2 -2) \right ] (\alpha_1 + \beta_1 - (\alpha_2 + \beta_2))$$ which is never $0.$ This concludes the proof of  the weak bracket condition.

For the strong bracket condition, under our assumption that $r$ is constant on a open set meeting $\mathbf{C}(\Sigma^1, \Sigma^2)$, it suffices to show that   $$Q(x) = \frac{1}{(x_1 x_2 x_3)^3} Det(U^1(x) - U^2(x),[U^2,U^1](x), [[U^2,U^1],U^1](x))$$
is  a non zero polynomial.
 Thanks again to Python/Sympy and Jean Baptiste Bardet,
the coefficient of $x_1^4 x_2^2$ in $Q$ is
 $$Q_{4,2,0} =  \alpha_1^2\beta_2^2 - 2\alpha_1^2\beta_2 - 2\alpha_1\alpha_2\beta_1\beta_2 + 2\alpha_1\alpha_2\beta_1 +
2\alpha_1\alpha_2\beta_2 - \alpha_1\beta_1\beta_2$$
$$ + 2\alpha_1\beta_1 + \alpha_1\beta_2^2 - 2\alpha_1\beta_2 + \alpha_2^2\beta_1^2 - 2\alpha_2^2\beta_1 + \alpha_2\beta_1^2 - \alpha_2\beta_1\beta_2 - 2\alpha_2\beta_1 + 2\alpha_2\beta_2.$$
 The coefficient of $x_1^2 x_2^3$ is
 $$Q_{2,3,0} = -2 \alpha_1^3 \beta_2 + 2 \alpha_1^2 \alpha_2 \beta_1 + 2 \alpha_1^2 \alpha_2 \beta_2 - 2 \alpha_1^2 \beta_1 + 2 \alpha_1^2 \beta_2 - 2 \alpha_1 \alpha_2^2 \beta_1$$ $$+ 2 \alpha_1 \alpha_2 \beta_1 - 2 \alpha_1 \alpha_2 \beta_2
  + 2 \alpha_1 \beta_1^2 \beta_2 - 2 \alpha_1 \beta_1^2 - 2 \alpha_1 \beta_1 \beta_2^2 + 2 \alpha_1 \beta_1 \beta_2
   - 4 \alpha_1 \beta_1 + 4 \alpha_1 \beta_2 - 2 \alpha_2 \beta_1^3 + 2 \alpha_2 \beta_1^2 \beta_2$$
    $$+ 2 \alpha_2 \beta_1^2 - 2 \alpha_2 \beta_1 \beta_2 + 4 \alpha_2 \beta_1 - 4 \alpha_2 \beta_2$$
  The  solutions  of the polynomial equation $Q_{4,2,0} = Q_{2,3,0} = 0$ are the sets
   $$\{\alpha_1= -1, \beta_1= 0\},
 \{\alpha_1= 0, \alpha_2= 0\},
 \{\alpha_1= 0, \beta_1= 2\},
 \{\alpha_1= -\beta_1 - 1, \alpha_2= -\beta_2 - 1\},$$
$$ \{\alpha_1= -\beta_1 + 2, \alpha_2= -\beta_2 + 2\},
 \{\alpha_1= \beta_1/2 - 1, \alpha_2= \beta_2/2 - 1\},
 \{\beta_1= 0, \beta_2= 0\}.$$
 None of these solutions is compatible with the constraints on the parameters. Hence $Q$ is non zero and the strong bracket condition holds true.
   \qed
   \section{$H$-Exponents and Lyapunov Exponents}
   This section discusses the relations between $H$-exponents as defined in Section \ref{sec:persistence} and classical Lyapunov exponents. For this purpose we consider the  situation where the process $(X_t)$ is  solution to the SDE (\ref{eq:sde}) on $\RR^n$ under the assumptions that  $\alpha_i = 0$ for all $i,$ and
   $$F(0) = \Sigma^1(0) = \ldots = \Sigma^m(0) = 0.$$ Such a situation has been considered by  Baxendale in   \cite{Bax90}. We will   retrieve here and (mildly extend)  some of his results. Note that similar results for ODEs with random switching have been recently obtained in \cite{BStr17}.

   In addition to the assumption that $0$ is a common equilibrium of the vector fields $F$ and $\Sigma^j$ we  assume that $F$ and $\Sigma^j$ are smooth, $\Sigma^j$ are bounded and that the conditions (\ref{lyapsde}, \ref{gammasde}) of Proposition \ref{prop:ecosde} are satisfied. For all $i = 0, \ldots, m$ we let
   $$A^i = DS^i(0)$$ denote the Jacobian matrix of $S^i$ at the origin.

   Recall (see Section \ref{sec:degensde}) that the SDE (\ref{eq:sde}) can be written as the Stratonovich SDE on $\RR^n$
   $$dx_t = F(x_t) dt + \sum_{j = 1}^m \Sigma^j(x_t) dB_t^j = S^0(x_t) dt + \sum_{j = 1}^m S^j(x_t) \circ dB_t^j$$ where $S^j = \Sigma^j$ for $j \geq 1$ and $$S^0(x) = F(x) - \frac{1}{2} \sum_{j = 1}^m DS^j(x) S^j(x).$$ Clearly  $\{0\}$ is invariant under this SDE, but the dynamics on $\{0\}$
    is  trivial and doesn't  convey any  information  on the behavior of the process away from $\{0\}.$ To circumvent this problem a useful trick, (whose idea goes back to Hasminskii for linear stochastic differential equations \cite{Has67}) is to replace the origin by the unit sphere by working in polar coordinates.

 Set $M_+ := \Rp^* \times S^{n-1}, M_0 = \{0\} \times S^{n-1}$ and $M = M_+ \cup M_0.$

 Let
 $\mathrm{P} : \RR^n \setminus \{0\} \mapsto M_+$ be the polar decomposition diffeomorphism defined as
$$\mathrm{P}(x) = (\|x\|, \frac{x}{\|x\|})$$  and let $\tilde{S}^j$  be the vector field  on $M_+$ defined as the pushforward of $S^i$ by $\mathrm{P}.$  Note that for all $(\rho, \theta) \in \Rp^* \times S^{n-1}$ and $u \in \RR^n$
$$D \mathrm{P}(\rho \theta). u = (\langle \theta, u \rangle, \frac{u - \langle \theta, u \rangle \theta}{\rho}).$$ Thus
$$\tilde{S}^j(\rho,\theta) = D\mathrm{P}(\rho \theta) S^j(\rho \theta) = (\langle \theta, S^j(\rho \theta) \rangle,
\frac{1}{\rho}(S^j(\rho \theta) - \langle S^j(\rho \theta), \theta \rangle \theta)).$$
Observe that $\tilde{S}^j$ extends smoothly to $M$   by
setting
$$\tilde{S}^j(0,\theta) = (0, \tilde{A^j}(\theta))$$
where $\tilde{A^j}$ is the vector field on $S^{n-1}$ defined as
\beq
\label{eq: deftildeA}
\tilde{A^j}(\theta) = A^j \theta - \langle A^j \theta, \theta \rangle \theta.\eeq The process $\tilde{X}_t:= (\rho_t, \theta_t) = \mathrm{P}(X_t)$ can then be viewed as the restriction to $M_+$ of the process on $M$ solution to the SDE
\beq
\label{eq:stratopolar}
d(\rho,\theta) = \tilde{S}^0(\rho,\theta) dt + \sum_j \tilde{S}^j(\rho,\theta) \circ dB_t^j.
\eeq
On $M_0$ naturally identified with $S^{n-1}$ the dynamics writes
\beq
\label{eq:stratopolartheta}
d \theta_t = \tilde{A^0}(\theta) dt + \sum_{j = 1}^m \tilde{A}^j(\theta) \circ dB_t^j.
\eeq
Let $\tilde{L}$  denote the formal generator  of (\ref{eq:stratopolar}). Then, for any smooth function $f : M_+ \mapsto \RR,$
$\tilde{L}(f) = L (f \circ \mathrm{P}) \circ \mathrm{P}^{-1}.$

Let $\ln_{*}  : \Rp^* \mapsto \RR_{-}$ denote a smooth function such that $\ln_{*}(t) = \ln(t)$ for $t <1/2$ and $\ln_{*}(t) = 0$ for $t \geq 1.$

\bprop Let $V: M_+ \mapsto \Rp$ and $H : M \mapsto \RR$ be the maps defined by
  $$V(\rho,\theta) = - \ln_{*}(\rho),$$
  $$H(\rho,\theta) = \tilde{L}(V)(\rho,\theta) \mbox{ for } \rho > 0,$$ and
   \beq \label{eq:HforSDEwith02}
  H(0,\theta)
  = - \left  \{ \langle A^0 \theta, \theta \rangle  + \frac{1}{2} \sum_{k = 1}^m  \left (\langle (A^k)^2 \theta, \theta \rangle +
  \|A^k \theta\|^2  -  2 \langle A^k \theta , \theta \rangle^2 \right ) \right \}
  \eeq Then $(V,H)$   satisfy hypothesis \ref{hyp:H+}.

\eprop
\prf
Recall that for all $x \in \RR^n$
$a_{ij}(x) = \sum_{k = 1}^m S^k_i(x) S^k_j(x).$
Then, for $\rho  < 1/2$
$$H(\rho,\theta) = L (-\ln(\|.\|))(\rho \theta) = -\frac{1}{2} [\frac{L (\|.\|^2)}{\rho^2}(\rho \theta) - \frac{1}{2 \rho^4}\Gamma_L(\|.\|^2)(\rho \theta)]$$
$$= \frac{1}{\rho^2} [ - \langle  F(\rho \theta), \rho \theta\rangle -  \frac{1}{2}\sum_i a_{ii}(\rho \theta) +  \sum_{ij} a_{ij}(\rho \theta) \theta_i \theta_j ]$$
$$
= \frac{1}{\rho^2} [ - \langle  F(\rho \theta), \rho \theta\rangle -  \frac{1}{2}\sum_{k = 1}^m\|S^k(\rho \theta)\|^2 +  \sum_{k = 1}^m \langle S^k(\rho \theta) , \theta \rangle^2\ ].$$ When $\rho \rar 0$ this latter expression converges to $H(0,\theta)$ defined by (\ref{eq:HforSDEwith02}).

We proceed now like in the proof of Theorem \ref{persSDE}. Let $B :\Rp \mapsto [0,2]$ be a smooth map such that $B(t) = t$
for $t \leq 1$ and $B(t) = 2$ for $t \geq 2.$ Set $V_n(\rho,\theta)  = n B(\frac{V(\rho,\theta}{n})$ for $\rho > 0$
 and
 $V_n(0,\theta) = 2 n.$ Then $V_n \in {\cal D}^2 ({\cal {\tilde L}})$ (because $V_n$ is smooth with compact support),
 $V_n = V$ and ${\cal L}(V_n) = H$ on $\{V \leq n\}.$
 Also $\Gamma(V_n) = B'(\frac{V}{n})^2 \Gamma_{\tilde{L}}(V)$ and, for $\rho > 0$
 $$\Gamma_{\tilde{L}}(V)(\rho,\theta) = \Gamma_L (V \circ P)(\rho \theta) =
 \sum_{k \geq 1}   (D (V \circ P)(\rho \theta) S^k(\rho \theta))^2 = \sum_{k \geq 1} (DV(\rho,\theta) \tilde{S}^k(\rho, \theta))^2
 .$$ Thus
  $\Gamma_{\tilde{L}}(V)(\rho,\theta) = \sum_{k \geq 1} [\frac{1}{\rho} \langle \theta,S^k(\rho \theta) \rangle]^2$
   for $\rho < 1/2.$ This shows that $\Gamma_{\tilde{L}}(V),$ hence $(\Gamma(V_n))_{n}$ is bounded and the results follows from
  Proposition \ref{martconv}.
\qed
We now make precise the links with Lyapunov exponents.
Consider the linear SDE on $\RR^n$
\beq
\label{eq:linearsde}
dy_t = A^0 y dt + \sum_{k = 1}^m A^k y \circ dB_t^k.
\eeq
Let $t \rar \varphi(t,\omega) y$ denote the  solution to (\ref{eq:linearsde}) with initial condition $y \in \RR^n.$
By the multiplicative ergodic theorem (see e.g~\cite{Arn98} and compare to \cite{Bax90}, Proposition 2.8) there exist $1 \leq d \leq n$ numbers
$$\lambda_d < \ldots < \lambda_1$$ called the {\em Lyapunov exponents} of (\ref{eq:linearsde}), a set $\Omega_0 \subset \Omega$ of full measure, and for all $\omega \in \Omega_0$  disctinct  vector spaces
$$\{0\} = V_{d+1}(\omega) \subset V_d(\omega) \subset \ldots \subset V_1(\omega) = \RR^n$$ (measurable in $\omega$) such that
$$\lim_{t \rar \infty} \frac{1}{t} \log \|\varphi(t,\omega) y\| = \lambda_i$$ for all $y \in V_i(\omega) \setminus V_{i+1}(\omega).$

\bprop One has
$\Lambda^+(H) = \lambda_1 \mbox{ and } \Lambda_{-}(H) \in \{\lambda_d, \ldots, \lambda_1\}$
\eprop
Recall (see section \ref{sec:degensde}) that given  a family ${\cal X}$ of smooths vector fields on $\RR^n$ we let $[{\cal X}]$ denote the family consisting of  ${\cal X}$ and all the Lie brackets obtained recursively from ${\cal X}.$ This definition extends obviously to the situation where the vector fields are defined on a manifold (such as $S^{n-1}$). The next result gives sufficient conditions ensuring that both  $\Lambda_{-}(H)$ and $\Lambda_{+}(H)$ coincide with $\lambda_1.$ Condition $(i)$ is the one assumed by Baxendale (\cite{Bax90}, condition (2.5)) and condition $(ii)$ is similar (albeit more general) to condition (\cite{Bax90}, condition (5.5)).
\bprop
Consider the two following conditions:
\bdes
\iti For all $x \in \RR^n \setminus \{0\} \, [\{A^1, \ldots, A^m\}](x)$ spans $\RR^n$
\itii There exists $\theta_0 \in S^{n-1}$ such that
\bdes
 \ita $[\{\tilde{A^0}, \ldots, \tilde{A^m}\}](\theta_0)$ has dimension $n-1$
\itb For every $\theta \in S^{n-1}$ and every neighborhood $O$ of $(\theta_0,-\theta_0)$ in $S^{n-1}$ there exists a (piecewize continuous) control $u = (u^1, \ldots, u^m)$ such that the solution to
$$\dot{\alpha} = \tilde{A^0}(\alpha) + \sum_{k = 1}^m \tilde{A^k}(\alpha) u^k(t)$$ with initial condition $\alpha(0) = \theta$ meets $O.$
\edes
\edes
Then $$(i) \Rightarrow (ii) \Rightarrow \lambda_+(H) = \Lambda_-(H) = \lambda_1.$$
\eprop
\prf
For $x  \in  \RR^n \setminus \{0\}$ let $\mathrm{P}_2(x) = \frac{x}{\|x\|}.$ If $A$ is a linear vector field on $\RR^n$, let $\tilde{A}$ denote the vector field on $S^{n-1}$ defined by (\ref{eq: deftildeA}) with $A$ instead of  $A^j.$ It is easy to check that
for every smooth map $f : S^{n-1} \mapsto \RR,$ $A(f \circ P_2) = \tilde{A}(f) \circ P_2.$ Thus $[A,B](f \circ P_2) = [\tilde{A},\tilde{B}](f) \circ P_2,$ and consequently, $Y(f \circ P_2) = \tilde{Y}(f) \circ P_2$ for all $Y \in [\{A^1, \ldots, A^m\}].$ This implies that for all $\theta_0 \in S^{n-1}$ $DP_2(\theta_0) [\{A^1, \ldots, A^m\}] = [\{\tilde{A^1}, \ldots, \tilde{A^m}\}](\theta_0).$
Thus $[\{\tilde{A^1}, \ldots \tilde{A^m}\}](\theta_0)$ has  dimension $n-1$ for all $\theta_0 \in S^{n-1}$ because $DP_2(\theta_0)$ has rank $n-1.$ This obviously implies condition $(ii), (a)$,  and condition $(ii), (b)$ by Chow's theorem.

 It remains to prove that under condition $(ii)$ $\Lambda_+(H) = \Lambda_{-}(H).$  Let $P\RR^{n-1}$ be the projective space (the quotient of $S^{n-1}$ by the equivalence relation identifying antipodal points) and $p : S^{n-1} \mapsto P\RR^{n-1}$ the projection map. If $(\theta_t)$ is solution to (\ref{eq:stratopolartheta}) on $M_0$ so is $(-\theta_t)$. This makes $(p(\theta_t))$  a Feller Markov process on $P\RR^{n-1}$ for which  $p(\theta_0)$ is  accessible  (by condition $(ii), (a)$); and  satisfies the Hörmander condition  (by condition $(ii), (b)$)). Therefore (see Corollary \ref{cor:hyposde} (i)), $(p(\theta_t)$ has a unique invariant probability measure $\nu.$
Let now $\mu$ be any invariant probability measure for $(\theta_t).$ Then  $\mu \circ p^{-1}$ (the image measure of $\mu$ by $p$) is invariant for  $(p(\theta_t)),$ hence $\mu \circ p^{-1}= \nu$ and $\mu H$ = $\nu h$ where $h$ denotes the real valued map on $P\RR^{n-1}$  defined by $h (p(\theta)) = H(0,\theta).$
\qed

Part $(i)$ of the next result is similar to Theorem 2.13 of Baxendale \cite{Bax90}. Part $(ii)$ gives some exponential rate of convergence
\bthm
Assume $0$ is an accessible point for the SDE (\ref{eq:sde}) (or equivalently $M_0 \subset \Gam_{M_+}$
for the process $((\rho_t,\theta_t))$ solution to (\ref{eq:stratopolar})),  $[\{A^1, \ldots, A^m\}](x)$ spans $\RR^n$ for all $x \neq 0$ and $\lambda_1 > 0.$
 Then \bdes
  \iti ${\cal P}_{inv}(\RR^n \setminus \{0\}) = \{\Pi\}$ where $\Pi << \lambda,$
   \itii $\lim_{t \rar \infty} \Pi_t^x f = \Pi f$ a.s
 for all $x \in \RR^n \setminus \{0\}$ and  $f \in L^1(\Pi);$
 \itiii
$$|P_t f(x) - \Pi f| \leq cst (1 + W_{\theta}(x)) e^{-\lambda t} \|f\|_{W_{\theta}}$$
for some positive constants $\lambda, \theta$ and $W_{\theta}(x) = U(x) (1 + \frac{1}{\|x\|^{\theta}}).$
\edes
\ethm
\prf
We claim that for some $r > 0$ and $0 < \|x\| < r$ $[\{F^1, \ldots, F^m\}](x)$ spans $\RR^n.$ Thus,
by Chow's theorem for every pair of point $x_0, x_1 \in B(0,r) \setminus \{0\}$ there is a solution $y(.)$ to the deterministic control system \ref{controlsde} with $y(0) = x_0$ and $y(t) = x_1$ for some $t \geq 0.$ Since  $0$ is accessible, it follows that all the points in  $B(0,r) \setminus \{0\}$ are  accessible Doeblin points and the result follows from Corollary \ref{cor:hyposde} $(ii).$

We now prove the claim. By the assumption that $[\{A^1, \ldots, A^m\}](x)$ spans $\RR^n$ for all $x \neq 0$ there exists a finite covering of $S^{n-1}$ by open set $(O_j)_{j \in J}, \eps >0$ and for each $j \in J$ vector fields $\{G^{1,j}, \ldots, G^{n,j}\} \subset [\{[F^1, \ldots, F^m\}]$ such that for all $\theta \in O_j$
$$|det(DG^{1,j}(0) \theta, \ldots, DG^{n,j}(0) \theta)| > \eps.$$
Since
$$\|G^{i,j}(x) - DG^{i,j}(0).x\| \leq O(\|x\|^2),$$ one can choose  $r$ small enough such that for  $x \in B(0,r) \setminus \{0\}$ and $\frac{x}{\|x\|} \in O_j$
$$|det(G^{1,j}(x) \theta, \ldots, G^{n,j}(x) \theta)| > \eps \|x\|^n > 0.$$ This proves the claim. \qed
\section{Return times and convergence rates}
\label{sec:rate}
Theorem \ref{posinvthm} shows that, under $H$ persistence, the process $(X_t^x)$ spends most of its time in a compact set far from the extinction set.
 Here we are interested in more quantitative  consequences of $H$ persistence. First we will estimate the mean time needed to reach such a compact set.Then we will give condition ensuring  that the rate of convergence in Theorem """ is exponential
Throughout the remainder of this section we assume,
 that the process $\{X_t^x : \: x \in M_+\}$ is $H$-persistent, that is $\Lambda^{-}(H) > 0,$
\subsection*{The case $M_0$ compact}
\label{sec:ratecompact}
We  assume here that $M_0$ is compact.
For $\delta > 0$ we let
$$M_0^{\delta} = \{x \in M_+ : \: d(x,M_0) < \delta\}.$$

\bprop \label{pakeslem}
Let $0 < \lambda < \Lambda^{-}(H).$
  For every $T_0 > 0$ (sufficiently large) and $T_1 > T_0,$   there exists  $\delta > 0$  such that for all $T \in [T_0,T_1]$
\begin{eqnarray*}
   P_T V(x) - V(x) & \leq &  - \lambda T  \mbox{ for } x \in  M_0^{\delta}.\\
\end{eqnarray*}
Given such a $T \in [T_0,T_1],$ let $$\tau_1  \defn \inf \{k \in  \NN^*  :  X^x_{kT} \in M_+ \setminus M_0^{\delta} \}$$ and
$$\tau_n = \inf \{k \in  \NN : k > \tau_{n-1}, X^x_{nT} \in M_+ \setminus M_0^{\delta}\}.$$  Then
$$
\E_x(\tau_1)  \leq \left \{
\begin{array}{c}
\displaystyle \frac{V(x)}{\lambda T} \mbox{ if } x \in M_0^{\delta}, \\
\\
\displaystyle 1 + \frac{P_T V(x) }{\lambda T} \mbox{ if } x \in M_+ \setminus  M_0^{\delta}
 \end{array} \right. $$ and

$$ \E_x(\tau_{n+1}) \leq \E_x(\tau_n) + 1 + \frac{v(\delta, T)}{\lambda T}$$ for all $n \geq 1,$
where
 $$v(\delta, T) := \sup \{ P_T V(x) : \: x \in M_+ \setminus  M_0^{\delta} \}.$$
\eprop
\prf
For all $x \in M, t \geq 0$ we let
$$\overline{H}(t,x) =  \int_0^t P_s H(x).$$
Recall, that by Lemma \ref{martconv}, $$P_T V(x) - V(x) =  \overline{H}(T,x)$$ for all $x \in M_+, T \geq 0.$
The first assertion then follows from  the two following facts (and compactness of $M_0$):
\bdes
\ita There exists $T_0 > 0$ (arbitrary large) such that all $x \in M_0$ and $t  \geq T_0$  $\overline{H}(t,x) < -\lambda t;$
\itb  $\overline{H}$ is continuous in $(t,x).$
\edes

 Proof of a). Suppose the contrary. Then for all $n \in \NN^*, \exists t_n \geq n, x_n \in M_0$ such that
$\mu_n H \geq - \lambda$ where $\mu_n$ stands for the measure defined by
$$\mu_n f = \frac{1}{t_n} \int_0^{t_n} P_s f(x_n) ds$$ for all $f \in {\cal M}_b(M).$
By Proposition \ref{tightlypnv}, $ \mu_n \tilde{W} \leq \frac{W(x_n)}{t_n} +  C.$ Thus,  by Lemma \ref{lem:tightpi}, $(\mu_n)$ is tight.  Let $\mu$ be a limit point of $(\mu_n)$ it is easily seen that $\mu \in {\cal P}_{inv}(M_0)$ (because for all $f \in C_b(M_0), |\mu_n f - \mu_n P_t f| \leq \frac{2 t \|f\|}{t_n} \rar 0$ as $n \rar \infty$).
By Hypotheses \ref{hyp:H} (ii),  \ref{hyp:tightpi}  and Lemma \ref{lem:tightpi} (ii) we get that $\mu H \geq - \lambda > -\Lambda^-(H).$ A contradiction.
\\

Proof of b). Let $(\mu_n)$ be defined like in the proof of $(a)$ but, this time with $t_n \rar t^*$ and $x_n \rar x_*.$  The sequence $(\mu_n)$ is tight,  and for every limit point $\nu$  of $(\mu_n)$ and  all $f \in C_b(M)$
$\nu f = \frac{1}{t^*} \int_0^{t*} P_s f(x_*)ds$ by $C_b(M)$ Feller continuity. Thus $\nu f = \frac{1}{t^*} \int_0^{t*}P_s f(x^*) ds$ for all $f \in L^1(\nu).$ In particular $\lim_{n \rar \infty} \overline{H}(t_n,x_n) = \overline{H}(t,x)$.

The last assertions follow from Pakes's criterion (see Theorem 9.1.2 and its proof in \cite{duf00})
\qed
\brem
\label{rem:pakes}
{\rm Although there is no evidence that the quantity $v(\delta,T)$ in Proposition \ref{pakeslem} is finite, we can always modify $V$ and $H$ outside a neighborhood of $M_0$ such that:
\bdes
\iti $H$ is bounded on $M,$ and
\itii $V$ is bounded on $M_+ \setminus M_0^{\delta}$ for all $\delta > 0;$
\edes
In particular  $$v(\delta,T) \leq \sup \{ V(x) \:  : x \in M_+ \setminus M_0^{\delta}\} + T \|H\|.$$
Indeed, let $C$ be a compact set such that   $M_0 \cup M_0^{\delta_0} \subset \mathrm{int}(C)$ for  some $\delta_0 > 0.$  The set $K = C \setminus (M_0 \cup M_0^{\delta_0})$ is a compact subset of $M_+.$ Set $\tilde{V}(x) = V(x)$ if $x \in C \setminus M_0,$ $\tilde{H}(x) = H(x)$ if $x \in C,$  $\tilde{V}(x) = V_K(x)$  and $\tilde{H}(x) = {\LA}(V_K)(x)$ for $x \not \in C.$ The map $(\tilde{V}, \tilde{H})$ coincide with $(V,H)$ on $C \setminus M_0 \times C$ and satisfies the required conditions.
}
\erem
\bprop \label{hajeklem}
Assume that the process is $H$-persistent (strong version) with $V$ and $H$  like in Remark \ref{rem:pakes}. Then,
 for every $T_0 > 0$ (sufficently large) and $T_1 > T_0,$   there exist positive numbers  $\theta, \delta, \kappa$   and $\rho < 1$ such that for all $T \in [T_0,T_1]$
$$P_T (e^{\theta V}) (x) \leq \left \{
\begin{array}{c}
 \displaystyle \rho e^{\theta V(x)}  \mbox{ on } M_0^{\delta}, \\
 \\
\displaystyle \kappa \mbox{ on } M_+ \setminus M_0^{\delta}
  \end{array} \right. $$

Furthermore, letting $b = 1/\rho,$
$$\E_x(b^{\tau}) \leq \left \{
\begin{array}{c}
\displaystyle e^{\theta V(x)} \mbox{ if } x \in M_0^{\delta}, \\
\\
\displaystyle b(1 + P_T e^{\theta V}(x))  \mbox{ if } x \in M_+ \setminus  M_0^{\delta}
 \end{array} \right. $$
and
\beq
\E_x(b^{\tau_{n+1}}) \leq \E_x(b^{\tau_{n}}) b (1 + \kappa)
\eeq
for $n \geq 1.$
\eprop
\prf
Let $\lambda, T_0,T_1, \delta$ be as in Proposition \ref{pakeslem}.
 For $x \in M_+$ and $T \in [T_0, T_1]$

\beq
\label{eq:expvtheta}
e^{\theta V(X_T)} = e^{\theta V(x)} e^{\theta (P_T V(x) - V(x))} e^{\theta (M^1_T(x) + M_T^V(x))}
\eeq
where
$$M^1_T(x) = \int_0^T H(X^x_s) - P_s H(x) ds$$ and $(M_t^V(x))$ is the Martingale defined in Lemma \ref{martconv}.

Observe that $\E (M^1_T(x)) = 0$ and $M^1_T(x)^2 \leq c_1 T^2$ with $c_1 = 4 \|H\|^2.$
Thus, by elementary properties of the log-laplace transform
\beq
\label{eq:loglapM1} \E (e^{\theta M^1_T(x)}) \leq e^{ c_1 T^2 \frac{\theta^2}{2}}.
\eeq
Let $$r(x) = e^{x} - x - 1 \leq x^2 e^x.$$ Using Lemma 26.19 of \cite{Kall},
the process \beq
\label{eq:expomartZ}
Z_t(\theta) = \exp{(\theta M_t^V(x) -  \frac{r(\theta \Delta V)}{  (\Delta V)^2} \langle M_t^V(x) \rangle)}
 \eeq
 is a supermartingale for all $\theta > 0$ (here and below we adopt the convention that $\frac{r(\theta \Delta V)}{  (\Delta V)^2} = \frac{\theta^2}{2}$ when $\Delta V = 0$).
Thus,
\beq
\label{eq:expmartin}
\E(\exp {(\theta M_T^V(x) -  \frac{r(\theta \Delta V)}{  (\Delta V)^2} \gamma T })) \leq 1
\eeq where $\gamma$ is the supremum in condition $(ii)$.
It follows, by Hölder inequality, that
$$\E( e^{\theta (M^1_T(x) + M_T^V(x))}) \leq e^{c_1 T^2 \theta^2} e^{ \gamma T \frac{r(2 \theta \Delta V)}{2   (\Delta V)^2}} \leq  e^{ c_2 \theta^2 T}$$
 with 
 $c_2 = c_1 T_1 + 2 \gamma e $
 and  $2 \theta \Delta V  \leq 1.$ This latter inequality combined with (\ref{eq:expvtheta}) and Proposition \ref{pakeslem} proves the first assertion with
$$\rho  =  e^{- \theta T_0 (\lambda - c_2 \theta)},
\kappa = e^{\theta \left ( T_1  \|H\|  + \sup \{ V(x) : \: x \in M_+ \setminus M_0^{\delta}  \}  + c_2 \theta T_1 \right )}$$ and $\theta$ small enough.
The last assertion follows by observing that
$$W_n = e^{\theta V(X_{n \wedge \tau})} b^{n \wedge \tau}$$ is a supermartingale with respect to $({\cal F}_{nT}).$
Hence, for $x \in M_0^{\delta}$ $$\E_x(b^{n \wedge \tau}) \leq \E_x (W_n) \leq W_0 = e^{\theta V(x)}$$ while for $x \not \in M_0^{\delta}$
  $$\E_x(b^{\tau}) =
   b \E_x( \Ind_{X_T \not \in M_0^{\delta}} + \E_{X_T} (b^{\tau}) \Ind_{X_T \in M_0^{\delta}}) \leq b (1 + P_T e^{\theta V}(x))$$ by the Markov property (compare to  the proof of Theorem 8.1.5(2) of \cite{duf00}) The bound for $ \E_x(b^{\tau_n})$ is obtained similarly by using the strong Markov property.
\qed

\subsection{The case $M_0$ non compact}
\label{sec:ratenoncompact}
The purpose of this section is to prove the following result, similar to Proposition \ref{hajeklem}, when $M_0$ is noncompact.  The proof is inspired by the proof of Proposition 4.1 in  \cite{HN18b}.
\bprop \label{hajeklemnoncompact} Assume that the process is  $H$-persistent (strong version') (see Hypothesis \ref{hyp:H+} and condition $(ii)$' in Hypothesis \ref{hyp:H}) and  {\em persistent at infinity}, meaning that $V$ is proper and
there exists a compact $C \subset M$ such that
$$\sup_{x \in M \setminus C} H(x) < 0.$$  Then
there exist positive numbers  $T_1 > T_0,$   $\theta,  \kappa$   and $\rho < 1$ such that for all $T \in [T_0,T_1]$
$$P_T (e^{\theta V}) (x) \leq \rho e^{\theta V(x)} + \kappa.$$
\eprop
From now on and throughout the section we assume  that the process is  $H-$ persistent (strong version) and
 without loss of generality\footnote{It suffices to multiply $V$ and $H$ by a sufficient large constant}, that  $$H(x) \leq -2$$ on $M \setminus C.$

We let $\theta_0$ denote a positive number small enough so that $$\gamma \theta_0 e^{\theta_0 \Delta V} \leq 1$$
 where $\gamma, \Delta V$ are like in Hypothesis  \ref{hyp:H}'.

\blem
\label{lem:croissmaj}
There exists $\omega_0 > 0$   such that for all $\theta \leq \theta_0, T \geq 0$ and $x \in M_+$
$$P_T e^{\theta V}(x) \leq e^{\theta V(x)} e^{\theta \omega_0 T},$$
and
$$ \E_x (e^{\theta ( V(T \wedge \tau) + T \wedge \tau) }) \leq e^{\theta V(x)}.$$
where $\tau = \inf \{t \geq 0: \: X_t^x \in C \};$
\elem
\prf
 Let
 \beq
 \label{eq:H*} H^* = \sup \{|H(x)|: \: x \in C\}.
 \eeq By persistence at infinity,  $H(x) \leq H^*$ for all $x \in M.$
From (\ref{vmartdefn})
$$P_T e^{\theta V}(x) \leq e^{\theta V(x)} e^{\theta H^* T} \E(e^{\theta M_T^V(x)}).$$
By (\ref{eq:expmartin}) $$\E(e^{\theta M_T^V(x)}) \leq e^{ \gamma T \frac{r( \theta \Delta V)}{   (\Delta V)^2}} \leq \exp{ (\gamma T \theta^2 e^{\theta \Delta V})} \leq e^{\theta T}.$$
It suffices to set $\omega_0 = H^* +  1$

From (\ref{vmartdefn}) again and the fact that $H \leq -2$ on $M \setminus C$
$$V(X^x_{t \wedge \tau}) \leq V(x) - 2(t \wedge \tau) + M_{t \wedge \tau}^V(x).$$ Thus
$$e^{\theta (V(X^x_{t \wedge \tau}) + t \wedge \tau)}  \leq e^{\theta V(x)} Z_{t \wedge \tau}(\theta) e^{(t \wedge \tau) (\gamma \theta^2 e^{\theta \Delta V} - \theta)} \leq e^{\theta V(x)} Z_{t \wedge \tau}(\theta)$$
where $(Z_t(\theta))$ is the supermartingale given by (\ref{eq:expomartZ}). This proves the second assertion.
\qed
\blem
\label{eq:loglaplaceweak}   Let $$M^1_T(x) = \int_0^T (H(X_s^x) - P_s H(x)) ds.$$
For all $\eps > 0$ there exists $c > 0$ such that  for all $x \in C, 0 \leq \theta \leq \theta_0$ and $T \geq 1$
$$\E(e^{\theta M^1_T(x) }) \leq  e^{ \theta T (\eps + c \theta T)}.$$
\elem
\prf follows from the two following claims.

{\em Claim 1}: Let ${\cal M}$ be a uniformly integrable family of random variables, centered (i.e $\E(M) = 0$ for all $M \in {\cal M}$) and bounded from above (i.e $M \leq c_0 < \infty $ for  all $M \in {\cal M}$).   Then for every $\eps > 0$ there exists $c > 0$ such that for all $\theta \geq 0$
 $$\E(e^{\theta M}) \leq e^{ \theta (\eps + c \theta)}.$$
{\em Proof of Claim 1:}
Write $\E(e^{\theta M}) = 1 + \E(r(\theta M))$ with $r(u) = e^u - u - 1.$ It is easily checked that $0 \leq r(u) \leq - u$ for $u \leq 0$  and  $0 \leq r(u) \leq \max(a^2, b^2 e^b)$ for $- a \leq u \leq b, a \geq 0, b \geq 0.$ Thus, for all $R > 0,$
$$\E(r(\theta M)) = \E(r(\theta M) \Ind_{M \leq -R}) + \E(r(\theta M) \Ind_{M \geq -R})$$
$$\leq - \theta \E(M  \Ind_{M \leq -R}) + \theta^2 \max (R^2, c_0^2 e^{\theta_0 c_0}))$$
By uniform integrability  choose $R$ large enough so that $\E(|M|  \Ind_{M \leq -R}) \leq \eps$ and set $c = \max (R^2, c_0^2 e^{\theta_0 c_0})).$ Then $\E(e^{\theta M}) \leq 1+  \eps \theta + c  \theta^2 \leq e^{\eps \theta + c  \theta^2}.$

{\em Claim 2}: The family ${\cal M} = \{\frac{M^1_T(x)}{T}  \: : x \in C, T \geq 1\}$ is uniformly integrable, centered and bounded from above.

{\em Proof of Claim 2:}  Set, for $p \geq 1,$  $\|H\|_{T,p}(x) : = (\frac{1}{T}\int_0^T P_s |H|^p(x) ds)^{1/p}.$ Then, $$\frac{M_1^T(x)}{T} \leq H^* +  \|H\|_{T,1}(x), \;  (\E(|\frac{M^1_T(x)}{T}|^q))^{1/q} \leq \|H\|_{T,q}(x) + \|H\|_{T,1}(x)$$ and,  by  Hypothesis \ref{hyp:H} (ii)'   and Theorem \ref{tightlypnv},
$$\|H\|_{T,1}(x) \leq \|H\|_{T,q}(x) \leq [cst (1 + W(x)/T)]^{1/q}.$$ This latter quantity being bounded for $x \in C, T \geq 1$ this proves the claim.

 \qed

\blem
\label{hajeklemnoncompact0}
 For every $T_0 > 0$ (sufficiently large) and $T_1 > T_0,$   there exist positive numbers  $\theta \leq \theta_0, \delta, \kappa$   and $\rho < 1$ such that for all $T \in [T_0,T_1]$
$$P_T (e^{\theta V}) (x) \leq \left \{
\begin{array}{c}
 \displaystyle \rho e^{\theta V(x)}  \mbox{ on } M_0^{\delta} \cap C
 , \\
 \\
\displaystyle \kappa \mbox{ on } (M_+ \cap C) \setminus M_0^{\delta}
  \end{array} \right. $$
\elem
\prf
By compactness of  $M_0 \cap C$, the first assertion of Proposition \ref{pakeslem} remains valid  if  $M_0^{\delta}$ is replaced by $M_0^{\delta} \cap C.$
The proof of the Lemma is then similar to the proof of Proposition \ref{hajeklem}. It suffices to replace inequality (\ref{eq:loglapM1})  by the inequality given in Lemma
\ref{eq:loglaplaceweak} and to set $\kappa = \sup_{x \in C \cap M_+ \setminus M_0^{\delta}} e^{\theta V(x)} e^{\theta \omega_0 T_1}$ with $\omega_0$ given by  Lemma \ref{lem:croissmaj}.
\qed

We now prove Proposition  \ref{hajeklemnoncompact}. Relying on Lemma \ref{hajeklemnoncompact0} fix $T_0, T_1$ such that
 $$T_1 \geq (2 \omega_0 + 1) T_0$$ where $\omega_0$ is given by Lemma \ref{lem:croissmaj}, and let $\theta, \rho, \kappa$ be given by Lemma \ref{hajeklemnoncompact0}.

  Let $T \in  [\frac{T_0 + T_1}{2}, T_1].$
 Using the strong Markov property, Lemma \ref{hajeklemnoncompact0} implies
 $$\E_x(e^{\theta V(X_T)} | {\cal F}_{\tau}) \leq \rho e^{\theta V(X_{\tau})} + \kappa \leq \rho e^{\theta (V(X_{\tau}) + \tau)} + \kappa  $$ on the event $\tau \leq T - T_0;$ and  Lemma  \ref{lem:croissmaj} $(i)$ implies that
$$E_x(e^{\theta V(X_T)} | {\cal F}_{\tau}) \leq e^{\theta V(X_{\tau})} e^{\omega_0 (T - \tau)}  \leq
 e^{\theta (V(X_{\tau}) + \tau)} e^{\theta \omega_0 T_0}  e^{-\theta (T -T_0)}$$ on the even $ T - T_0 < \tau \leq T.$
 Thus, by Lemma  \ref{lem:croissmaj} $(i)$
 $$E_x(e^{\theta V(X_T)} \Ind_{ \tau \leq T }) \leq \max (\rho, e^{\theta ((\omega_0 + 1) T_0 - T)} )e^{\theta V(x)} + \kappa.$$
 Also, by Lemma  \ref{lem:croissmaj} $(ii)$
  $$E_x(e^{\theta V(X_T)} \Ind_{ \tau > T }) \leq e^{- \theta T} e^{\theta V(x)}.$$

  Replacing $\rho$ by $\max \{\rho, e^{- \theta (T_0 + T_1)/2} , e^{- \theta (T_1 - T_0 (2 \omega_0 + 1))/2}\}$ and $T_0$ by $(T_0 + T_1)/2$ proves the result.

 \paragraph{Proof of  Theorem  \ref{th:expoconvnoncompact}} The proof of Theorem  \ref{th:expoconvnoncompact}  follows from Proposition \ref{hajeklemnoncompact}. The argument is verbatim the same as in the proof of Theorem \ref{th:expoconvcompact}.
\section{Appendix}
\subsection{Proof of Proposition \ref{extendeddomain}}
\label{sec:appendmartconv}
For any function $f \in \DA(\LA)$ and $x \in M$ recall that  the process $(M_t^f(x))$ (defined by
(\ref{martdefn}))
is a $({\cal F}_t)$ Martingale.  We let $\langle M^f (x)\rangle_t$ denote its {\em predictable quadratic variation}, defined as the compensator of $(M^f_t(x))^2.$

\blem \label{characGam}
Let $f \in \DA^2(\LA).$  Then
\beq
\langle M^f (x)\rangle_t=\int_0^t(\Gamma f)(X^x_s) ds
\eeq
\elem

\prf
The map $t \rar \int_0^t (\Gamma f)(X^x_s) ds$ is nondecreasing and continuous (hence predictable). It then remains  to show that
$((M_t^f(x))^2 - \int_0^t (\Gamma f)(X^x_s) ds)$ is a Martingale. This is a folklore result, for which we provide a proof.
 Let $M_t = f(x) + M_t^f(x)$ and  $N_t = f^2(x) + M_t^{f^2}(x).$
 Then $\{M_t\}_{t \geq 0}$ and $(N_t)_{t \geq 0}$ are both martingales. It then suffices to prove that $(Z_t)$ is a martingale, where
 $Z_t = M_t^2 - \int_0^t \Gamma f(X^x_s) ds - N_t.$
 Set $g_t = \LA f(X_t^x)$ and  $G_t = \int_0^t g_s ds.$ Then
 $$Z_t =  (f(X_t^x) - G_t)^2 - \int_0^t (\Gamma f)(X_s^x) ds - (f^2(X_t^x) - \int_0^t \LA (f^2) (X_s^x) ds)$$
 $$ = 2 \int_0^t f(X_s^x) g_s ds+ G_t^2 - 2 f(X_t) G_t = 2 \int_0^t (G_s + M_s) g_s ds + G_t^2 - 2(G_t + M_t) G_t.$$
 By Fubini formulae $G_t^2 = 2 \int_0^t G_s g_s ds.$ Thus
 $Z_t = 2 \int_0^t M_s g_s - G_t M_t$
 and
 $$Z_{t+u} - Z_t = 2 \int_t^{t+u} (M_s - M_{t+u}) g_s ds  + (M_t - M_{t+u}) G_t.$$ From this expression it is clear that $\E(Z_{t+u} - Z_t| \F_t) = 0$ for
 all $t, u \geq 0.$
\qed
\blem
\label{qlc0}
 Let $\tau$ be a stopping time. Then for all $x \in M$ $$\Ind_{M_0}(X_{\tau}^x) = \Ind_{M_0}(x) \mbox{ a.s on } \tau < \infty.$$
\elem
\prf
If $x \in M_0$ the event $E = \cap_{t \in \QQ^+} \{X^x_t \in M_0\}$ has probability one by Hypothesis \ref{hyp:standing}. By right continuity of paths and closeness of $M_0$ $E \subset  \cap_{t \in \RR^+} \{X^x_t \in M_0\}.$ In particular $X^x_{\tau} \in M_0$ a.s on $\{\tau < \infty\}.$

Suppose now $x \in M_+.$
Thus, using successively Hypothesis \ref{hyp:standing}, the strong Markov property (see Remark \ref{rem:strongmark}) and Hypothesis \ref{hyp:standing} we get
$$\E(\Ind_{M_0} (X^x_{\tau}) \Ind_{\{\tau \leq N\}}) =
\E(P_{N-\tau} \Ind_{M_0} (X^x_\tau) \Ind_{\{\tau \leq N\}})$$ $$
 = \E( \Ind_{M_0}(X^x_N) \Ind_{\{\tau \leq N\}}) \leq P_N  \Ind_{M_0}(x) = 0$$
 \qed

Recall that  $\{K_n\}_{n\geq 1}$ is the sequence of compact sets  as defined in (\ref{defKn}).
\blem \label{stoptimelem}
Let $x \in M_+$ and $\tau_n (x) \defn \inf\{t \geq 0: X_t^x \in K_n^c\}.$ Then $\{\tau_n (x)\}_{n \geq 1}$ is a localizing sequence. That is $\tau_n(x)$ is a stopping time and $\lim_{n \rar \infty} \tau_n(x) = \infty.$
\elem

\prf Fix $x \in M_+$ and set $\tau_n = \tau_n (x).$ Then $\tau_n$ is stopping time as $K_n^c$ is open and the filtration right continuous. Obviously, $\tau_n \leq \tau_{n+1}.$  Hence, $\lim_{n \rar \infty} \tau_n  = \tau \in \RR^+ \cup \{\infty\}$ exists a.s. and is a stopping time. Furthermore,
by Lemma \ref{qlc0}, $\tau_n < \tau$ a.s. on $\tau < \infty$ (since on $\{\tau_n = \tau; \, \tau < \infty\} \, X^x_{\tau_n} = X^x_{\tau} \in \cap_{m \geq n} K_m^c = M_0$).  The fact that $\tau_n < \tau$ implies that $(X_t)$ is almost surely left continuous at $\tau$ (i.e $X_{\tau^-} = X_{\tau}$)  on $\tau < \infty.$ This later property knows as a the quasi left continuity property is often proved for Feller processes but the proof only requires the cad-lag continuity of paths and the strong Markov property (see Remark \ref{rem:strongmark}).  Since $X^x_{\tau^-} \in M_0$ on $\tau <\infty$ we get that $X^x_{\tau} \in M_0$ and the conclusion follows from Lemma \ref{qlc0}.
\qed
We now prove  Proposition \ref{martconv}. Without loss of generality we assume that ${\cal M} = M_+,$ the proof for $M = {\cal M}$ being similar.
For all $x \in M_+, t \rar X_t^x \in M_+$ and has cad-lag paths. Thus $\{M^f_t(x)\}_{t \geq 0}$ is well-defined .
 Let  $\{\tau_n (x)\}_{n \geq 1}$ be as defined in Lemma \ref{stoptimelem}. Set $\tau_n = \tau_n(x).$ Then, by
assumption $(a)$ of the proposition and Lemma \ref{stoptimelem},
  \[M^{f}_{t \wedge \tau_n} (x) = f(X^x_{t \wedge \tau_n}) - f(x) - \int_0^{t \wedge \tau_n} g(X^x_s)ds = f_{K_n}(X^x_{t \wedge \tau_n}) - f(x) - \int_0^{t \wedge \tau_n} (\LA f_{K_n})(X^x_s)ds\]
is a martingale. Then,
$\{M^f_t (x)\}_{t \geq 0}$ is a local martingale. Now, by Lemma \ref{characGam}
\[\E(\langle M^{f} (x) \rangle_{t \wedge \tau_n}) = \E(\int_0^{t \wedge \tau_n} (\Gamma f_{K_n})(X^x_s)ds) \leq \int_0^t P_s \Gamma (f_{K_n})(x) ds  \leq C_x t\] for some constant $C_x.$
 Hence,
\beq
\label{bracbnd}
\E(\langle  M^f (x) \rangle_t) \leq C_x t < \infty
\eeq

This makes $\{M_t^f\} := \{M_t^f(x)\}$ a (true) $L^2$ martingale and $\{(M_t^f)^2 - \langle M^f \rangle_t\}$  a martingale. A proof can be found in  \cite{Legall}, theorem 4.3 for continuous martingales. The proof extends verbatim for right continuous martingales (provided we replace the quadratic variation by the predictable quadratic variation).

The last part of the proposition follows from the following standard argument.  For all integer $n$ and $\eps > 0,$  Doob's inequality for right continuous martingales implies that
 $$\Pr( \sup_{ 2^n \leq t \leq 2^{n+1}} \frac{|M_t^f|}{t} \geq \eps) \leq \Pr (\sup_{  t \leq 2^{n+1}} |M_t^f| \geq \eps 2^n) \leq \frac{1}{\eps^2 2^{2n}} \langle M^f \rangle_{2^{n+1}} \leq \frac{2 C_x}{\eps^2 2^n}.$$
 Thus, $\frac{M_t^f}{t} \rar 0$ a.s by Borel Cantelli.
 \qed
 \subsection{Proof of Theorem \ref{tightlypnv}}
\label{sec:appendtight}
The following Lemma is folklore and will be used repeatedly.
\blem
\label{lem:tightpi} Let $W$ be a nonnegative proper map, $C \geq 0$ and let $(\mu_n) \subset  {\cal P}(M)$ be such that $\limsup_{n \rar \infty} \mu_n W \leq C.$ Then,
\bdes
\iti The sequence $(\mu_n)$ is tight and every limit point $\mu$ of $(\mu_n)$ verifies $\mu W \leq C.$
\itii  Let $H : M \mapsto \RR$ be a continuous function  such that $\frac{W}{1 + |H|}$ is proper.  If $\mu_n \Rightarrow \mu$ then  $\mu_n H \rightarrow \mu H.$
\edes
\elem
\prf Assertion $(i)$ easily follows from Markov inequality and monotone convergence.

$(ii).$
Let $G = \frac{W}{1 + |H|}.$
 For all $R \in \RR \setminus D_G$ with $D_G$ at most countable, $\mu \{G = R\} = 0$ and, therefore,
 $$\lim_{n \rar \infty} \mu_n (H \Ind_{G \leq R}) =  \mu (H \Ind_{G \leq R}).$$
 On the other hand
 $\mu_n (|H| \Ind_{G > R}) \leq \mu_n (\frac{W}{G} \Ind_{G > R}) \leq \frac{1}{R} \mu_n(W).$
 Thus $$\lim_{R \rar \infty} \limsup_{n \rar \infty}   \mu_n (|H| \Ind_{G > R}) = 0$$ and, similarly, $$\lim_{R \rar \infty}    \mu (|H| \Ind_{G > R}) = 0.$$ This proves the result \qed

We now pass to the proof of Theorem \ref{tightlypnv}.

 $(i).$ Assumption $(i)$ of Hypothesis \ref{hyp:tightpi} makes the process
 $$M_t = W(X^x_t) - W(x) - \int_0^t LW(X^x_s) ds, t \geq 0$$  a  square integrable martingale satisfying the strong law of large numbers:
  $\lim_{t \rar \infty} \frac{M_t}{t} = 0$  a.s.
Thus, using condition $(ii)$ of Hypothesis \ref{hyp:tightpi},
$$0 \leq W(X_t^x) + \int_0^t \tilde{W}(X_s^x) ds \leq W(x) + Ct + M_t.$$ Taking the expectation and using Tonelli's Theorem proves assertion $(i).$

$(ii).$ Dividing by $t$ and letting $t \rar \infty$ proves that
 $\limsup_{t \rar \infty} \Pi_t^x \tilde{W} \leq C$ $\Pr$ a.s. Tightness follows from Lemma
\ref{lem:tightpi}.

It remains to show that limit points of $(\Pi_t^x)$ are invariant probabilities. For Feller discrete time Markov chains, this is a classical result (see e.g~ \cite{duf00}, Proposition 6.1.8). The proof easily adapts to the present setting as follows.

We claim that for each $f \in C_b(M)$ and $r > 0$ there exists a full measure set $\Omega_{f,r} \in {\cal F}$ such that for all $\omega \in \Omega_{f,r}$ $\lim_{t \rar \infty} \Pi_t^x(\omega) f - \Pi_t^x(\omega) P_r f = 0.$

Assume the claim is proved. Let ${\cal S} \subset C_0(M)$ be a countable dense subset of $C_0(M)$ (recall that $C_0(M)$ is separable) and $\Omega' = \bigcap_{f \in {\cal S}, r \geq 0, r \in \mathbb{Q}} \Omega_{f,r}$. Then, by density of ${\cal S},$  continuity of $r \mapsto  P_r f(x)$ (Hypothesis \ref{hyp:feller}) and dominated convergence,
$\mu(\omega) P_r f = \mu(\omega)  f$ for all $f \in C_0(M), r \geq 0, \omega \in \Omega'$ and $\mu(\omega)$ a limit point of $\{\Pi_t^x(\omega)\}_{t \geq 0}.$ This proves the result.

We now prove the claim. Replacing $(X_t)$ by  with $(X_{t r})$ we can  always assume that $r = 1.$
Set $$Qf(x) = \int_0^1 P_s f(x) ds, \; U_{k+1} = \int_{k}^{(k+1)} f(X_s) ds,$$
$$M_n = \sum_{k = 0}^{n-1} (U_{k+1} - Qf(X_k)), \; N_n = \sum_{k = 0}^{n-1} (Qf(X_{k+1}) - P_1 Qf(X_k)).$$   The sequences $(M_n)$ and $(N_n)$ are  martingales with bounded increments with respect to $\{\F_n\}.$ Thus, by the strong law of large number for martingales, $\lim_{n \rar \infty} \frac{1}{n} M_n = \lim_{n \rar \infty} \frac{1}{n} N_n = 0 $ $\Pr$ a.s.
Thus $$\lim_{n \rar \infty} \Pi_n^x f - \tilde{\Pi}^x_n Q f = \lim_{n \rar \infty} \tilde{\Pi}^x_n Q f  - \tilde{\Pi}^x_n P_1 Q f = 0$$  $\Pr$ a.s,  where $\tilde{\Pi}_n^x = \frac{1}{n} \sum_{k = 0}^n \delta_{X_k^x}.$  Replacing $f$ by $P_1 f$ also gives $$\lim_{n \rar \infty} \Pi_n^x P_1 f - \tilde{\Pi}^x_n Q P_1 f = 0$$ $\Pr$ a.s.
Since $P_1 Q f = Q P_1 f$ we then get that $$\lim_{n \rar \infty} \Pi_n^x f - \Pi^x_n P_1 f = 0.$$ $\Pr$ a.s. The claim is proved.

Probability $\mu$ is invariant if and only if $\mu P_t f = \mu f$ for all $t$ and $f \in C_b(M).$ Thus, by Feller continuity, ${\cal P}_{inv}(M)$ is closed and  compactness  equates  tightness. The latter  will follow from Lemma \ref{lem:tightpi} once we have proved that $\mu \tilde{W} \leq C$ for all $\mu \in {\cal P}_{inv}(M).$ Let $\mu \in {\cal P}_{inv}(M).$ First assume $\mu$ ergodic. Then, by Birkhoff ergodic Theorem, $\Pi_t^x \Rightarrow \mu$ for $\mu$ almost every $x$ and $\mathbb{P}_x$ almost surely. Thus, $\mu \tilde{W} \leq C$ by Lemma \ref{lem:tightpi} $(i).$ If now $\mu$ is invariant, the ergodic decomposition theorem, implies that $\mu \tilde{W} \leq C.$ This concludes the proof of assertion $(ii).$

$(iii).$ Set $w(t) = P_t W(x).$ Using the semigroup property and Fubini-Tonelli, we get that
\begin{eqnarray}
\label{ineqdiff}
\label{ineq1}
w(t+s) - w(t) & \leq & - \alpha \int_t^{t+s} w(r) dr + C s \leq   C s\\ \label{ineq2}
w(t) - w(t-u) & \leq &- \alpha \int_{t-u}^{t} w(r) dr + C s \leq  C u
                                                          \end{eqnarray}
                                                          for all $t \geq 0, s \geq 0$ and $0 \leq u   \leq t.$
On the other hand,
by Fatou Lemma and right continuity of $t \rar W(X_t^x)$ $$\liminf_{s \rar 0, s > 0} w(t+s) = \liminf_{s \rar 0, s > 0} \E( W(X_{t+s}^x)) \geq \E(W(X_t^x)) = w(t).$$ Combined with (\ref{ineq1})
this shows that $t \rar w(t)$ is right-continuous. From (\ref{ineq2}) we also get that $t \rar w(t)$ is lower semi continuous.
Set $\Delta^+ w (t) = \limsup_{s \rar 0, s> 0} \frac{w(t+s) - w(t)}{s}$ and  $\Delta^- w (t) = \limsup_{s \rar 0, s> 0} \frac{w(t) - w(t-s)}{s}.$
 Using (\ref{ineq1}) and right continuity, we get that
$$\Delta^+ w(t) \leq - \alpha w(t) + C.$$ Using (\ref{ineq2}) and lower semi continuity we get that
$$\Delta^- w(t) \leq -\alpha w(t) + C.$$
Set now $\tilde{w}(t) = e^{\alpha t} (w(t) - \frac{C}{\alpha}) - \epsilon t$ for some $\eps > 0.$ Then, defining $\Delta^{+,-} \tilde{w}$ like $\Delta^{+,-} w$  with $\tilde{w}$ in place of $w$ we get that
$$\Delta^+ \tilde{w}(t) \leq -\eps \mbox{ and } \Delta^- \tilde{w}(t) \leq -\eps.$$
This implies that for all $t \geq 0$ there exists an open subset of $\RR^+,$  $I_t$ containing $t$ such that $w(s) \leq w(t)$ for all $s \in I_t.$
 In particular the set $\{t \geq 0 \: : \tilde{w}(t) \leq \tilde{w}(0)\}$ is open in $\RR^+$. By lower semi continuity of $\tilde{w},$
  it also closed. Being nonempty it equals $\RR^+$ by connectedness.
Thus $\tilde{w}(t) \leq \tilde{w}(0)$ for all $t.$ Since $\eps$ is arbitrary this leads to
$$P_t W(x) = w(t) \leq e^{-\alpha t} (w(0) - \frac{C}{\alpha}) + \frac{C}{\alpha}.$$

\qed

\subsection{Proof of Proposition \ref{prop:ecosde}}
\label{sec:append}
$(i).$ By local Lipschitz continuity  and classical results on stochastic differential equations,
 there exists for any $x \in \RR^n$  a unique continuous process $(X_t^x)$
 defined on some interval $[0, \tau^x[$ solution to (\ref{eq:sde}), with initial condition $X_0^x = x$ and such that
  $t < \tau^x \Leftrightarrow \|X_t^x\| < \infty$ (see e.g~\cite{RY} Chapter IX, exercise 2.10).
  Furthermore, it is easily checked (by Ito formula and uniqueness of the solutions) that, if $\alpha_i \neq 0$
  $$X_{t,i} =
  x_i \exp{ \left ( \int_0^t [X_{s,i}^{\alpha_i - 1} F_i(X_s)- \frac{1}{2} X_{s,i}^{2(\alpha_i -1)} a_{ii}(X_s)] ds + \sum_j \int_0^t  X_{s,i}^{\alpha_i -1}
    \Sigma_i^j(X_s) dB_s^j \right )}$$
    where, to shorten notation, $X_t$ stands for $X_t^x.$
    Thus
    \beq
    \label{invsde1}
    x_i > 0 \Rightarrow X_{t,i}^x > 0  \mbox{ for all } t \in [0,\tau^x[
    \eeq and
    \beq
    \label{invsde2}
    x_i = 0 \Rightarrow X_{t,i}^x = 0  \mbox{ for all } t \in [0,\tau^x[ \eeq
  We shall now prove that $\tau^x = \infty.$

For any $C^2$ function $\psi : M \mapsto \RR,$ by Ito formulae,
\beq
\label{eq:ito}
\psi(X_t^x) - \psi(x) - \int_0^t L\psi(X_s^x) ds = \sum_i \int_0^t
\frac{\partial \psi}{\partial x_i}(X_s^x) \left[(X_{s,i}^x)^{\alpha_i}\sum_{j = 1}^m  \Sigma_i^j(X_s^x)\right]dB^j_s,
\eeq
 Let $\tau_k^x = \inf\{ t \geq  0 :\: U(X_t^x) \geq k \}$ for all $k \in \NN.$
  By the assumption on $U,$   for all $x \in M,$
 $$LU(x)  \leq - \alpha U(x) + \beta.$$
Thus
\begin{eqnarray}
\label{eq:nonexplos}
 k \Pr(\tau_k^x \leq t) & = & \E(U(X^x_{\tau_k^x}) \Ind_{\tau_k^x \leq t}) \\ \label{eq1PW}
& \leq & \E(U(X^x_{t \wedge \tau_k^x} ))  =  U(x) +
\E( \int_0^{t\wedge\tau_k^x} L  U (X^x_s) ds) \\ \label{eq2PW}
   & \leq & U(x) - \alpha \E (\int_0^{t\wedge\tau_k^x} U(X_s^x) ds ) + \beta t \\  \label{eq3PW}
    & \leq &  U(x) + \beta t
\end{eqnarray}
Hence
$$\Pr(\tau^x \leq t) = \Pr( \cap_{k \geq 0} \{\tau_k^x \leq t) ) = \lim_{k \rar \infty} \Pr(\tau_k^x \leq t) = 0$$
proving that $\tau^x = \infty$ almost surely.

We now let $(P_t)$ denote the semigroup acting on bounded (respectively non-negative) measurable functions
$f : M \mapsto \RR,$ by $P_t f(x) = \E(f(X_t^x)).$ $C_b(M)-$ Feller continuity just follows from Lebesgue dominated convergence theorem and the continuity in $x$
of the solution $X_t^x.$

$(ii).$
Inequalities (\ref{eq2PW}, \ref{eq3PW}) and monotone convergence imply
that $$P_t U(x) = \E(U(X_t^x)) = \lim_{k \rar \infty} \E(U(X_t^x) \Ind_{\tau_k^x \geq t}) \leq
U(x) - \alpha \E(\int_0^t U(X_s^x) ds) + \beta t $$
$$ = U(x) - \alpha \int_0^t P_s U(x) ds + \beta t  \leq U(x) + \beta t$$  where the last equality follows from Fubini-Tonelli theorem.
Thus, reasoning exactly like in the proof of Theorem \ref{tightlypnv} $(iii)$ we get that
$$P_t U(x) \leq e^{-\alpha t} (U(x) - \beta/\alpha) + \beta/\alpha.$$
$(iii).$ Let $\psi \in C^2_c(M).$ By Ito formulae
$\psi(X_t^x) - \psi(x) - \int_0^t L\psi(X_s^x) ds$ is a Martingale. Thus, taking the expectation, $P_t \psi(x) - \psi(x) = \int_0^t P_s (L\psi)(x) ds.$
Thus $|P_t(\psi)(x) - \psi(x)| \leq t \|L \psi\|$ and
$$\lim_{t \rar 0} \frac{P_t \psi(x) - \psi(x)}{t} = L\psi(x).$$ This proves that $\psi \in \DA(\LA)$ and $L\psi = \LA \psi.$  Replacing $\psi$ by $\psi^2$ shows that $\psi \in \DA^2(\LA)$ and $\Gamma(\psi) = \Gamma_L(\psi).$

$(iv)$ is immediate from (\ref{invsde1}) and (\ref{invsde2}).

$(v).$
 For any smooth function $h : \RR^+ \mapsto \RR$
 $$L (h (U)) = h'(U) L U + \frac{1}{2}h''(U) \Gamma_L (U).$$ If $h$ is concave and nondecreasing, this gives
  $$L (h (U)) \leq h'(U) L U \leq - \alpha h'(U)  U (1 + \varphi) + \beta h'(1).$$
Set  $h(t) = t^{\frac{1-\eta}{2}}$  and $W = h(U).$ Then $h'(t) t  = \frac{1-\eta}{2} h(t).$ Thus
$$L (W ) \leq  \frac{1-\eta}{2}(- \alpha W (1 + \varphi) + \beta).$$
Now
$$\Gamma_L(W) = h'(U)^2 \Gamma_L (U) = (\frac{1-\eta}{2})^2 U^{-\eta - 1} \Gamma_L(U).$$   Thus
 $$\Gamma_L(W) \leq cst(1 + U).$$
Let $B : \RR \mapsto \RR$ be a smooth function such that $B(t) = t$ for $t \leq 1,$ $B(t) = 2$
for $t \geq 3,$ and
$0 \leq B'(t) \leq 1.$ Set $W_n = n B (W/n)$  Then $W_n \in \DA^2(\LA)$ (since $W_n - 2n \in C_c^2(M)$), $W_n(x) = W(x)$ and  $\LA W_n(x) = L W(x)$ whenever $W(x) \leq n.$
On the other hand $\Gamma(W_n)(x) =   B'^2(W/n) \Gamma_L(W)(x) \leq \Gamma_L(W)(x).$ Thus
 $$\sup_{t \geq 0,n} P_t (\Gamma(W_n))(x) \leq \sup_{t \geq 0} P_t \Gamma_L(W) (x) \leq cst(1 +  \sup_{t \geq 0} P_t U(x)) < \infty$$
 where the last inequality follows from $(ii).$ Hypothesis \ref{hyp:tightpi}  then follows from Proposition \ref{martconv}.

\bibliographystyle{amsplain}
\bibliography{persistence}

\section*{Acknowledgments}
This work is supported by the SNF grant $200021_15728$. I thank Jean Baptiste Bardet, Patrick Cattiaux, Alex Hening, Tobias Hurt, Eva Locherbach, Florent Malrieu, Janusz Mierczynski, Sebastian Schreiber,  Edouard Strickler, Pierre André Zitt for valuable discussions on different topics related to this paper. Special thanks to Edouard Strickler for his help with the simulation of the May Leonard process with Scilab and to Jean Baptiste Bardet for his help with the computations of Lie Brackets using the  formal software Python/Sympy.

\end{document}